\documentclass[12pts]{amsart}

\usepackage{amscd,latexsym,amsthm,amsfonts,amssymb,amsmath,amsxtra}
\usepackage[mathscr]{eucal}
\newcommand{\BA}{{\mathbb {A}}}

\newcommand{\BC}{{\mathbb {C}}}

\newcommand{\RO}{{\mathrm {O}}}

\newcommand{\diag}{{\mathrm{diag}}}

\newcommand{\GL}{{\mathrm{GL}}}

\newcommand{\Ind}{{\mathrm{Ind}}}

\newcommand{\I}{{\mathrm{I}}}

\newcommand{\Res}{{\mathrm{Res}}}

\newcommand{\Sp}{{\mathrm{Sp}}}

\newtheorem{thm}{Theorem}[section]
\newtheorem{cor}[thm]{Corollary}
\newtheorem{lem}[thm]{Lemma}
\newtheorem{prop}[thm]{Proposition}
\newtheorem {conj}[thm]{Conjecture}

\begin{document}
	\renewcommand{\theequation}{\arabic{equation}}
	\numberwithin{equation}{section}

	\title{A new regularized Siegel-Weil type formula\\ part I}
	
	\author{David Ginzburg}
	
	\address{School of Mathematical Sciences, Sackler Faculty of Exact Sciences, Tel-Aviv University, Israel
		69978} \email{ginzburg@tauex.tau.ac.il}

	
	\author{David Soudry}
	\address{School of Mathematical Sciences, Sackler Faculty of Exact Sciences, Tel-Aviv University, Israel
		69978} \email{soudry@tauex.tau.ac.il}
	
	
	
	
	
	\maketitle

	\section{The proposed new type of the Siegel-Weil formula: main theorems and conjectures}
	
	In this paper, we propose a formula relating certain residues of Eisenstein series on  symplectic groups. These Eisenstein series are attached to parabolic data coming from Speh representations. The proposed formula bears a strong similarity to the regularized Siegel-Weil formula, established by Kudla and Rallis, \cite{KR94}, for symplectic-orthogonal dual pairs. Their work was later generalized by Ikeda, Moeglin, Ichino, Yamana, Gan-Qiu-Takeda and others. See \cite{GQT14} and the references therein. 
	
	\subsection{The work of Kudla-Rallis on the Siegel-Weil formula.}
	We start by reviewing the work of Kudla-Rallis, mainly from \cite{KR94}, focusing and connecting, as in the introduction of \cite{GQT14}, the theta correspondence, Rallis inner product formula, the Siegel-Weil formula, the doubling method and $L$-functions. Let $F$ be a number field and $\BA$ its ring of adeles. Consider a dual pair $(\Sp_{2n},\RO_{2m})$ inside $\Sp_{4mn}$, where $\Sp_{2n}$ denotes the symplectic group of rank $n$, regarded as an algebraic group over $F$, and $\RO_{2m}$ is an orthogonal group corresponding to a quadratic space over $F$, $(V,Q)$, where $V$ is a $2m$ dimensional vector space over $F$, and $Q$ is an $F$- nondegenerate, symmetric bilinear form on $V\times V$, with Witt index $r$. We assume that $\RO_{2m}$ is not binary and split. Denote by $\chi_V$ the corresponding quadratic character of $F^*\backslash \BA^*$. Fix a nontrivial character $\psi$ of $F\backslash \BA$. Let $\pi$ be an irreducible, cuspidal, automorphic representation of $\Sp_{2n}(\BA)$, and consider $\theta_{\psi,2m}(\pi)$ - the theta lift (with respect to $\psi$) of $\pi$ to $\RO_{2m}(\BA)$. Its space is spanned by the functions on $\RO_{2m}(\BA)$,
	\begin{equation}\label{0.1}
		\theta_{\psi,2m}^\phi(\varphi_\pi)(h)=\int\limits_{\Sp_{2n}(F)\backslash \Sp_{2n}(\BA)}\theta_{\psi,4mn}^\phi(g,h)\varphi_\pi(g)dg,
	\end{equation}
	where $\varphi_\pi$ is a cusp form in the space of $\pi$, $\phi\in \mathcal{S}(V(\BA)^n)$, and $\theta_{\psi,4mn}^\phi(g,h)$ is the restriction to the image of the dual pair $\Sp_{2n}(\BA)\times \RO_{2m}(\BA)$, inside the double cover $\widetilde{\Sp}_{4mn}(\BA)$, of the theta series
	$$
	\theta_{\psi,4mn}^\phi(g,h)=\sum\limits_{x\in V(F)^n}\omega_{\psi,4mn}(g,h)\phi(x).
	$$
	Here, $\omega_{\psi,4mn}$ denotes the Weil representation of $\widetilde{\Sp}_{4mn}(\BA)$, corresponding to $\psi$, acting in $\mathcal{S}(V(\BA)^n)$. To study the question of non-vanishing of $\theta_{\psi,2m}(\pi)$, consider, formally, the $L^2$-inner product of two theta lifts of the form \eqref{0.1}. Applying one more formal manipulation, we get\\
	\\
	$	(\theta_{\psi,2m}^{\phi_1}(\varphi_\pi), \theta_{\psi,2m}^{\phi_2}(\varphi'_\pi))=$
	\begin{equation}\label{0.2}
		=\int\limits_{[\Sp_{2n}\times \Sp_{2n}]}\varphi_\pi(g_1)\bar{\varphi}'_\pi(g_2)(\int\limits_{[\RO_{2m}] }\theta_{\psi,4mn}^{\phi_1}(g_1,h)\theta_{\psi^{-1},4mn}^{\bar{\phi_2}}(g_2,h)dh)dg_1dg_2.
	\end{equation}
	We used the shorthand notation $[G]=G(F)\backslash G(\BA)$. The goal is to make sense out of the integral \eqref{0.2}. The first thing is to note that the $dh$-integrand in \eqref{0.2} can be expressed as one theta series. This is the multiplicative property of theta series, which is easy to establish, namely, for $g_1,g_2\in \Sp_{2n}(\BA)$, $h\in \RO_{2m}(\BA)$,
	\begin{equation}\label{0.3}
		\theta_{\psi,4mn}^{\phi_1}(g_1,h)\theta_{\psi^{-1},4mn}^{\bar{\phi_2}}(g_2,h)=\theta_{\psi,8mn}^{\phi_1\otimes\bar{\phi}_2}((g_1,g_2),h).
	\end{equation}
	The r.h.s. of \eqref{0.3} is a theta series on $\widetilde{\Sp}_{8mn}(\BA)$, restricted, first, to the image of $\Sp_{4n}(\BA)\times \RO_{2m}(\BA)$, and then to $(\Sp_{2n}(\BA)\times \Sp_{2n}(\BA))\times \RO_{2m}(\BA)$. As we shall see, the analog, in our case, of \eqref{0.3} turns out to be quite involved. The next step in figuring out the meaning of \eqref{0.2} is to interpret the $dh$-inner integral, using \eqref{0.3}, that is, for $g\in \Sp_{4n}(\BA)$, $\Phi\in \mathcal{S}(V(\BA)^{2n})$ ($\Phi=\phi_1\otimes \bar{\phi}_2$ in \eqref{0.2}),
	\begin{equation}\label{0.4}
		I(\Phi,g)=\int\limits_{\RO_{2m}(F)\backslash \RO_{2m}(\BA)}\theta_{\psi,8mn}^\Phi(g,h)dh.
	\end{equation}
	The integral \eqref{0.4} is absolutely convergent when $r=0$, or when $2m-r>2n+1$. In this range, we have the Siegel-Weil formula, proved by Weil and Kudla-Rallis,
	\begin{equation}\label{0.5}
		I(\Phi,g)=\kappa E(f_{\Phi,s},g)\Big|_{s=m-n-\frac{1}{2}},
	\end{equation}
	where $\kappa=1,2$ and $E(f_{\Phi,s})$ is the Eisenstein series on $\Sp_{4n}(\BA)$ attached to the Siegel-Weil section $f_{\Phi,s}$ of the parabolic induction $\Ind_{Q_{2n}(\BA)}^{\Sp_{4n}(\BA)}\chi_V\circ \det |\det\cdot|^s$,
	\begin{equation}\label{0.5'}
	f_{\Phi,s}(g)=\omega_{\psi,8mn}(g,I_{2m})\Phi(0)a(g)^{s-m+n+\frac{1}{2}}.
	\end{equation}
	Here, $Q_{2n}$, with Levi decomposition $Q_{2n}=L_{2n}\ltimes U_{2n}$ is the Siegel parabolic subgroup, and $a(g)$ is obtained by writing the Iwasawa decomposition $g=\hat{m}_gu_gk_g$, where $k_g\in K_{\Sp_{4n}}$, the maximal compact subgroup of $\Sp_{4n}(\BA)$, $u_g\in U_{2n}(\BA)$, $\hat{m}_g=diag(m_g,m_g^*)\in L_{2n}(\BA)$, where $m_g\in \GL_{2n}(\BA)$. Then $a(g)=|\det(m_g)|$. Note that $\omega_{\psi,8mn}(x)\Phi(0)$ is the constant term at $x$ of $\theta_{\psi,8mn}^\Phi$ along the Siegel radical of $\Sp_{8mn}$,
	\begin{equation}\label{0.5''}
		\int\limits_{U_{4mn}(F)\backslash U_{4mn}(\BA)}\theta_{\psi,8mn}^\Phi(ux)du= \omega_{\psi,8mn}(x)\Phi(0),
	\end{equation}	
	where $U_{4mn}$ is the unipotent radical of the Siegel parabolic subgroup $Q_{4mn}$ of $\Sp_{8mn}$. Indeed, the last integral is equal
	$$
	\int\limits_{S_{4mn}(F)\backslash S_{4mn}(\BA)}\sum_{\xi\in F^{4mn}}\psi(tr(\xi\cdot z\cdot{}^t\xi))\omega_{\psi,8mn}(x)\Phi(\xi)dz=\omega_{\psi,8mn}(x)\Phi(0).
	$$
	Here, $S_{4mn}$ denotes the space of $4mn\times 4mn$ symmetric matrices.	
	
	When we substitute \eqref{0.5}, with $\Phi=\phi_1\otimes\bar{\phi}_2$, in \eqref{0.2}, we get\\
	\\
	$	(\theta_{\psi,2m}^{\phi_1}(\varphi_\pi), \theta_{\psi,2m}^{\phi_2}(\varphi'_\pi))=$
	\begin{equation}\label{0.6}
		=\kappa \int\limits_{[\Sp_{2n}\times \Sp_{2n}]}\varphi_\pi(g_1)\bar{\varphi}'_\pi(g_2)E(f_{\Phi,s},(g_1,g_2))\Big|_{s=m-n-\frac{1}{2}}dg_1dg_2.
	\end{equation}
	This is the global integral of the doubling method of Piatetski-Shapiro and Rallis \cite{P-SR87}, at $s=m-n-\frac{1}{2}$. This global integral represents $L(\pi\times \chi_V,s+\frac{1}{2})$, up to normalization. In our case, the analogous integrals of the generalized doubling method of \cite{CFGK17} show up.\\
	\\ 
	{\bf Regularization:} In the range $2m-r\leq 2n+1$, $r\geq 1$, Kudla and Rallis found elements $z\in \mathfrak{z}_{\scriptstyle{sp_{4n}(F_v)}}$, $z'\in \mathfrak{z}_{\scriptstyle{o_{2m}(F_v)}}$, in the centers of the enveloping algebras of the Lie algebras of $\Sp_{4n}(F_v)$, $\RO_{2m}(F_v)$, at one archimedean place $v$ satisfying
	\begin{thm}\label{thm 0.1} 
		For $\Phi\in \mathcal{S}(V(\BA)^{2n})$, 
		$$
		\omega_{\psi,8mn}(z)\Phi=\omega_{\psi,8mn}(z')\Phi,
		$$
		and $\theta_{\psi,8mn}^{\omega_{\psi,8mn}(z)\Phi}(g,h)$ is rapidly decreasing in $h\in \RO_{2m}(F)\backslash \RO_{2m}(\BA)$.
	\end{thm}
	Ichino obtained a similar result by convolving $\Phi$ against a function in the spherical Hecke algebra of $\Sp_{4n}(F_{v'})$, at one nonarchimedean place $v'$. See \cite{I01}, Sec. 1. We will obtain an analogous theorem in our case. Next, Kudla and Rallis took an Eisenstein series $E(h,\zeta)$ on $\RO_{2m}(\BA)$, attached to the maximal parabolic subgroup with Levi part isomorphic to $\GL_r\times \RO_{2(m-r)}$ and the character of its adele points given by $|\det_{\GL_r}\cdot|^\zeta$. This Eisenstein series has a simple pole at $\zeta=m-\frac{r+1}{2}$, with constant residue. Then they introduce 
	\begin{equation}\label{0.7}
		\mathcal{E}(g,\Phi,\zeta)=\frac{1}{P(\zeta)}\int\limits_{\RO_{2m}(F)\backslash \RO_{2m}(\BA)}\theta_{\psi,8mn}^{\omega_{\psi,8mn}(z)\Phi}(g,h)E(h,\zeta)dh,
	\end{equation}
	where $P(\zeta)$ is the polynomial obtained by the action of $z'$ (from Theorem \ref{thm 0.1}) on $E(h,\zeta)$. Note that, from Theorem \ref{thm 0.1}, the integral \eqref{0.7} converges absolutely, away from the poles of $E(h,\zeta)$. They prove
	\begin{thm}\label{thm 0.2}
		$\mathcal{E}(g,\Phi,\zeta)$ is an Eisenstein series on $\Sp_{4n}(\BA)$ attached to the maximal parabolic subgroup with Levi part isomorphic to $\GL_r\times \Sp_{2(2n-r)}$ and the representation of its adele points given by $|\det_{\GL_r}\cdot|^\zeta \otimes \theta_{\psi, 2(2n-r)}(1_{\RO(V_{an})})$, where $\theta_{\psi, 2(2n-r)}(1_{\RO(V_{an})})$ denotes the theta lift to $\Sp_{2(2n-r)}(\BA)$ of the trivial representation of the adele points of the orthogonal group of the anisotropic kernel $V_{an}$ of $V$.
	\end{thm}
	We will prove an analogous theorem in our case. Kudla and Rallis compute $P(\zeta)$ explicitly and find out that when $m\leq n$ (and then $2m-r\leq 2n+1$), $P(m-\frac{2r+1}{2})\neq 0$, so that $\mathcal{E}(g,\Phi,\zeta)$ has at most a simple pole at $\zeta=m-\frac{r+1}{2}$, and then
	\begin{equation}\label{0.8}
		Res_{\zeta=m-\frac{r+1}{2}}	\mathcal{E}(g,\Phi,\zeta)=\frac{c}{P(m-\frac{r+1}{2})}\int\limits_{\RO_{2m}(F)\backslash \RO_{2m}(\BA)}\theta_{\psi,8mn}^{\omega_{\psi,8mn}(z)\Phi}(g,h)dh,
	\end{equation}	
	where $c=Res_{\zeta=m-\frac{r+1}{2}}E(h,\zeta)$. This is the interpretation, or the regularization of the integral \eqref{0.4}. Denote $B_{-1}(g,\Phi)=Res_{\zeta=m-\frac{r+1}{2}}	\mathcal{E}(g,\Phi,\zeta)$. When $m>n$ and $2m-r\leq 2n+1$, $P(m-\frac{r+1}{2})=0$, so that $\mathcal{E}(g,\Phi,\zeta)$ has at most a double pole at $\zeta=m-\frac{r+1}{2}$. Denote by $B_{-2}(g,\Phi)$ the leading term of the Laurent expansion of $\mathcal{E}(g,\Phi,\zeta)$ around $\zeta=m-\frac{r+1}{2}$.
	The generalization, of the Siegel-Weil formula \eqref{0.5} in the convergence range, is
	\begin{thm}\label{thm 0.3}(The regularized Siegel-Weil formula: first term identity)\\
		1. Assume that $m\leq n$. Then
		$$
		2Res_{\zeta=m-\frac{r+1}{2}}\mathcal{E}(g,\Phi,\zeta)=Value_{s=m-n+\frac{1}{2}}E(f_{\Phi,s},g).
		$$
		If $m<n$, then we also have
		$$
		2Res_{\zeta=m-\frac{r+1}{2}}\mathcal{E}(g,\Phi,\zeta)=Res_{s=n-m-\frac{1}{2}}E(f_{\Phi',s},g),
		$$
		for an appropriate $\Phi'\in \mathcal{S}(V'(\BA)^{2n})$, where $V'$ is the complementary space of $V$ (see \cite{KR94}, p. 4).\\
		2. Assume that $2n+2\leq 2m\leq 2n+r+1$. Then, with similar notation,
		$$
		B_{-2}(g,\Phi)=B_{-1}(g,\Phi')=Res_{s=m-n-\frac{1}{2}}E(f_{\Phi,s},g).
		$$
		
	\end{thm}
	Now, as in \eqref{0.6}, assume that $\chi$ is a quadratic character of $F^*\backslash \BA^*$ is such that the partial $L$-function $L^S(\pi\times\chi,s)$ has a pole at $s=k$, a positive integer. Kudla and Rallis show that these are the only possible poles, and that necessarily $k\leq [\frac{n}{2}]+1$. Let $\tilde{m}=n+k$. Using the doubling integral, as in \eqref{0.6}, to represent this $L$-function, they show that there is a quadratic space $V'$ over $F$, of dimension $2\tilde{m}$, with $\chi=\chi_{V'}$, and $\Phi\in \mathcal{S}(V'(\BA)^{2n})$, such that the following integral is not (identically) zero,
	$$
	\int\limits_{[\Sp_{2n}\times \Sp_{2n}]}\varphi_\pi(g_1)\bar{\varphi}'_\pi(g_2)\Res_{s=m-n-\frac{1}{2}}(E(f_{\Phi,s},(g_1,g_2)))dg_1dg_,
	$$
	and hence, by Theorem \ref{thm 0.3}(2), so is the integral
	$$
	\int\limits_{[\Sp_{2n}\times \Sp_{2n}]}\varphi_\pi(g_1)\bar{\varphi}'_\pi(g_2)B_{-1}(\Phi',(g_1,g_2)))dg_1dg_2,
	$$	
	where $\Phi'\in \mathcal{S}(V^{2n}(\BA))$, and $V$ is the quadratic space, complementary to $V'$, $dim(V)=4n+2-2\tilde{m}=2n+2-2k:=2m$. Kudla and Rallis show that there are $\phi_1, \phi_2\in \mathcal{S}(V(\BA)^n)$, such that we may replace in the last integral remains $\Phi'$ by $\phi_1\otimes \bar{\phi}_2$. Using \eqref{0.7}, \eqref{0.3}, we conclude that the theta lift of $\pi$ to $\RO_{2m}(\BA)$ (corresponding to $V$) is nontrivial. This was one of the main goals of Kudla and Rallis, that is
	\begin{thm}\label{thm 0.4}
		Let $\pi$ be an irreducible, cuspidal, automorphic representation of $\Sp_{2n}(\BA)$, and $\chi$, a quadratic character of $F^*\backslash \BA^*$. Assume that $L^S(\pi\times\chi,s)$ has a pole at $s=k$, a positive integer. Then $k\leq [\frac{n}{2}]+1$, and there is a quadratic space $V$, over $F$, of dimension $2m$, $m=n+1-k$, and $\chi_V=\chi$, such that $\theta_{\psi,2m}(\pi)$, to $\RO_{2m}(\BA)$, corresponding to $V$, is nontrivial.
	\end{thm}
	
	Our program is to follow a similar itinerary, guided by the poles of the $L$-functions for $\Sp_{2n}(\BA)\times \GL_d(\BA)$, $L(\pi\times\tau,s)$, where $\pi, \tau$ are irreducible, cuspidal, automorphic representations of $\Sp_{2n}(\BA)$, $\GL_d(\BA)$, respectively. We now know the generalized doubling integrals of Cai, Friedberg, Ginzburg and Kaplan \cite{CFGK17}. Then we ask whether there is a analogous new ``theta correspondence'' characterizing the poles of $L(\pi\times\tau,s)$, and then is there a related new Siegl-Weil formula?
	
	\subsection{Notation:} Before we continue, we set up some notation. We will write the symplectic group $\Sp_{2k}$ as the subgroup of $\GL_{2k}$ of matrices $g$ satisfying
	$$
	{}^tgJ_{2k}g=J_{2k},
	$$
	where $J_{2k}=\begin{pmatrix}&w_k\\-w_k\end{pmatrix}$, and $w_k$ is the $k\times k$ permutation matrix with $1$ along the anti-diagonal. Let $1\leq r\leq k$ be an integer. We denote by $Q_r$ the standard parabolic subgroup of $\Sp_{2k}$, with Levi decomposition $Q_r=L_r\ltimes U_r$, where $L_r\cong \GL_r\times \Sp_{2(k-r)}$. The Siegel parabolic subgroup of $\Sp_{2k}$ is $Q_k$. The elements of its Levi part $L_k$ are
	$$
	\hat{a}=\begin{pmatrix}a\\&a^*\end{pmatrix},\ a^*=w_k{}^t a^{-1}w_k,\ a\in \GL_k;
	$$
	The elements of $U_k$ are
	$$
	u_k(x)=\begin{pmatrix} I_k&x\\&I_k\end{pmatrix},\ {}^t(w_kx)=w_kx.
	$$
	For $a\in \GL_r, r\leq k$, we will also denote
	$$
	\hat{a}=\begin{pmatrix}a\\&I_{2(k-r)}\\&&a^*\end{pmatrix}\in \Sp_{2k},
	$$
	when $k$ is understood.
	
	More generally, for positive integers $\underline{i}=(i_1,...,i_\ell)$, such that $i=i_1+\cdots i_\ell\leq k$, we denote by $Q_{\underline{i}}$ the standard parabolic subgroup of $\Sp_{2k}$ with Levi part $L_{\underline{i}}\cong \GL_{i_1}\times\cdots\times\GL_{i_\ell}\times\Sp_{2(k-i)}$. We denote its unipotent radical by $U_{\underline{i}}$. When we want to recall that these are subgroups of $\Sp_{2k}$, we denote $Q^{2k}_{\underline{i}}, L^{2k}_{\underline{i}}, U^{2k}_{\underline{i}}$. For positive integers $\underline{j}=(j_1,...,j_\ell)$, such that $j_1+\cdots j_\ell=k$, we denote by $P_{\underline{j}}=M_{\underline{j}}\ltimes V_{\underline{j}}$ the standard parabolic subgroup of $\GL_k$, with unipotent radical $V_{\underline{j}}$ and Levi part $=M_{\underline{j}}\cong \GL_{j_1}\times\cdots\times\GL_{j_\ell}$. We will denote the standard Borel subgroups of $\Sp_{2k}$, $\GL_k$ by $B_{\Sp_{2k}}$, $B_{\GL_k}$. We will denote the corresponding diagonal subgroups by $T_{\Sp_{2k}}, T_{\GL_k}$. We will sometimes denote $T_{\GL_k}=T_k$
	
	Let $v$ be a place of $F$. We denote by $K_{2m,v}$ the standard maximal compact subgroup of $\Sp_{2m}(F_v)$. Similarly, we denote by $K_{\GL_m,v}$ the standard maximal compact subgroup of $\GL_m(F_v)$. When $v$ is finite, we denote by $\mathcal{O}_v$ the ring of integers of $F_v$, and by $\mathcal{P}_v$ its maximal ideal. Denote by $q_v$ the number of elements in the residue field $\mathcal{O}_v/\mathcal{P}_v$, and by $p_v$ a agenerator of $\mathcal{P}_v$. Then $K_{2m,v}=\Sp_{2m}(\mathcal{O}_v)$, $K_{\GL_m,v}=\GL_m(\mathcal{O}_v)$. We denote $K_{2m}=\prod_vK_{2m,v}$, $K_{\GL_m}=\prod_vK_{\GL_m,v}$.

	\subsection{A new theta correspondence, conjectural new Siegel-Weil formulas and applications.} We start with the correspondence constructed by Ginzburg in \cite{G03}. It is different from the classical theta correspondence. For a given irreducible, self-dual, cuspidal, automorphic representation $\tau$ of $\GL_d(\BA)$, there is a space of theta kernel functions on the adele points of a commuting pair of symplectic groups inside a larger symplectic group. This pair is not a reductive dual pair. We restrict ourselves to irreducible, cuspidal, automorphic representations $\tau$ of $\GL_2(\BA)$, with trivial central character, and such that $L(\tau,\frac{1}{2})\neq 0$. We keep this assumption throughout the paper. This case is already deep and challenging. We can formulate our program for any self-dual cuspidal $\tau$ and any $d$. This will be done elsewhere, but we will comment on this more general case in the end of this introduction. 
	
	Let $\Delta(\tau,\ell)$ ($\ell$, a positive integer) denote the Speh representation of $\GL_{4\ell}(\BA)$, attached to $\tau$. See \cite{MW89}. This is the representation spanned by the (multi-)
	residues of Eisenstein series corresponding to the parabolic induction
	from
	$$
	\tau|\det\cdot|^{s_1}\times
	\tau|\det\cdot|^{s_2}\times\cdots\times
	\tau|\det\cdot|^{s_\ell},
	$$
	at the point
	$$
	(\frac{\ell-1}{2},\frac{\ell-3}{2},...,\frac{1-\ell}{2}).
	$$
	Consider Eisenstein series, induced from $\Delta(\tau,\ell)$, on the adelic symplectic group $\Sp_{8\ell}(\BA)$. We will write $\Sp_{8\ell}$ as a matrix group in a standard form, so that the standard Borel subgroup consists of upper triangular matrices. Let $f_{\Delta(\tau,\ell),s}$ be a smooth, holomorphic section of
	\begin{equation}\label{0.9}
		\rho_{\Delta(\tau,\ell),s}=\Ind_{Q_{4\ell}(\BA)}^{\Sp_{8\ell}(\BA)}\Delta(\tau,\ell)|\det\cdot |^s.
	\end{equation}
	We denote the corresponding Eisenstein series by $E(f_{\Delta(\tau,\ell),s})$, and sometimes also by $E^{\Sp_{4\ell}}(f_{\Delta(\tau,\ell),s})$. In \cite{JLZ13}, Theorem 6.2, the poles of the normalized
	Eisenstein series $E^*(f_{\Delta(\tau,\ell),s})$, in $Re(s)\geq 0$, are determined, and they are simple. The largest pole is at $s=\frac{\ell}{2}$ and the remaining poles are $\frac{\ell}{2}-1, \frac{\ell}{2}-2,...$, up to $1$, or $\frac{1}{2}$, according to whether $\ell$ is even or odd, respectively. It is a simple pole of $E(f_{\Delta(\tau,\ell),s})$ (unnormalized), as the section varies. Denote by $\Theta_{\Delta(\tau,\ell)}$ the automorphic representation of $\Sp_{8\ell}(\BA)$ generated by the residues $Res_{s=\frac{\ell}{2}}E(f_{\Delta(\tau,\ell,s})$. We note that
	\begin{prop}\label {prop 1.1}
		The automorphic representation $\Theta_{\Delta(\tau,\ell)}$ is irreducible and square-integrable.
	\end{prop}
	The square-integrability is proved in \cite{JLZ13}, Theorem 6.1. The irreducibility is proved in \cite{L13}, Theorem 7.1.
	The elements $\theta_{\Delta(\tau,\ell)}\in \Theta_{\Delta(\tau,\ell)}$ will be our new ``theta series''. Let $n\leq 2\ell$ (integers). Restrict $\theta_{\Delta(\tau,\ell)}$ to $\Sp_{2n}(\BA)\times \Sp_{4\ell-2n}(\BA)$, where we use the following direct sum embedding $\Sp_{2n}(\BA)\times \Sp_{4\ell-2n}(\BA)\hookrightarrow \Sp_{4\ell}(\BA)$. Let $g\in \Sp_{2n}(\BA)$, $h\in \Sp_{4\ell-2n}(\BA)$. Write $g$ as
	$$
	g=\begin{pmatrix}g_1&g_2\\g_3&g_4\end{pmatrix},
	$$
	where $g_i$ are $n\times n$ matrices. Then
	\begin{equation}\label{0.10}
		i(g,h)=\begin{pmatrix}g_1&&g_2\\&h\\g_3&&g_4\end{pmatrix}.
	\end{equation}
	We will usually simply write $(g,h)$ instead of $i(g,h)$. We use the functions $\theta_{\Delta(\tau,\ell)}(i(g,h))$ as kernel functions. Let $\pi$ be an irreducible, cuspidal, automorphic representation of $\Sp_{2n}(\BA)$. Define, for $h\in \Sp_{4\ell-2n}(\BA)$, 
	\begin{equation}\label{0.11}
		T_\tau^{4\ell-2n}(\varphi_\pi,\theta_{\Delta(\tau,\ell)})(h)=\int\limits_{\Sp_{2n}(F)\backslash \Sp_{2n}(\BA)}\theta_{\Delta(\tau,\ell)}(g,h)\varphi_\pi(g)dg.
	\end{equation}
	We get representations $\Theta_{\Delta(\tau,\ell)}(\pi)$ of $\Sp_{4\ell-2n}(\BA)$. These representations satisfy the tower property (Theorem 5.3 in \cite{G03}).
	\begin{thm}\label{thm 0.6}
		At the first $\ell\geq \frac{n}{2}$, where $\Theta_{\Delta(\tau,\ell)}(\pi)$ is nontrivial, $\Theta_{\Delta(\tau,\ell)}(\pi)$ is cuspidal, in the sense that for all $1\leq r\leq 2\ell-n$, the constant term along the unipotent radical $U^{4\ell-2n}_r$ is identically zero on all elements of $\Theta_{\Delta(\tau,\ell)}(\pi)$.
		
	\end{thm}
	We call the index $\ell$ in Theorem \ref{thm 0.6} the first $\tau$- occurrence of $\pi$. Computations of the correspondence above, at the unramified level, as in Sec. 6 in \cite{G03} show 
	\begin{thm}\label{thm 0.7}
		1. If the first $\tau$-occurrence of $\pi$ is at $\frac{n}{2}\leq \ell<n$, and $\sigma$ is an irreducible (cuspidal) subrepresentation of $\Theta_{\Delta(\tau,\ell)}(\pi)$, then $\pi$ is CAP with respect to
		$$
		\Ind^{\Sp_{2n}(\BA)}_{Q_{2(n-\ell)}(\BA)}\Delta(\tau,n-\ell)|\det\cdot|^{\frac{n-\ell}{2}}\otimes \sigma.
		$$
		2. If the first $\tau$-occurrence of $\pi$ is at $\ell>n$, then any irreducible summand of  $\Theta_{\Delta(\tau,\ell)}(\pi)$ is a CAP representation with respect to  	
		$$
		\Ind^{\Sp_{2n}(\BA)}_{Q_{2(\ell-n)}(\BA)}\Delta(\tau,\ell-n)|\det\cdot|^{\frac{\ell-n}{2}}\otimes \pi.
		$$
		3. If the first $\tau$-occurrence of $\pi$ is at $\ell=n$, then any irreducible summand of  $\Theta_{\Delta(\tau,\ell)}(\pi)$ is nearly equivalent to $\pi$.
		
	\end{thm}
	Thus, the $\Theta_{\Delta(\tau,\ell)}$-correspondence helps detect CAP representations on symplectic groups. In the first case of Theorem \ref{thm 0.7}, the fuctorial lift of $\pi$ to $\GL_{2n+1}(\BA)$ is $\Delta(\tau, 2(n-\ell))\boxplus L(\sigma)$, where $L(\sigma)$ is the functorial lift of $\sigma$ to $\GL_{4\ell-2n+1}(\BA)$. In the second case of the theorem, the functorial lift to $\GL_{4\ell-2n+1}(\BA)$ of each irreducible summand of $\Theta_{\Delta(\tau,\ell)}(\pi)$ is  $\Delta(\tau, 2(\ell-n))\boxplus L(\pi)$.
	
	Consider the question of nonvanishing of $\Theta_{\Delta(\tau,\ell)}(\pi)$. As in \eqref{0.2}, we consider, formally, the inner product of two functions of the form \eqref{0.11},\\
	\\
	$(	T_\tau^{4\ell-2n}(\varphi_\pi,\theta_{\Delta(\tau,\ell)}),T_\tau^{4\ell-2n}(\varphi'_\pi,\theta'_{\Delta(\tau,\ell)}))=$
	\begin{equation}\label{0.12}
		=\int\limits_{[\Sp_{2n}\times \Sp_{2n}]}\varphi_\pi(g_1)\bar{\varphi}'_\pi(g_2)(\int\limits_{[\Sp_{4\ell-2n}]}\theta_{\Delta(\tau,\ell)}(g_1,h)\overline{\theta'_{\Delta(\tau,\ell)}(g_2,h)}dh)dg_1dg_2.
	\end{equation}
	Of course, all of \eqref{0.12} is formal, and we want to make sense out of the r.h.s. of \eqref{0.12}. For this, we need to interpret the inner $dh$-integration. Thus, we would like to have an analog of the multiplicative property \eqref{0.3} of classical theta series, and then we would like to find an analog of the regularized Siegel-Weil formula, Theorem \ref{thm 0.3}, which will interpret and relate the inner product \eqref{0.12}, as in \eqref{0.6} and the proof of Theorem \ref{thm 0.4}, to the generalized doubling integrals \cite{CFGK17}, representing $L(\pi\times\tau,s)$. 
	
	Assume that $\pi$ is as in the first case of Theorem \ref{thm 0.7}. In particular, $\frac{n}{2}\leq \ell<n$. Let $S$ be a finite set of places of $F$, containing the archimedean places, outside which $\pi$ is unramified. Assume also that $L^S(\sigma\times \tau,s)$ is holomorphic and nonzero at $s=n-\ell+\frac{1}{2}$, for example, when $\sigma$ is generic. Then $L^S(\pi\times\tau,s)$ has a simple pole at $s=n-\ell+\frac{1}{2}$. Let us represent $L(\pi\times\tau,s+\frac{1}{2})$ by the generalized doubling integrals. These have the form
	\begin{equation}\label{0.13}
		\mathcal{L}(\varphi_\pi,\varphi'_\pi,f_{\Delta(\tau,2n),s})=	\int\limits_{[\Sp_{2n}\times \Sp_{2n}]}\varphi_\pi(g_1)\bar{\varphi}'_\pi(g_2)E^{\psi_{U_{2n}}}(f_{\Delta(\tau,2n),s})(t(g_1,g_2))dg_1dg_2,
	\end{equation}
	where $E^{\psi_{U_{2n}}}(f_{\Delta(\tau,2n),s})$ denotes the following Fourier coefficient along $U_{2n}$,
	$$
	E^{\psi_{U_{2n}}}(f_{\Delta(\tau,2n),s})(x)=\int\limits_{U_{2n}(F)\backslash U_{2n}(\BA)}E(f_{\Delta(\tau,2n),s})(ux)\psi^{-1}_{U_{2n}}(u)du,
	$$
	and $\psi_{U_{2n}}$ is the following character of $U_{2n}(\BA)$. Let
	$$
	u=\begin{pmatrix}I_{2n}&y&z\\&I_{4n}&y'\\&&I_{2n}\end{pmatrix}\in U_{2n}(\BA).
	$$
	Then, when we write $y=(y_1,y_2,y_3)$, $y_1,y_3\in M_{2n\times n}(\BA)$, 
	\begin{equation}\label{0.14}
		\psi_{U_{2n}}(u)=\psi(tr(y_1+y_3)).
	\end{equation}
	Finally, for $g_i\in \Sp_{2n}(\BA)$, $i=1,2$, and $g_1=\begin{pmatrix}a_1&b_1\\c_1&d_1\end{pmatrix}$, with $n\times n$ blocks,
	$$
	t(g_1,g_2)=\begin{pmatrix}g_1\\&a_1&&b_1\\&&g_2\\&c_1&&d_1\\&&&&g_1^*\end{pmatrix}.
	$$ 
	We conclude that there exist data such that
	\begin{equation}\label{0.15}
		\int\limits_{[\Sp_{2n}\times \Sp_{2n}]}\varphi_\pi(g_1)\bar{\varphi}'_\pi(g_2)Res_{s=n-\ell}(E(f_{\Delta(\tau,2n),s}))^{\psi_{U_{2n}}}(t(g_1,g_2))dg_1dg_2\neq 0
	\end{equation}
	This line of thought suggests that the inner $dh$-integral in the r.h.s. of \eqref{0.12} should be related to the residue inside the integral \eqref{0.15}, namely, for some choice of data
	\begin{equation}\label{0.16}
		\int\limits_{[\Sp_{4\ell-2n}]}\theta_{\Delta(\tau,\ell)}(g_1,h)\overline{\theta'_{\Delta(\tau,\ell)}(g_2,h)}dh=Res_{s=n-\ell}(E(f_{\Delta(\tau,2n),s}))^{\psi_{U_{2n}}}(t(g_1,g_2)).
	\end{equation}
	We are still at the formal level since the l.h.s. of \eqref{0.16} may diverge. 	Let us apply, still formally, the $\Theta_{\Delta(\tau,2n+\ell)}$-correspondence \eqref{0.11} to the non-cuspidal representation $\Theta_{\Delta(\tau,\ell)}$, that is, for $ h\in \Sp_{8n}(\BA)$,
	\begin{equation}\label{0.17}
		T^{8n}(\theta_{\Delta(\tau,\ell)}, \theta_{\Delta(\tau,2n+\ell)})(h)=\int\limits_{\Sp_{4\ell}(F)\backslash \Sp_{4\ell}(\BA)}\theta_{\Delta(\tau,2n+\ell)}(g,h)\theta_{\Delta(\tau,\ell)}(g)dg.
	\end{equation}
	The proof of Theorem \ref{thm 0.7} shows that the unramified parameters of the representation of $\Sp_{8n}(\BA)$ generated by the functions $T^{8n}(\theta_{\Delta(\tau,\ell)}, \theta_{\Delta(\tau,2n+\ell)})$, that is $\Theta_{\Delta(\tau,2n+\ell)}(\Theta_{\Delta(\tau,\ell)})$, are identical to those of the representation generated by the residues $Res_{s=n-\ell}(E(f_{\Delta(\tau,2n),s}))$. Thus, we expect the following crude, formal type of a  Siegel-Weil formula,
	\begin{conj}\label{conj 0.8}{\bf (Siegel-Weil formula, crude form)}
		Assume that $1\leq \ell<n$. Given $\theta_{\Delta(\tau,\ell)}, \theta_{\Delta(\tau,2n+\ell)}$, there is a section $f'_{\Delta(\tau,2n),s}$, such that, for $h\in \Sp_{8n}(\BA)$,
		$$
		\int\limits_{\Sp_{4\ell}(F)\backslash \Sp_{4\ell}(\BA)}\theta_{\Delta(\tau,2n+\ell)}(g,h)\theta_{\Delta(\tau,\ell)}(g)dg	=Res_{s=n-\ell}(E(f'_{\Delta(\tau,2n),s}))(h).
		$$
	\end{conj}

	One question that immediately arises is how does $f'_{\Delta(\tau,2n),s}$ depend on $\theta_{\Delta(\tau,\ell)}$, $\theta_{\Delta(\tau,2n+\ell)}$? Is there an analog of the Siegel-Weil section as in \eqref{0.5}? We propose the following analog which can be tracked down to a certain term in the Fourier expansion of the constant term along the Siegel radical of the Eisenstein series \eqref{2.10} in Theorem \ref{thm 2.4}, which turns out to be an analog of the Eisenstein series \eqref{0.7}. Consider the constant term of $\theta_{\Delta(\tau,2n+\ell)}$ along the unipotent radical $U_{4n}=U_{2n}^{4n+2\ell}$,
	$$
	\theta^{U_{4n}}_{\Delta(\tau,2n+\ell)}(x)=\int\limits_{U_{4n}(F)\backslash U_{4n}(\BA)}\theta_{\Delta(\tau,2n+\ell)}(ux)dx.
	$$
	Define
	\begin{equation}\label{0.17'}
		\Phi(\theta_{\Delta(\tau,\ell)}, \theta_{\Delta(\tau,2n+\ell)})(h)=\int\limits_{\Sp_{4\ell}(F)\backslash \Sp_{4\ell}(\BA)}\theta^{U_{4n}}_{\Delta(\tau,2n+\ell)}(i^*(g,h)w)\theta_{\Delta(\tau,\ell)}(g)dg,
	\end{equation}
    where, in the notation of \eqref{0.10}, $i^*(g,h)$ is obtained from $i(g,h)$ by switching $g$ and $h$, and $w$ is a Weyl element which does this switch by conjugating $i(g,h)$. Thus, writing $h=\begin{pmatrix}h_1&h_2\\h_3&h_4\end{pmatrix}$,
    $$
    i^*(g,h)=\begin{pmatrix}h_1&&h_2\\&g\\h_3&&h_4\end{pmatrix}.
    $$
    The integral \eqref{0.17'} converges absolutely. This follows from Propositions \ref{prop 1.2}, \ref{prop 1.3}. Define
    \begin{equation}\label{0.17''}
    	\Phi(\theta_{\Delta(\tau,\ell)}, \theta_{\Delta(\tau,2n+\ell)},s)(h)=a(h)^{s+n-\ell}\Phi(\theta_{\Delta(\tau,\ell)}, \theta_{\Delta(\tau,2n+\ell)})(h),
    \end{equation}
where $a(h)$ is as in the definition of the Siegel-Weil section \eqref{0.5}. One can check that $\Phi(\theta_{\Delta(\tau,\ell)}, \theta_{\Delta(\tau,2n+\ell)},s)$ is a section of $\rho_{\Delta(\tau,2n),s}$. This is our proposed analog of a Siegel-Weil section. Note the analogy with the Siegel-Weil section \eqref{0.5'}. It depends on the constant term \eqref{0.5''} of $\theta^\Phi_{\psi,8mn}$ along $U_{4mn}$ and the section \eqref{0.17''} depends on the constant term $\theta^{U_{4n}}_{\Delta(\tau,2n+\ell)}$.

Consider the Eisenstein series on $\Sp_{8n}(\BA)$ attached to the section\\
 $\Phi(\theta_{\Delta(\tau,\ell)}, \theta_{\Delta(\tau,2n+\ell)},s)$, $E(\Phi(\theta_{\Delta(\tau,\ell)}, \theta_{\Delta(\tau,2n+\ell)},s))$. Then a more precise form of the last conjecture, still at the formal level, is
\begin{conj}\label{conj 0.8'}
	1. Assume that $1\leq \ell\leq n$. Then $E(\Phi(\theta_{\Delta(\tau,\ell)}, \theta_{\Delta(\tau,2n+\ell)},s))$ is holomorphic at $s=\ell-n$, and, up to a possible constant, for all $h\in \Sp_{8n}(\BA)$,
	$$
	\int\limits_{\Sp_{4\ell}(F)\backslash \Sp_{4\ell}(\BA)}\theta_{\Delta(\tau,2n+\ell)}(g,h)\theta_{\Delta(\tau,\ell)}(g)dg	=Value_{s=\ell-n}E(\Phi(\theta_{\Delta(\tau,\ell)}, \theta_{\Delta(\tau,2n+\ell)},s))(h).
	$$
	2. Let $M_s$ denote the intertwining operator from  $\rho_{\Delta(\tau,2n),s}$ to  $\rho_{\Delta(\tau,2n),-s}$, and let $M^*_s$ denote the normalized intertwining operator. Then $M_s(\Phi(\theta_{\Delta(\tau,\ell)}, \theta_{\Delta(\tau,2n+\ell)},s))$ has a zero at $s=\ell-n$. Denote
	$$
	\Phi^*(\theta_{\Delta(\tau,\ell)}, \theta_{\Delta(\tau,2n+\ell)})=M^*_s(\Phi(\theta_{\Delta(\tau,\ell)}, \theta_{\Delta(\tau,2n+\ell)},s))\Big|_{s=\ell-n},
	$$
	$$
	\Phi^*(\theta_{\Delta(\tau,\ell)}, \theta_{\Delta(\tau,2n+\ell)},s)(h)=a(h)^{s+\ell-n}\Phi^*(\theta_{\Delta(\tau,\ell)}, \theta_{\Delta(\tau,2n+\ell)})(h).
	$$
	This is a section of $\rho_{\Delta(\tau,2n),s}$. Up to a possible constant, we have, for all\\ $h\in \Sp_{8n}(\BA)$,
	$$
	\int\limits_{\Sp_{4\ell}(F)\backslash \Sp_{4\ell}(\BA)}\theta_{\Delta(\tau,2n+\ell)}(g,h)\theta_{\Delta(\tau,\ell)}(g)dg	=Res_{s=n-\ell}E(\Phi^*(\theta_{\Delta(\tau,\ell)}, \theta_{\Delta(\tau,2n+\ell)},s))(h).
	$$
	3. Assume that $n<\ell\leq 2n$. Then, up to a possible constant, for all $h\in \Sp_{8n}(\BA)$,
	$$
	\int\limits_{\Sp_{4\ell}(F)\backslash \Sp_{4\ell}(\BA)}\theta_{\Delta(\tau,2n+\ell)}(g,h)\theta_{\Delta(\tau,\ell)}(g)dg	=Res_{s=\ell-n}E(\Phi(\theta_{\Delta(\tau,\ell)}, \theta_{\Delta(\tau,2n+\ell)},s))(h).
	$$
	
\end{conj}
	In view of \eqref{0.16}, we expect
	\begin{conj}\label{conj 0.9}
		Assume that $\frac{n}{2}\leq \ell<n$. Given $\theta_{\Delta(\tau,2n+\ell)}, \theta_{\Delta(\tau,\ell)}$, there exist $\theta'_{\Delta(\tau,2n+\ell)}, \theta''_{\Delta(\tau,\ell)}$, such that
		$$
		\int\limits_{U_{2n}^{8n}(F)\backslash U_{2n}^{8n}(\BA)}\int\limits_{\Sp_{4\ell}(F)\backslash \Sp_{4\ell}(\BA)}\theta_{\Delta(\tau,2n+\ell)}(g,u\cdot t(g_1,g_2))\theta_{\Delta(\tau,\ell)}(g)\psi^{-1}_{U_{2n}}(u)dgdu=
		$$
		$$
		=\int\limits_{\Sp_{4\ell-2n}(F)\backslash \Sp_{4\ell-2n}(\BA)}\theta'_{\Delta(\tau,\ell)}(g_1,h)\overline{\theta''_{\Delta(\tau,\ell)}(g_2,h)}dh.
		$$
	\end{conj}
	We view Conjecture \ref{conj 0.9} as an analog of the multiplicative property \eqref{0.3} of theta series. As an application, we will have the following analog of Theorem \ref{thm 0.4}.
	\begin{thm}\label{thm 0.10}
		Let $\pi$ be an irreducible, cuspidal, automorphic representation of $\Sp_{2n}(\BA)$. Assume that $L^S(\pi\times \tau,s)$ has its largest pole at $s=n-\ell+\frac{1}{2}$, where $\ell<n$. Then $\ell\geq \frac{n}{2}$ and $\Theta_{\Delta(\tau,\ell)}(\pi)$ is nonzero and cuspidal, so that by Theorem \ref{thm 0.7}, $\pi$ is CAP with respect to 
		$$
		\Ind^{\Sp_{2n}(\BA)}_{Q_{2(n-\ell)}(\BA)}\Delta(\tau,n-\ell)|\det\cdot|^{\frac{n-\ell}{2}}\otimes \sigma,
		$$
		where $\sigma$ is an irreducible (cuspidal) subrepresentation of $\Theta_{\Delta(\tau,\ell)}(\pi)$.
	\end{thm}
	\begin{proof}
		We sketch the main steps, based on the conjectures above. By \cite{GS22}, Prop. 3.3, we have $n-\ell\leq \frac{n}{2}$, and hence $\ell\geq \frac{n}{2}$. As we explained before, we conclude that \eqref{0.15} is satisfied. Assuming 	that the section $f'_{\Delta(\tau,2n),s}$ obtained in Conjecture \ref{conj 0.8} is sufficiently general, we conclude from Conjecture \ref{conj 0.8} and then Conjecture \ref{conj 0.9}, that there exist $\theta'_{\Delta(\tau,2n+\ell)}, \theta''_{\Delta(\tau,\ell)}$, such that
		$$
		\int\limits_{[\Sp_{2n}\times \Sp_{2n}]}\int\limits_{[\Sp_{4\ell-2n}] }\varphi_\pi(g_1)\bar{\varphi}'_\pi(g_2)\theta'_{\Delta(\tau,\ell)}(g_1,h)\overline{\theta''_{\Delta(\tau,\ell)}(g_2,h)}dhdg_1dg_2\neq 0.
		$$
		This is \eqref{0.12}, and hence the $\Theta_{\Delta(\tau,\ell)}$-correspondence of $\pi$ to $\Sp_{4n-2\ell}(\BA)$ is nontrivial. Since $s=n-\ell+\frac{1}{2}$ is the largest pole of $L^S(\pi\times\tau,s)$, one can show that $\Theta_{\Delta(\tau,\ell)}(\pi)$ is cuspidal. Now we apply Theorem \ref{thm 0.7}.	
		
	\end{proof}

	Our {\bf starting point} in realizing the program outlined above will be a regularization of the integral in Conjecture \ref{conj 0.8}. We will show that the regularized integral is equal to the residue of an Eisenstein series on $\Sp_{8n}(\BA)$ corresponding to an Eisenstein series induced from $\Delta(\tau,\ell)|\det\cdot|^\zeta\otimes \Theta(\tau,2n-\ell)$ at $\zeta=\frac{\ell}{2}$. This is an analog of Theorem \ref{thm 0.2}. This will enable us to formulate a precise version of Conjecture \ref{conj 0.8'}. We will state our main theorems in the next section.
	
	\subsection{Generalizations.}
	We comment on generalizations of the above conjectures. First, we may consider any irreducible, cuspidal, automorphic representation $\tau$ of $GL_{2d}({\bf A})$, such that $L(\tau,1/2)\ne 0$ and $L(\tau,\wedge^2,s)$ has a pole at $s=1$. We write down the analog of the integral in Conjecture \ref{conj 0.8}. Consider the automorphic representations $\Theta_{\Delta(\tau,\ell)}$, $\Theta_{\Delta(\tau,2n+(2d-1)\ell)}$ of $\Sp_{4d\ell}(\BA)$,  $\Sp_{4d(2n+(2d-1)\ell)}(\BA)$, respetively. Let $U=U_{(4d\ell)^{d-1}}$ be the unipotent radical of the parabolic subgroup of $\Sp_{4d(2n+(2d-1)\ell)}$, whose Levi part is isomorphic to $GL_{4d\ell}^{d-1}\times Sp_{4d(2n+\ell)}$. Consider the unipotent orbit of  $\Sp_{4d(2n+(2d-1)\ell)}$ corresponding to the partition $((2d-1)^{4d\ell}1^{8dn})$. Then one can define a character $\psi_U$ of $U(\BA)$, trivial on $U(F)$, corresponding to this orbit, such that the stabilizer of this character, inside $GL_{4d\ell}^{d-1}(\BA)\times Sp_{4d(2n+\ell)}(\BA)$, is isomorphic to $Sp_{4d\ell}(\BA)\times Sp_{8dn}(\BA)$. For $\theta_{\Delta(\tau,2n+(2d-1)\ell)}\in \Theta_{\Delta(\tau,2n+(2d-1)\ell)}$, denote
	$$
	\theta_{\Delta(\tau,2n+(2d-1)\ell)}^{\psi_U}(x)=\int\limits_{U(F)\backslash U(\BA)}\theta_{\Delta(\tau,2n+(2d-1)\ell)}(ux)\psi^{_1}_U(u)du.
	$$	
	Recall that the Eisenstein series on $Sp_{8dn}({\BA})$, $E(f_{\Delta(\tau,2n),s})$ has its positive poles at $s=1,2,...,n$.	Our analog of Conjecture \ref{conj 0.8} (crude form and formal) is
	\begin{conj}\label{conj 0.10}
	Let $\theta_{\Delta(\tau,\ell)}\in \Theta_{\Delta(\tau,\ell)}$,  $\theta_{\Delta(\tau,2n+(2d-1)\ell)}\in \Theta_{\Delta(\tau,2n+(2d-1)\ell)}$.\\ 
	1. Assume that $1\leq \ell\leq n$. Then there is a section $f_{\Delta(\tau,2n),s}$, such that $E(f_{\Delta(\tau,2n),s})$ is holomorphic at $s=\ell-n$ and
	$$
			\int\limits_{Sp_{4d\ell}(F)\backslash Sp_{4d\ell}({\BA})}
				\theta_{\Delta(\tau,2n+(2d-1)\ell)}^{\psi_U}((h,g))\theta_{\Delta(\tau,\ell)}(g)dg=Value_{s=\ell-n}E(f_{\Delta(\tau,2n),s},h).
		$$
		2. Assume that $n<\ell\leq 2n$. Then there is a section $f_{\Delta(\tau,2n),s}$, such that 
		$$
		\int\limits_{Sp_{4d\ell}(F)\backslash Sp_{4d\ell}({\BA})}
		\theta_{\Delta(\tau,2n+(2d-1)\ell)}^{\psi_U}((h,g))\theta_{\Delta(\tau,\ell)}(g)dg=Res_{s=\ell-n}E(f_{\Delta(\tau,2n),s},h).
		$$ 
	\end{conj}
The previous case corresponds to $d=1$. The technical difference and difficulty here is the presence of the unipotent integration.	
	
	The next generalization to consider is for an irreducible, cuspidal, automorphic representation $\tau$ of $\GL_d(\BA)$, such that $L(\tau,1/2)\ne 0$ and $L(\tau,sym^2,s)$ has a pole at $s=1$. In this case, we replace our ``theta series'' corresponding to $\tau$ by the following residual representations of the double covers $\Sp^{(2)}_{2dr}(\BA)$ of $\Sp_{2dr}(\BA)$. Consider the representation
	$$
	\rho^{(2)}_{\Delta(\tau,r),s;\psi}=\Ind_{Q^{(2)}_{rd}(\BA)}\Delta(\tau,r)|\det\cdot|^s\gamma_\psi,
	$$
	where $\gamma_\psi$ is the Weil factor composed with the determinant. Let $	f_{\Delta(\tau,r),s;\psi}$ be a $K_{2dr}$-finite, holomorphic section of $\rho^{(2)}_{\Delta(\tau,r),s;\psi}$, and denote by $E^{(2)}(	f_{\Delta(\tau,r),s;\psi})$ the corresponding Eisenstein series on $\Sp^{(2)}_{2dr}(\BA)$. Then, as the section varies its positive poles are at $s=\frac{r}{2},\frac{r}{2}-1,\frac{r}{2}-2,...$. We let $\Theta^{(2)}_{\Delta(\tau,r),\psi}$ denote the automorphic representation of $\Sp^{(2)}_{2dr}(\BA)$ generated by the residues at the largest pole $Res_{s=\frac{r}{2}}E^{(2)}(f_{\Delta(\tau,r),s;\psi})$. See \cite{GS22}. Similarly, we can write theta representations corresponding to $\tau$, self-dual and cuspidal, on split orthogonal groups, general linear groups and also on higher covers of all the groups above, and outline a similar program. We will deal with these in future works.

	We finish this introduction with two facts that we need on the representations $\Theta_{\Delta(\tau,m)}$. Now, we go back to our initial case of study, namely that $\tau$ is an irreducible, cuspidal, automorphic representation of $\GL_2(\BA)$, such that $L(\tau,\wedge^2,s)$ has a pole at $s=1$ and $L(\tau,\frac{1}{2})\neq 0$. The representation $\Theta_{\Delta(\tau,m)}$ is of $\Sp_{4m}(\BA)$. 
	For $f\in \Theta_{\Delta(\tau,m)}$, consider the constant term of $f$ along the unipotent radical $U_r$, $r\leq 4m$,
	\begin{equation}\label{1.1}
		f^{U_r}(g)=\int\limits_{U_r(F)\backslash U_r(\BA)}f(ug)du.
	\end{equation}

	\begin{prop}\label{prop 1.2}
		The constant term \eqref{1.1} is zero on $\Theta_{\Delta(\tau,m)}$, unless $r=2i$, $1\leq i\leq m$,  and then, for each $f\in \Theta_{\Delta(\tau,m)}$, the constant term $f^{U_{2i}}$, as a function on $\Sp_{4m}(\BA)$, lies in the space of 
		$$
		\Ind_{Q_{2i}(\BA)}^{\Sp_{4m}(\BA)}\Delta(\tau,i)|\det\cdot|^{-m+\frac{i}{2}}\otimes \Theta_{\Delta(\tau,m-i)}.
		$$
	\end{prop}
	
	This is Lemma 2.3 in \cite{L13} (in the special case of $\GL_2$).
	
	Recall that Fourier coefficients supported by automorphic forms correspond to nilpotent orbits. See \cite{GRS03}.
	In the case of the automorphic forms in the space of $\Theta_{\Delta(\tau,m)}$, we have
	
	\begin{prop}\label{prop 1.3}
		There is a unique maximal nilpotent orbit, in the Lie algebra of $\Sp_{4m}$ over the algebraic closure of $F$, attached to Fourier coefficients admitted by $\Theta_{\Delta(\tau,m)}$. This is the orbit which corresponds to the partition $((2)^{2m})$.
	\end{prop}
	This was proved (in general) by Ginzburg in \cite{G08}. See \cite{L13} for a detailed proof. See also \cite{GS22}.

	\section{Statement of the main theorems}
	
	As we explained in the end of the introduction, we start with the integral in Conjecture \ref{conj 0.8}. We will replace $2n$ by $m$ and assume that $1\leq \ell\leq m$. Consider $\Theta_{\Delta(\tau,m+\ell)}$ and restrict it to $\Sp_{4\ell}(\BA)\times \Sp_{4m}(\BA)$, where, as in \eqref{0.10}, for $g\in \Sp_{4\ell}(\BA)$, $h\in \Sp_{4m}(\BA)$, writing $g$ as
	$$
	g=\begin{pmatrix}g_1&g_2\\g_3&g_4\end{pmatrix},
	$$
	where $g_i$ are $2\ell\times 2\ell$ matrices, we embed $(g,h)$ inside $\Sp_{4(m+\ell)}(\BA)$ by
	\begin{equation}\label{2.1}
		i(g,h)=\begin{pmatrix}g_1&&g_2\\&h\\g_3&&g_4\end{pmatrix},
	\end{equation}
	and when convenient, we simply denote $i(g,h)=(g,h)$.
	
	Let $\theta_{\Delta(\tau,m+\ell)}\in \Theta_{\Delta(\tau,m+\ell)}$ and $\theta_{\Delta(\tau,\ell)}\in \Theta_{\Delta(\tau,\ell)}$. Define, for $h\in \Sp_{4m}(\BA)$, the following integral whenever convergent
	\begin{equation}\label{2.2}
		E(\theta_{\Delta(\tau,m+\ell)}, \theta_{\Delta(\tau,\ell)};h)=\int\limits_{\Sp_{4\ell}(F)\backslash \Sp_{4\ell}(\BA)}\theta_{\Delta(\tau,m+\ell)}(i(g,h))\theta_{\Delta(\tau,\ell)}(g)dg.
	\end{equation}
	We will show how to regularize it and turn it into an absolutely convergent integral. The regularization is similar to the one carried out by Ichino for the regularization of the Siegel-Weil formula \cite{I01}. We will choose a finite place $v$, where $\Theta_{\Delta(\tau,m+\ell),v}$ is unramified and apply an element $\xi_v\in \mathcal{H}(\Sp_{4m}(F_v)// K_{4m,v})$ in the spherical Hecke algebra of $\Sp_{4m}(F_v)$. Denote
	\begin{equation}\label{2.3}
		((1\otimes\xi_v)\ast \theta_{\Delta(\tau,m+\ell)})(x)=\int\limits_{\Sp_{4m}(F_v)}\xi_v(h_v)\theta_{\Delta(\tau,m+\ell)}(x(1,h_v))dh_v.
	\end{equation}
	We will prove
	\begin{thm}\label{thm 2.1}
		There is a function $\xi_{m,\ell,v}\in \mathcal{H}(\Sp_{4m}(F_v)// K_{4m,v})$, depending on $\tau_v,m,\ell$, such that the function on $\Sp_{4\ell}(\BA)$,
		$$
		g\mapsto ((1\otimes\xi_{m,\ell,v})\ast\theta_{\Delta(\tau,m+\ell)})(g,h)
		$$
		is rapidly decreasing, uniformly in $h$ varying in bounded subsets of $\Sp_{4m}(F)\backslash\Sp_{4m}(\BA)$. Let $m\geq 2\ell$. There is a nonzero complex number $c_{m,\ell,v}$, depending on $\tau_v,m,\ell$, such that the following holds. Assume that $\theta_{\Delta(\tau,m+\ell)}$ is right $i(K_{4\ell,v}\times K_{4m,v})$-invariant and the integral \eqref{2.2} converges absolutely. Then 
		\begin{equation}\label{2.4}
			E((1\otimes\xi_{m,\ell,v})\ast\theta_{\Delta(\tau,m+\ell)}, \theta_{\Delta(\tau,\ell)};h)=c_{m,\ell,v}E(\theta_{\Delta(\tau,m+\ell)}, \theta_{\Delta(\tau,\ell)};h).
		\end{equation}
	\end{thm}
	Thus, we regularize \eqref{2.2} by
	\begin{equation}\label{2.5}
		E_{reg}(\theta_{\Delta(\tau,m+\ell)}, \theta_{\Delta(\tau,\ell)};h)=\frac{1}{c_{m,\ell,v}}E((1\otimes\xi_{m,\ell,v})\ast\theta_{\Delta(\tau,m+\ell)}, \theta_{\Delta(\tau,\ell)};h).
	\end{equation}	
	Note that we only need that the finite place $v$ is such that $\theta_{\Delta(\tau,m+\ell)}$ is right $i(K_{4\ell,v}\times K_{4m,v})$-invariant. We will define the function $\xi_{m,\ell,v}$ in Section 3, and prove the rapid decrease in Section 5.
	
	In Section 6, we will identify the regularized integral \eqref{2.5} as a residual Eisenstein series on  $\Sp_{4m}(\BA)$.
	\begin{thm}\label{thm 2.2}
		Assume that $m\geq 2\ell$, and assume that $\theta_{\Delta(\tau,m+\ell)}$ is $i(K_{4\ell}\times K_{4m})$-finite. There is a $K_{4m}$-finite, holomorphic section $\varphi_{\Delta(\tau,\ell)|\det\cdot|^s\otimes \Theta_{\Delta(\tau,m-\ell)}}$ of    
		$$
		\Ind_{Q_{2\ell}^{4m}(\BA)}^{\Sp_{4m}(\BA)}\Delta(\tau,\ell)|\det\cdot|^s\otimes \Theta_{\Delta(\tau,m-\ell)}, 
		$$
		depending on $\theta_{\Delta(\tau,m+\ell)}, \theta_{\Delta(\tau,\ell)}$, such that
		\begin{equation}\label{2.6}
			E_{reg}(\theta_{\Delta(\tau,m+\ell)}, \theta_{\Delta(\tau,\ell)};h)=Res_{s=\frac{\ell}{2}}E(\varphi_{\Delta(\tau,\ell)|\det\cdot|^s\otimes \Theta_{\Delta(\tau,m-\ell)}},h),
		\end{equation}
		where $	E(\varphi_{\Delta(\tau,\ell)|\det\cdot|^s\otimes \Theta_{\Delta(\tau,m-\ell)}},h)$ denotes the Eisenstein series on $\Sp_{4m}(\BA)$ corresponding to $\varphi_{\Delta(\tau,\ell)|\det\cdot|^s\otimes \Theta_{\Delta(\tau,m-\ell)}}$.
	\end{thm}
	Let $f_{\Delta(\tau,\ell),s}$ be a smooth, holomorphic section of $\rho_{\Delta(\tau,\ell),s}$ (\eqref{0.9}), and consider the Eisenstein series on $\Sp_{4\ell}(\BA)$, $E(f_{\Delta(\tau,\ell),s})$, corresponding to $f_{\Delta(\tau,\ell),s}$.
	We know that it has a simple pole at $s=\frac{\ell}{2}$, as the section varies. Let us take
	\begin{equation}\label{2.7}
		\theta_{\Delta(\tau,\ell)}=Res_{s=\frac{\ell}{2}}E(f_{\Delta(\tau,\ell),s}).
	\end{equation}
	Then we can rewrite \eqref{2.5} as
	\begin{equation}\label{2.8}
		E_{reg}(\theta_{\Delta(\tau,m+\ell)}, \theta_{\Delta(\tau,\ell)};h)=\frac{1}{c_{m,\ell,v}}E((1\otimes\xi_{m,\ell,v})\ast\theta_{\Delta(\tau,m+\ell)}, Res_{s=\frac{\ell}{2}}E(f_{\Delta(\tau,\ell),s});h)=
	\end{equation}	
	$$
	\frac{1}{c_{m,\ell,v}}Res_{s=\frac{\ell}{2}}(\int\limits_{\Sp_{4\ell}(F)\backslash \Sp_{4\ell}(\BA)}(1\otimes\xi_{m,\ell,v})\ast\theta_{\Delta(\tau,m+\ell)}(i(g,h))E(f_{\Delta(\tau,\ell),s};g)dg).
	$$	
	Note that Theorem \ref{thm 2.1} implies that the last integral converges absolutely, whenever $E(f_{\Delta(\tau,\ell),s}$ is holomorphic.	Theorem \ref{thm 2.2} will follow from
	\begin{thm}\label{thm 2.3}
		Assume that $\theta_{\Delta(\tau,m+\ell)}$ is $i(K_{4\ell}\times K_{4m})$-finite. There is a $K_{4m}$-finite, holomorphic section $\varphi_{\Delta(\tau,\ell)|\det\cdot|^s\otimes \Theta_{\Delta(\tau,m-\ell)}}$ of    
		$$
		\Ind_{Q_{2\ell}^{4m}(\BA)}^{\Sp_{4m}(\BA)}\Delta(\tau,\ell)|\det\cdot|^s\otimes \Theta_{\Delta(\tau,m-\ell)}, 
		$$
		depending on $\theta_{\Delta(\tau,m+\ell)}, f_{\Delta(\tau,\ell),\frac{\ell}{2}}$, such that
		\begin{equation}\label{2.9}
			\int\limits_{\Sp_{4\ell}(F)\backslash \Sp_{4\ell}(\BA)}(1\otimes\xi_{m,\ell,v})\ast\theta_{\Delta(\tau,m+\ell)}(i(g,h))E(f_{\Delta(\tau,\ell),s};g)dg=
		\end{equation}
		$$
		=E(\varphi_{\Delta(\tau,\ell)|\det\cdot|^s\otimes \Theta_{\Delta(\tau,m-\ell)}},h).
		$$
	\end{thm}
	
	We prove that Theorem \ref{thm 2.3} is valid, for any $\ell\leq m$. In fact, we also get the analogs of the previous two theorems for $m\leq 2\ell-1$. We just need to specify a little more our choice of the place $v$. Since $\tau_v$ is unramified and self-dual, $\tau_v$ is induced from a character of the Borel subgroup
	$$
	\begin{pmatrix}a&x\\0&b\end{pmatrix}\mapsto \chi(ab^{-1}), \ \ a, b\in F_v^*, x\in F_v,
	$$
	where $\chi=\chi_v$ is an unramified character of $F_v^*$. We assume that $\chi^2\neq 1$. We can choose such a place $v$, since otherwise, at all finite places $v$, where $\tau_v$ is unramified, we would get that $L(\tau_v\times\tau_v,s)=(1-q_v^{-s})^{-4}$, and hence the partial $L$-function, away from the archimedean places, $L^{S_\infty}(\tau\times\tau,s)=L^{S_\infty}(1,s)^4$ has a pole of order four at $s=1$, which is impossible, since $\tau$ is cuspidal. We prove
	
	\begin{thm}\label{thm 2.4}
		Assume that the unramified character $\chi$, corresponding to $\tau_v$, as above, is not quadratic. Then the function $\xi_{m,\ell,v}$ satisfies the following property. There is a polynomial $P(x, y)\in \BC[x,y]$, depending on $\tau_v, \ell,m$, such that, for $m\geq 2\ell$, $P(q_v^{-\frac{\ell}{2}},q_v^{\frac{\ell}{2}})=c_{m,\ell,v}\neq 0$, and for $\ell\leq m\leq 2\ell-1$, $P(q_v^{-s},q_v^s)$ has a simple zero at $s=\frac{\ell}{2}$. In the notation of Theorem \ref{thm 2.3}, denote\\
		\\
		$\mathcal{E}(\theta_{\Delta(\tau,m+\ell)},f_{\Delta(\tau,\ell),s};h)=$
		\begin{equation}\label{2.10}
			\frac{1}{P(q_v^{-s},q_v^s)}\int\limits_{\Sp_{4\ell}(F)\backslash \Sp_{4\ell}(\BA)}(1\otimes\xi_{m,\ell,v})\ast\theta_{\Delta(\tau,m+\ell)}(i(g,h))E(f_{\Delta(\tau,\ell),s};g)dg.
		\end{equation}
		Then $\mathcal{E}(\theta_{\Delta(\tau,m+\ell)},f_{\Delta(\tau,\ell),s};h)$ is an Eisenstein series on $\Sp_{4m}(\BA)$, corresponding to $	\Ind_{Q_{2\ell}^{4m}(\BA)}^{\Sp_{4m}(\BA)}\Delta(\tau,\ell)|\det\cdot|^s\otimes \Theta_{\Delta(\tau,m-\ell)}$. At $s=\frac{\ell}{2}$, it has at most a simple pole, when $m\geq 2\ell$, and at most a double pole, when $\ell\leq m\leq 2\ell-1$. 	
		
	\end{thm}
	
	The polynomial $P(x,y)$ is written in \eqref{5.22} and the properties of $P(q_v^{-s},q_v^s)$ at $s=\frac{\ell}{2}$ are shown in Lemma \ref{lem 5.5}. This together with Theorem \ref{thm 2.3} will yield Theorem \ref{thm 2.4}.
	
	The series \eqref{2.10} is the analog of \eqref{0.7}. As in Theorem \ref{thm 0.3}, we regularize the integral \eqref{2.2} as follows. Consider the Laurent expansion of \eqref{2.10} at $s=\frac{\ell}{2}$. The leading term when $m\geq 2\ell$ is the residue $B_{-1}(h, \theta_{\Delta(\tau,m+\ell)},f_{\Delta(\tau,\ell),\frac{\ell}{2}})$. When $m\leq 2\ell-1$, the leading term is the coefficient of $(s-\frac{\ell}{2})^{-2}$, which we denote by $B_{-2}(h, \theta_{\Delta(\tau,m+\ell)},f_{\Delta(\tau,\ell),\frac{\ell}{2}})$. Thus, we generalize \eqref{2.8}, for all $1\leq \ell\leq m$, by	
	\begin{equation}\label{2.11}
		E_{reg}(\theta_{\Delta(\tau,m+\ell)}, \theta_{\Delta(\tau,\ell)};h)=\begin{cases}B_{-1}(h, \theta_{\Delta(\tau,m+\ell)},f_{\Delta(\tau,\ell),\frac{\ell}{2}}),\ \  2\ell\leq m\\
			B_{-2}(h, \theta_{\Delta(\tau,m+\ell)},f_{\Delta(\tau,\ell),\frac{\ell}{2}}),\ \  m\leq 2\ell-1 \end{cases}.
	\end{equation}
    We can now give a precise formulation of Conjecture \ref{conj 0.8'}.
    \begin{conj}\label{conj 2.1}
    	1. Assume that $1\leq \ell\leq n$. Then $E(\Phi(\theta_{\Delta(\tau,\ell)}, \theta_{\Delta(\tau,2n+\ell)},s))$ is holomorphic at $s=\ell-n$, and, up to a possible constant, for all $h\in \Sp_{8n}(\BA)$,
    	$$
    	B_{-1}(h, \theta_{\Delta(\tau,m+\ell)},f_{\Delta(\tau,\ell),\frac{\ell}{2}})	=Value_{s=\ell-n}E(\Phi(\theta_{\Delta(\tau,\ell)}, \theta_{\Delta(\tau,2n+\ell)},s))(h).
    	$$
    	2. Assume that $n<\ell\leq 2n$. Then, up to a possible constant, for all $h\in \Sp_{8n}(\BA)$,
    	$$
    	B_{-2}(h, \theta_{\Delta(\tau,m+\ell)},f_{\Delta(\tau,\ell),\frac{\ell}{2}})	=Res_{s=\ell-n}E(\Phi(\theta_{\Delta(\tau,\ell)}, \theta_{\Delta(\tau,2n+\ell)},s))(h).
    	$$
    	
    \end{conj}

	\section{ An analog of Howe duality of spherical Hecke algebras}
	
	Fix a finite place $v$ of $F$, where $\Theta_{\Delta(\tau,m+\ell),v}$ is unramified. Assume that $\tau_v=\Ind_{B_{\GL_2}(F_v)}^{\GL_2(F_v)}\chi\otimes\chi^{-1}$, where $B_{\GL_2}$ is the standard Borel subgroup of $\GL_2$, and $\chi$ is an unramified character of $F_v^*$. Recall that the central charcater of $\tau$ is trivial. Later on, we will need to assume that $\chi$ is not quadratic. This will be needed only in Lemma \ref{lem 5.5}, in order that the polynomial $P(x,y)$ (of Theorem \ref{thm 2.4}) is such that $P(q_v^{-s},q_v^s)$ has a simple zero at $s=\frac{\ell}{2}$, when $m\leq 2\ell-1$. 
	\begin{lem}\label{lem 3.1}
		$\Theta_{\Delta(\tau,m+\ell),v}$ is the unramified constituent of
		$$
		\rho_\chi=\rho^{4(m+\ell)}_\chi=\Ind_{Q_{2(m+\ell)}(F_v)}^{\Sp_{4(m+\ell)}(F_v)}\chi\circ\det.
		$$
	\end{lem}
	\begin{proof}
		This is a special case of Lemma 3.1 in \cite{G03}. We bring it for convenience. Denote, for this proof, $k=m+\ell$. The representation $\Theta_{\Delta(\tau,k),v}$ is the unramified constituent of the representation of 	$\Sp_{4k}(F_v)$ parabolically induced from the standard Borel subgroup $B_{\Sp_{4k}}(F_v)$ and the character of the diagonal subgroup given by 
		\begin{equation}\label{3.1}
			(\chi|\cdot|^{\frac{k-1}{2}+\frac{k}{2}}\otimes \chi^{-1}|\cdot|^{\frac{k-1}{2}+\frac{k}{2}})\otimes (\chi|\cdot|^{\frac{k-3}{2}+\frac{k}{2}}\otimes \chi^{-1}|\cdot|^{\frac{k-3}{2}+\frac{k}{2}})\otimes\cdots
		\end{equation}
		$$
		\cdots\otimes (\chi|\cdot|^{\frac{1-k}{2}+\frac{k}{2}}\otimes \chi^{-1}|\cdot|^{\frac{1-k}{2}+\frac{k}{2}}).
		$$
		Conjugating by an appropriate Weyl element, the representation induced from \eqref{3.1} shares the same unramified constituent with the representation parabolically induced from $B_{\Sp_{4k}}(F_v)$  and the character
		\begin{equation}\label{3.2}
			(\chi|\cdot|^{k-\frac{1}{2}}\otimes \chi|\cdot|^{-k+\frac{1}{2}})\otimes (\chi|\cdot|^{k-\frac{3}{2}}\otimes \chi|\cdot|^{-k+\frac{3}{2}})\otimes\cdots\otimes (\chi|\cdot|^{\frac{1}{2}}\otimes \chi|\cdot|^{-\frac{1}{2}}).
		\end{equation}
		We may permute the characters of $F_v^*$ in \eqref{3.2}. Hence $\Theta_{\Delta(\tau,k),v}$ is the unramified constituent of the representation of $\Sp_{4k}(F_v)$ parabolically induced from $B_{\Sp_{4k}}(F_v)$ and the character
		\begin{equation}\label{3.3}
			\chi|\cdot|^{k-\frac{1}{2}}\otimes \chi|\cdot|^{k-\frac{3}{2}}\otimes\cdots\otimes \chi|\cdot|^{\frac{1}{2}}\otimes \chi|\cdot|^{-\frac{1}{2}}\otimes\cdots\otimes \chi|\cdot|^{-k+\frac{1}{2}}.
		\end{equation}
		The character \eqref{3.3}, when viewed as a character of the diagonal subgroup of\\
		$\GL_{2k}(F_v)$, is the product of $\delta_{B_{\GL_{2k}}}^{\frac{1}{2}}$ times the restriction of $\chi\circ \det$. Since the trivial representation of $GL_{2k}(F_v)$ is a quotient of $\Ind_{B_{\GL_{2k}(F_v)}}^{GL_{2k}(F_v)}\delta_{B_{\GL_{2k}}}^{\frac{1}{2}}$, we see that $\Theta_{\Delta(\tau,k),v}$ is the unramified constituent of $\rho_\chi$.
	\end{proof}	
	
	We have the following homomorphism of spherical Hecke algebras,
	$$
	\eta^\chi_{m,\ell}: \mathcal{H}(\Sp_{4m}(F_v)// K_{4m,v})\longrightarrow \mathcal{H}(\Sp_{4\ell}(F_v)// K_{4\ell,v}),	
	$$	
	
	\begin{equation}\label{3.4}
		\eta^\chi_{m,\ell}(\xi)(g)=\int\limits_{\GL_{2(m-\ell)}(F_v)}\xi^{U_{2(m-\ell)}}(\begin{pmatrix}c\\&g\\&&c^*\end{pmatrix})\chi(\det(c))|\det(c)|^{-(m+\ell+\frac{1}{2})}dc=
	\end{equation}
	$$
	=\int\limits_{\GL_{2(m-\ell)}(F_v)}(\delta^{-\frac{1}{2}}_{Q_{2(m-\ell)}}\cdot \xi^{U_{2(m-\ell)}})(\begin{pmatrix}c\\&g\\&&c^*\end{pmatrix})\chi(\det(c))dc.
	$$
	Here, for $h\in \Sp_{4m}(F_v)$, we denote
	$$
	\xi^{U_{2(m-\ell)}}(h)=\int\limits_{U_{2(m-\ell)}(F_v)}\xi(uh)du.
	$$
	For a function $f$ in the space of $\rho_\chi$, we denote by $(1\otimes\xi)\ast f$ and $(\eta^\chi_{m,\ell}(\xi)\otimes 1)\ast f$ the functions in the space of $\rho_\chi$ given by
	$$
	(1\otimes\xi)\ast f(x)=\int\limits_{\Sp_{4m}(F_v)}\xi(h)f(x\cdot i(I_{4\ell},h))dh,
	$$
	$$
	(\eta^\chi_{m,\ell}(\xi)\otimes 1)\ast f(x)=\int\limits_{\Sp_{4\ell}(F_v)}\eta^\chi_{m,\ell}(\xi)(g)f(x\cdot i(g,I_{4m}))dg.
	$$
	In this section, we prove the following analog of (spherical) Howe duality.
	\begin{thm}\label{thm 3.1}
		For an $i(K_{4\ell,v}\times K_{4m,v})$-fixed function $f$ in the space of $\rho^{}_\chi$, and for
		$\xi\in	\mathcal{H}(\Sp_{4m}(F_v)// K_{4m,v})$, we have
		$$
		(1\otimes\xi)\ast f=(\eta^\chi_{m,\ell}(\xi)\otimes 1)\ast f.
		$$
	\end{thm}
	We start with analyzing the restriction of $\rho_\chi$ to $i(\Sp_{4\ell}(F_v)\times \Sp_{4m}(F_v))$. By Mackey theory, we need to consider $Q_{2(m+\ell)}(F_v)\backslash \Sp_{4(m+\ell)}(F_v)/i(\Sp_{4\ell}(F_v)\times \Sp_{4m}(F_v))$. By Lemma 2.2 in \cite{GS21a}, we may take the following (slightly modified) representatives. For $0\leq e\leq 2\ell$,
	\begin{equation}\label{3.5}
		\gamma_e=\gamma_e^1\gamma_e^2,
	\end{equation}
	$$
	\gamma_e^1=	\begin{pmatrix}I_e\\&I_{2(2\ell-e)}\\
		&&I_{2(e+2(m-\ell))}\\&I_{2(2\ell-e)}&&I_{2(2\ell-e)}\\&&&&I_e\end{pmatrix},
	$$
	\\
	$$	
	\gamma_e^2=	\begin{pmatrix}I_{2\ell}\\&&I_{2\ell-e}\\&I_{e+2(m-\ell)}\\&&&&I_{e+2(m-\ell)}\\&&&I_{2\ell-e}\\&&&&&I_{2\ell}\end{pmatrix}.
	$$ 
	Denote, for $g\in \Sp_{2r}(F_v)$,
	$$
	g^\iota=\begin{pmatrix}&I_k\\I_k\end{pmatrix}g\begin{pmatrix}&I_k\\I_k\end{pmatrix}.
	$$
	\begin{prop}\label{prop 3.1}
		Up to semi-simplification,\\
		\\
		$\Res_{i(\Sp_{4\ell}(F_v)\times \Sp_{4m}(F_v))}\rho_\chi\equiv$
		
		$$
		\bigoplus_{e=0}^{2\ell}ind_{Q_e^{4\ell}(F_v)\times Q_{e+2(m-\ell)}^{4m}(F_v)}^{\Sp_{4\ell}(F_v)\times \Sp_{4m}(F_v)}(((\chi\circ\det)|\det|^{m+\ell+\frac{1}{2}}\otimes (\chi\circ\det)|\det|^{m+\ell+\frac{1}{2}}))\cdot (\lambda^\iota\rho)_{\Sp_{2(2\ell-e)}},
		$$	
		where $ind$ denotes a non-normalized compact induction and $(\lambda^\iota\rho)_{\Sp_{2(2\ell-e)}}$ denotes the representation of $\Sp_{2(2\ell-e)}(F_v)\times \Sp_{2(2\ell-e)}(F_v)$ in the Schwartz space $\mathcal{S}(\Sp_{2(2\ell-e)}(F_v))$,
		$$
		(\lambda^\iota\rho)_{\Sp_{2(2\ell-e)}}(g_1,g_2)\varphi(x)=\varphi((g_1^\iota)^{-1}xg_2).
		$$
	\end{prop}
	\begin{proof}
		For $0\leq e\leq 2\ell$, let
		$$
		Q_{e,m,\ell}(F_v)=\gamma_e^{-1}Q_{2(m+\ell)}\gamma_e\cap i(\Sp_{4\ell}(F_v)\times \Sp_{4m}(F_v)).
		$$
		By Mackey theory, the semi-simplification of $\Res_{i(\Sp_{4\ell}(F_v)\times \Sp_{4m}(F_v))}\rho_\chi$ is the direct sum of the representations
		\begin{equation}\label{3.6}
			\rho_{\chi,e}=	ind_{Q_{e,m,\ell}(F_v)}^{i(\Sp_{4\ell}(F_v)\times \Sp_{4m}(F_v))}(\delta^{\frac{1}{2}}_{Q_{2(m+\ell)}}\chi\circ\det)^{\gamma_e},
		\end{equation}
		where, for $x\in Q_{e,m,\ell}(F_v)$, writing $\gamma_ex\gamma_e^{-1}=\hat{a}_xu_x$, with $a_x\in \GL_{2(m+\ell)}(F_v)$, $u_x\in U_{2(m+\ell)}(F_v)$,
		$$
		(\delta^{\frac{1}{2}}_{Q_{2(m+\ell)}}\chi\circ\det)^{\gamma_e}(x)=\delta^{\frac{1}{2}}_{Q_{2(m+\ell)}}(\gamma_ex\gamma_e^{-1})\chi(\det(a_x)).
		$$	
		Denote this character of $Q_{e,m,\ell}(F_v)$ by $\alpha_{\chi;m,\ell,e}$. The elements of $Q_{e,m,\ell}(F_v)$ have the form
		\begin{equation}\label{3.7}
			x=i(\begin{pmatrix}a&\ast&\ast\\&g&\ast\\&&a^*\end{pmatrix}, \begin{pmatrix}b&\ast&\ast&\\&g^\iota&\ast\\&&b^*\end{pmatrix})\in i(\Sp_{4\ell}(F_v)\times \Sp_{4m}(F_v)),
		\end{equation}
		where $a\in \GL_e(F_v)$, $b\in \GL_{e+2(m-\ell)}(F_v)$, $g\in \Sp_{2(2\ell-e)}(F_v)$.\\
		The character $\alpha_{\chi;m,\ell,e}$ applied to the element $x$ \eqref{3.7}, gives 
		\begin{equation}\label{3.8}
			\chi(\det(a)\det(b))|\det(a)\det(b)|^{m+\ell+\frac{1}{2}}.
		\end{equation}
		Thus, up to the identification $i$,
		\begin{equation}\label{3.9}
			\rho_{\chi,e}=	ind_{Q_e^{4\ell}(F_v)\times Q_{e+2(m-\ell)}^{4m}(F_v)}^{\Sp_{4\ell}(F_v)\times \Sp_{4m}(F_v)}(ind_{Q_{e,m,\ell}(F_v)}^{Q_e^{4\ell}(F_v)\times Q_{e+2(m-\ell)}^{4m}(F_v)}\alpha_{\chi;m,\ell,e}).
		\end{equation}
		By \eqref{3.7}, a function $\alpha$ in the space of the inner induction in \eqref{3.9} is determined by the function on $\Sp_{2(2\ell-e)}(F_v)$,
		$$
		\varphi_\alpha(h)=\alpha(I_{4\ell}, \begin{pmatrix} I_{e+2(m-\ell)}\\&h\\&&I_{e+2(m-\ell)}\end{pmatrix}).
		$$
		Using \eqref{3.8}, the left action of an element \eqref{3.7} on $\varphi_\alpha$ at $h\in \Sp_{2(2\ell-e)}(F_v)$ gives
		$$
		\alpha(\begin{pmatrix}a&\ast&\ast\\&g&\ast\\&&a^*\end{pmatrix}, \begin{pmatrix}b&\ast&\ast&\\&g^\iota h&\ast\\&&b^*\end{pmatrix})=\chi(\det(a)\det(b))|\det(a)\det(b)|^{m+\ell+\frac{1}{2}}\varphi_\alpha(h).
		$$ 	
		The right action on $\varphi_\alpha$ at $h$ of an element 
		$$
		(\begin{pmatrix}a&\ast&\ast\\&g_1&\ast\\&&a^*\end{pmatrix}, \begin{pmatrix}b&\ast&\ast&\\&g_2&\ast\\&&b^*\end{pmatrix})\in 	Q_e^{4\ell}(F_v)\times Q_{e+2(m-\ell)}^{4m}(F_v)
		$$
		gives\\
		\\
		$\alpha(\begin{pmatrix}a&\ast&\ast\\&g_1&\ast\\&&a^*\end{pmatrix}, \begin{pmatrix}b&\ast&\ast&\\&hg_2 &\ast\\&&b^*\end{pmatrix})=$
		$$
		=\chi(\det(a)\det(b))|\det(a)\det(b)|^{m+\ell+\frac{1}{2}}\varphi_\alpha((g_1^\iota)^{-1}hg_2)=
		$$
		$$
		=\chi(\det(a)\det(b))|\det(a)\det(b)|^{m+\ell+\frac{1}{2}}(\lambda^\iota\rho)_{\Sp_{2(2\ell-e)}}(g_1,g_2)\varphi_\alpha(h).
		$$ 
		The proposition follows.	
		
	\end{proof}	
	
	By the Iwasawa decomposition and the last proposition, a function $f$ in the space of $\rho_{\chi,e}$, which is fixed by $K_{4\ell,v}\times K_{4m,v}$, is determined by its restriction to $Q_e^{4\ell}(F_v)\times Q_{e+2(m-\ell)}^{4m}(F_v)$, and then, using the same notation as in the last proof,\\
	\\	
	$f(\begin{pmatrix}a&\ast&\ast\\&g_1&\ast\\&&a^*\end{pmatrix}, \begin{pmatrix}b&\ast&\ast&\\&hg_2 &\ast\\&&b^*\end{pmatrix})=$
	\begin{equation}\label{3.10}
		=\chi(\det(a)\det(b))|\det(a)\det(b)|^{m+\ell+\frac{1}{2}}\varphi_f((g_1^\iota)^{-1}hg_2),
	\end{equation}
	where $a\in \GL_e(F_v)$, $b\in \GL_{e+2(m-\ell)}(F_v)$, $g_1,g_2, h\in \Sp_{2(2\ell-e)}(F_v)$, and
	$\varphi_f\in \mathcal{H}(\Sp_{2(2\ell-e)}(F_v)// K_{2(2\ell-e),v})$. Thus, as vector spaces,
	\begin{equation}\label{3.11}	
		\rho_{\chi,e}^{K_{4\ell,v}\times K_{4m,v}}\cong \mathcal{H}(\Sp_{2(2\ell-e)}(F_v)// K_{2(2\ell-e),v}).
	\end{equation}
	Let us prove now Theorem \ref{thm 3.1}.	Let $\xi\in	\mathcal{H}(\Sp_{4m}(F_v)// K_{4m,v})$. It is enough to prove the theorem for $f$ in the space of $\rho_{\chi,e}^{K_{4\ell,v}\times K_{4m,v}}$, for each $0\leq e\leq 2\ell$. Thus, let $f$ be in the space of
	$$
	ind_{Q_e^{4\ell}(F_v)\times Q_{e+2(m-\ell)}^{4m}(F_v)}^{\Sp_{4\ell}(F_v)\times \Sp_{4m}(F_v)}(((\chi\circ\det)|\det|^{m+\ell+\frac{1}{2}}\otimes (\chi\circ\det)|\det|^{m+\ell+\frac{1}{2}}))\cdot (\lambda^\iota\rho)_{\Sp_{2(2\ell-e)}},
	$$
	and assume that $f$ is right $K_{4\ell,v}\times K_{4m,v}$ - invariant. By the Iwasawa decomposition and \eqref{3.10}, it is enough to show that for $h\in \Sp_{2(2\ell-e)}(F_v)$, the functions $(1\otimes\xi)\ast f$, $(\eta^\chi_{m,\ell}(\xi)\otimes 1)\ast f$ take the same value on $(I_{4l},\begin{pmatrix}I_{e+2(m-\ell)}\\&h\\&&I_{e+2(m-\ell)}\end{pmatrix})$. By the Iwasawa decomposition in $\Sp_{4m}(F_v)$, with respect to $Q_{e+2(m-\ell)}^{4m}(F_v)$, and \eqref{3.10}, we have\\
	\\
	$(1\otimes\xi)\ast f((I_{4l},\begin{pmatrix}I_{e+2(m-\ell)}\\&h\\&&I_{e+2(m-\ell)}\end{pmatrix}))=$
	\begin{equation}\label{3.12}
		=\int (\delta^{-\frac{1}{2}}_{Q_{e+2(m-\ell)}^{4m}}\cdot \xi^{U_{e+2(m-\ell)}})(\begin{pmatrix}b\\&x\\&&b^*\end{pmatrix})|\det(b)|^{\frac{e}{2}}\chi(\det(b))\varphi_f(hx)dbdx,
	\end{equation}
	where the integration is over $\GL_{e+2(m-\ell)}(F_v)\times \Sp_{2(2\ell-e)}(F_v)$. Using the Iwasawa decomposition in $\GL_{e+2(m-\ell)}(F_v)$ with respect to the Borel subgroup, the integral \eqref{3.12} becomes
	\begin{equation}\label{3.13}
		\int (\delta^{-\frac{1}{2}}_{Q_{1^{e+2(m-\ell)}}^{4m}}\cdot \xi^{U_{e+2(m-\ell)}})(\begin{pmatrix}t\\&x\\&&t^*\end{pmatrix})\delta_{B_{\GL_{e+2(m-\ell)}}}^{-\frac{1}{2}}(t)|\det(t)|^{\frac{e}{2}}\chi(\det(t))\varphi_f(hx)dtdx,
	\end{equation}
	where $t$ is integrated over $T_{\GL_{e+2(m-\ell)}}(F_v)$, the diagonal subgroup of $\GL_{e+2(m-\ell)}(F_v)$, and $x$ is integrated over $\Sp_{2(2\ell-e)}(F_v)$. In the same way we have,\\
	\\
	$(\eta^\chi_{m,\ell}(\xi)\otimes 1)\ast f((I_{4l},\begin{pmatrix}I_{e+2(m-\ell)}\\&h\\&&I_{e+2(m-\ell)}\end{pmatrix}))=$
	\begin{equation}\label{3.14}
		=\int (\delta_{Q_e^{4\ell}}^{-\frac{1}{2}}\cdot (\eta^\chi_{m,\ell}(\xi))^{U_e^{4\ell}})(\begin{pmatrix}a\\&x\\&&a^*\end{pmatrix})|\det(a)|^{m-\ell+\frac{e}{2}}\chi(\det(a))\varphi_f((x^\iota)^{-1}h)dadx,
	\end{equation}
	where the integration is over $\GL_e(F_v)\times \Sp_{2(2\ell-e)}(F_v)$. It is easy to see that for $\varphi\in 
	\mathcal{H}(\Sp_{2r}(F_v)\backslash\backslash K_{2r,v})$, we have, for all $g\in \Sp_{2r}(F_v)$, 
	$$
	\varphi(g)=\varphi(g^{-1})=\varphi(g^\iota).
	$$
	Then the integral \eqref{3.14} becomes 
	$$
	\int (\delta_{Q_e^{4\ell}}^{-\frac{1}{2}}\cdot (\eta^\chi_{m,\ell}(\xi))^{U_e^{4\ell}})(\begin{pmatrix}a\\&x\\&&a^*\end{pmatrix})|\det(a)|^{m-\ell+\frac{e}{2}}\chi(\det(a))\varphi_f(hx)dadx=
	$$	
	\begin{equation}\label{3.15}
		\int (\delta_{Q_{2(m-\ell),e}^{4m}}^{-\frac{1}{2}}\cdot \xi^{U^{4m}_{2(m-\ell),e}})(\begin{pmatrix}c\\&a\\&&x\\&&&a^*\\&&&&c^*\end{pmatrix})|\det(a)|^{m-\ell+\frac{e}{2}}\chi(\det(c)\det(a))
	\end{equation}		
	$$
	\varphi_f(hx)dcdadx,
	$$	
	where $c$, $a$, $x$ are integrated over $\GL_{2(m-\ell)}(F_v)$,$\GL_e(F_v)$, $\Sp_{2(2\ell-e)}(F_v)$, respectively. Using the Iwasawa decomposition in $\GL_{2(m-\ell)}(F_v)$,$\GL_e(F_v)$ with respect to the Borel subgroups, the integral \eqref{3.15} becomes exactly the integral \eqref{3.13}. This completes the proof of Theorem \ref{thm 3.1}.

	Recall the Satake isomorphism  
	\begin{equation}\label{3.16}
		\mathcal{H}(\Sp_{4m}(F_v)// K_{4m,v})\cong \BC[Z_1^{\pm 1},...,Z_{2m}^{\pm 1}]^{W_{\Sp_{4m}}},
	\end{equation}
	where $W_{\Sp_{4m}}$ is the Weyl group of $\Sp_{4m}$. It is given as follows.
	Let, for $t\in T_{\Sp_{4m}}(F_v)$,
	\begin{equation}\label{3.17}
		\mathcal{S}(\xi)(t)=\delta_{B_{\Sp_{2m}}}^{-\frac{1}{2}}(t)\xi^{U_{1^{2m}}}(t)=\delta_{B_{\Sp_{2m}}}^{-\frac{1}{2}}(t)\int\limits_{U_{1^{2m}}(F_v)}\xi(ut)du.
	\end{equation}	
	This is the Satake transform of $\xi$. It defines an isomorphism  
	$$
	\mathcal{H}(\Sp_{4m}(F_v)// K_{4m,v})\cong \mathcal{H}(T_{\Sp_{4m}}(F_v)\backslash T_{\Sp_{4m}}(\mathcal{O}_v))^{W_{\Sp_{4m}}}.
	$$ 
	Let $\mu$ be an unramified character of $T_{\Sp_{4m}}(F_v)$, 
	$$
	\mu(\begin{pmatrix}t\\&t^*\end{pmatrix})=\mu_1(t_1)\cdot \ldots\cdot \mu_{2m}(t_{2m}), \ t=\begin{pmatrix}t_1\\&\ddots\\&&t_{2m}\end{pmatrix},
	$$
	where $\mu_1,...,\mu_{2m}$ are unramified characters of $F_v^*$. Let $f_\mu$ be a spherical vector in the
	one dimensional subspace $(\Ind_{B_{\Sp_{4m}}(F_v)}^{\Sp_{4m}(F_v)}\mu)^{K_{4m,v}}$ . Then 
	$$
	\xi\ast f_\mu=(\mathcal{S}(\xi),\mu)f_\mu,
	$$	
	where
	$$
	(\mathcal{S}(\xi),\mu)=\int\limits_{T_{\Sp_{4m}}(F_v)}\mathcal{S}(\xi)(t)\mu(t)dt.
	$$
	There is a unique element $\tilde{\mathcal{S}}(\xi)\in \BC[Z_1^{\pm 1},...,Z_{2m}^{\pm 1}]^{W_{\Sp_{4m}}}$, such that
	$$
	(\mathcal{S}(\xi),\mu)=	\tilde{\mathcal{S}}(\xi)(\mu_1(p_v)^{\pm 1},...,\mu_{2m}(p_v)^{\pm 1}),
	$$
	where $p_v$ is a generator of the maximal ideal of $\mathcal{O}_v$. The isomorphism \eqref{3.16} is given by $\xi\mapsto \hat{\mathcal{S}}(\xi)$. It will be convenient to denote
	$$
	\tilde{\mathcal{S}}(\xi)(Z_1^{\pm 1},...,Z_{2m}^{\pm 1})=\hat{\mathcal{S}}(\xi)(Z_1,...,Z_{2m}).
	$$
	Note that\\
	\\
	$\hat{\mathcal{S}}(\eta^\chi_{m,\ell}(\xi))(Z_1,...,Z_{2\ell})=$
	\begin{equation}\label{3.18}
		=\hat{\mathcal{S}}(\xi)(\chi(p_v)q_v^{-(m-\ell-\frac{1}{2})},\chi(p_v)q_v^{-(m-\ell-\frac{3}{2})},...,\chi(p_v)q_v^{m-\ell-\frac{1}{2}},Z_1,...,Z_{2\ell}).
	\end{equation}
	We define the function $\xi_{m,\ell,v}\in \mathcal{H}(\Sp_{4m}(F_v)// K_{4m,v})$ to be such that
	\begin{equation}\label{3.19}
		\hat{\mathcal{S}}(\xi_{m,\ell,v})(Z_1,...,Z_{2m})=\prod_{i=1}^{2m}(Z_i-\chi(p_v)q_v^{-m+\ell-\frac{1}{2}})(Z_i^{-1}-\chi(p_v)q_v^{-m+\ell-\frac{1}{2}}).	
	\end{equation}	
	By \eqref{3.18},
	\begin{equation}\label{3.20}
		\hat{\mathcal{S}}(\eta^\chi_{m,\ell}(\xi_{m,\ell,v}))(Z_1,...,Z_{2\ell})=\alpha_{m,\ell,v}\prod_{i=1}^{2\ell}(Z_i-\chi(p_v)q_v^{-m+\ell-\frac{1}{2}})(Z_i^{-1}-\chi(p_v)q_v^{-m+\ell-\frac{1}{2}}),
	\end{equation}
	where
	$$
	\alpha_{m,\ell,v}=\prod_{j=1}^{2(m-\ell)}[(q_v^{-m+\ell+j-\frac{1}{2}}-q_v^{-m+\ell-\frac{1}{2}})(q_v^{m-\ell-j+\frac{1}{2}}-\chi^2(p_v)q_v^{-m+\ell-\frac{1}{2}})]\neq 0.
	$$	
	When $\ell=m$, we define $\alpha_{m,m,v}=1$. The reason that $\alpha_{m,\ell,v}$ is nonzero also for $\ell<m$ is that, for $1\leq j\leq 2(m-\ell)$, $q_v^{m-\ell-j+\frac{1}{2}}-\chi^2(p_v)q_v^{-m+\ell-\frac{1}{2}}\neq 0$. Otherwise, writing $|\chi(p_v)|=q_v^\alpha$, we get that $2\alpha=2(m-\ell)-j+1$, and hence $\frac{1}{2}\leq \alpha\leq m-\ell$. Recall that $\chi$ is obtained by considering our cuspidal representation $\tau$ at the finite place $v$, where $\tau_v$ is unramified and induced from the character $\chi\otimes \chi^{-1}$ of $B_{\GL_2}(F_v)$. In particular, we know that $-\frac{1}{2}<\alpha <\frac{1}{2}$, and so we get a contradiction.

	\section{The Fourier expansion of $((1\otimes\xi_{m,\ell,v})\ast \theta_{\Delta(\tau,m+\ell)})$ along $i(U_{2\ell}^{4\ell}\times I_{4m})$}
	
	We go back to the notation of Theorem \ref{thm 2.1}. In this section we write down the Fourier expansion of $((1\otimes\xi_{m,\ell,v})\ast \theta_{\Delta(\tau,m+\ell)})$ along $i(U_{2\ell}^{4\ell}\times I_{4m})$. Note 
	that the elements of the last subgroup have the form
	\begin{equation}\label{4.1}
		i(u^{4\ell}_{2\ell}(x),I_{4m})=\begin{pmatrix}I_{2\ell}&0&x\\&I_{4m}&0\\&&I_{2\ell}\end{pmatrix},\ \ u^{4\ell}_{2\ell}(x)=\begin{pmatrix}I_{2\ell}&x\\&I_{2\ell}\end{pmatrix},\ {}^t(w_{2\ell}x)=w_{2\ell}x.
	\end{equation}
	Denote by $S_{2\ell}(F)$ the subspace of matrices $x\in M_{2\ell}(F)$, such that $w_{2\ell}x$ is symmetric. Fix a nontrivial character $\psi$ of $F\backslash \BA$. We start with the Fourier expansion of $\theta_{\Delta(\tau,m+\ell)}$, along $i(U_{2\ell}^{4\ell}\times I_{4m})$, viewed, first, as a function of $b\in \Sp_{4(m+\ell)}(\BA)$,
	\begin{equation}\label{4.2}
		\theta_{\Delta(\tau,m+\ell)}(b)=\sum_{A\in S_{2\ell}(F)}\theta_{\Delta(\tau,m+\ell)}^{\psi_A}(b),
	\end{equation}	
	where
	$$
	\theta_{\Delta(\tau,m+\ell)}^{\psi_A}(b)=\int\limits_{S_{2\ell}(F)\backslash S_{2\ell}(\BA)}\theta_{\Delta(\tau,m+\ell)}((u_{2\ell}^{4\ell}(x),I_{4m})b)\psi^{-1}(tr(Ax))dx.
	$$	
	We denote by $\psi_A$ the character of $U^{4\ell}_{2\ell}(\BA)$ given by $\psi_A(u^{4\ell}_{2\ell}(x))=\psi(tr(Ax))$. Consider the sum of the Fourier coefficients $\theta_{\Delta(\tau,m+\ell)}^{\psi_A}$ over all $A$ with rank $c$, $0\leq c\leq 2\ell$. Consider the action of $\GL_{2\ell}(F)$, $\gamma\cdot A=\gamma^*A\gamma^{-1}$, $\gamma\in \GL_{2\ell}(F)$, $\gamma^*=w_{2\ell}{}^t\gamma^{-1}w_{2\ell}$.  Then there is a diagonal matrix $\delta'=diag(\delta_c,...,\delta_1)$, $\delta_i\in F^*$, such that $\diag(\delta',0,...,0)w_{2\ell}=\gamma^*A\gamma^{-1}$ is in the orbit of $A$. In this case,\\
	\\
	$\theta_{\Delta(\tau,m+\ell)}^{\psi_A}(b)=$
	\begin{equation}\label{4.3}
		=\int\limits_{S_{2\ell}(F)\backslash S_{2\ell}(\BA)}\theta_{\Delta(\tau,m+\ell)}((u_{2\ell}^{4\ell}(x),I_{4m})b)\psi^{-1}(tr(\begin{pmatrix}\delta'\\&0\end{pmatrix}w_{2\ell}\gamma x(\gamma^*)^{-1}))dx=
	\end{equation}
	$$
	=\int\limits_{S_{2\ell}(F)\backslash S_{2\ell}(\BA)}\theta_{\Delta(\tau,m+\ell)}((u_{2\ell}^{4\ell}(x)\hat{\gamma},I_{4m})b)\psi^{-1}(tr(\begin{pmatrix}&w_c\delta\\0\end{pmatrix}x))dx,
	$$	
	where $\delta=w_c\delta'w_c=diag(\delta_1,...,\delta_c)$. The stabilizer of $\diag(\delta',0,...,0)w_{2\ell}$ inside $\GL_{2\ell}(F)$ is
	\begin{equation}\label{4.4}
		P^\delta_{2\ell-c,c}(F)=\{\begin{pmatrix}\gamma_1&\gamma_2\\0&\gamma_4\end{pmatrix}\in \GL_{2\ell}(F)\ |\ \gamma_1\in \GL_{2\ell-c}(F),\ \gamma_4\in \RO_{c,\delta}(F) \},
	\end{equation}
	where $\RO_{c,\delta}$ denotes the $F$-orthogonal group in $c$ variables corresponding to $\delta^{w_c}=w_c\delta w_c$. Using \eqref{4.3}, the expansion \eqref{4.2} becomes
	\begin{equation}\label{4.5}
		\theta_{\Delta(\tau,m+\ell)}(b)=\sum_{c=0}^{2\ell}\sum_{[\delta]\in [T_c(F)]}\sum_{\gamma\in P^\delta_{2\ell-c,c}(F)\backslash \GL_{2\ell}(F)}\theta^{\psi_\delta}_{\Delta(\tau,m+\ell)}((\hat{\gamma},I_{4m})b),
	\end{equation}
	where\\
	\\
	$\theta^{\psi_\delta}_{\Delta(\tau,m+\ell)}(b)=$
	\begin{equation}\label{4.6}	
		=\int\limits_{S_{2\ell}(F)\backslash S_{2\ell}(\BA)}\theta_{\Delta(\tau,m+\ell)}(\begin{pmatrix}I_{2\ell-c}&0&0&x_1&x_2\\&I_c&0&x_3&x'_1\\&&I_{4m}&0&0\\&&&I_c&0\\&&&&I_{2\ell-c}\end{pmatrix}b)\psi^{-1}(tr(w_c\delta x_3))dx.
	\end{equation}
	In \eqref{4.5}, $[\delta]$ varies over the equivalence classes of $c\times c$ diagonal, invertible matrices (over $F$), representing quadratic forms $\delta_1x_1^2+\cdots+\delta_cx_c^2$. We re-denote, for short, by $\psi_\delta$, the character of $U_{2\ell}^{4\ell}(\BA)$, which we denoted before by $\psi_{diag(w_c\delta w_c,0,...,0)w_{2\ell}}$. For fixed $b$, consider the smooth function on the compact abelian group\\ $M_{(2\ell-c)\times 4m}(F)\backslash M_{(2\ell-c)\times 4m}(\BA)$,\\
	\\
	$\theta^{\psi_\delta}_{\Delta(\tau,m+\ell)}(b)(y)=$
	$$
	=\int\limits_{S_{2\ell}(F)\backslash S_{2\ell}(\BA)}\theta_{\Delta(\tau,m+\ell)}(\begin{pmatrix}I_{2\ell-c}&0&y&x_1&x_2\\&I_c&0&x_3&x'_1\\&&I_{4m}&0&y'\\&&&I_c&0\\&&&&I_{2\ell-c}\end{pmatrix}b)\psi^{-1}(tr(w_c\delta x_3))dx,
	$$	
	and write its Fourier expansion at $y=0$. Then
	\begin{equation}\label{4.7}
		\theta^{\psi_\delta}_{\Delta(\tau,m+\ell)}(b)=\sum_{B\in M_{4m\times (2\ell-c)}(F)}\theta^{\psi_{\delta,B}}_{\Delta(\tau,m+\ell)}(b),
	\end{equation}
	where
	\begin{equation}\label{4.8}
		\theta^{\psi_{\delta,B}}_{\Delta(\tau,m+\ell)}(b)=\int\limits_{M_{(2\ell-c)\times 4m}(F)\backslash M_{(2\ell-c)\times 4m}(\BA)}\theta^{\psi_\delta}_{\Delta(\tau,m+\ell)}(b)(y)\psi^{-1}(tr(yB))dy.
	\end{equation}
	\begin{prop}\label{prop 4.1}
		Assume that $0\leq c<2\ell$ and that $\theta^{\psi_{\delta,B}}_{\Delta(\tau,m+\ell)}$ is nontrivial. Then the column space of $B$ is a totally isotropic subspace of $F^{4m}$ with respect to the synplectic form corresponding to $J_{4m}$.
	\end{prop}
	\begin{proof}
		Consider the smooth function on the compact abelian group\\
		$M_{(2\ell-c)\times c}(F)\backslash M_{(2\ell-c)\times c}(\BA)$,
		$$
		\theta^{\psi_{\delta,B}}_{\Delta(\tau,m+\ell)}(b)(z)=\theta^{\psi_{\delta,B}}_{\Delta(\tau,m+\ell)}((\widehat{\begin{pmatrix}I_{2\ell-c}&z\\&I_c\end{pmatrix}},I_{4m})b).
		$$	
		We write its Fourier expansion at $z=0$,
		\begin{equation}\label{4.9}
			\theta^{\psi_{\delta,B}}_{\Delta(\tau,m+\ell)}(b)=\sum_{D\in M_{c\times (2\ell-c)}(F)}\theta^{\psi_{\delta,B,D}}_{\Delta(\tau,m+\ell)}(b),
		\end{equation}	
		where
		\begin{equation}\label{4.10}
			\theta^{\psi_{\delta,B,D}}_{\Delta(\tau,m+\ell)}(b)=\int\limits_{M_{(2\ell-c)\times c}(F)\backslash M_{(2\ell-c)\times c}(\BA)}\theta^{\psi_{\delta,B}}_{\Delta(\tau,m+\ell)}(b)(z)\psi^{-1}(tr(zD))dz.
		\end{equation}
		By our assumption, for the given $B$, one of the Fourier coefficients \eqref{4.10} is nontrivial. Then there is $D\in M_{c\times (2\ell-c)}(F)$, such that the following Fourier coefficient, being an inner integral of \eqref{4.10}, is nontrivial, for some automorphic form $\theta'_{\Delta(\tau,m+\ell)}$ in the space of $\Theta_{\Delta(\tau,m+\ell)}$, 
		\begin{equation}\label{4.11}
			\int\theta'_{\Delta(\tau,m+\ell)}(\begin{pmatrix}I_{2\ell-c}&e&t\\&I_{4m+2c}&e'\\&&I_{2\ell-c}\end{pmatrix})\psi^{-1}(tr(e\begin{pmatrix}D\\B\\0\end{pmatrix}))d(e,t)\neq 0,
		\end{equation}
		where the integration is over $U^{4(m+\ell)}_{2\ell-c}(F)\backslash U^{4(m+\ell)}_{2\ell-c}(\BA)$. Consider the right action of $\GL_{2\ell-c}(F)\times \Sp_{4m+2c}(F)$ on $M_{(4m+2c)\times (2\ell-c)}(F)$, given by $L\cdot (\alpha,\beta)=\beta^{-1}L\alpha$. See Lemma 9.1 in \cite{GS21a} for representatives of this action. Exactly as in the proof of Prop. 9.2 in \cite{GS21a}, using Prop. \ref{prop 1.3}, we see that the $\GL_{2\ell-c}(F)\times \Sp_{4m+2c}(F)$-orbit of $\begin{pmatrix}D\\B\\0\end{pmatrix}$ contains a matrix of the form 
		$$
		\begin{pmatrix}I_{\ell_1}&0\\0&0\end{pmatrix},\ \ \ell_1\leq 2\ell-c.
		$$
		Thus, the column space of $\begin{pmatrix}D\\B\\0\end{pmatrix}$ is a totally isotropic subspace of $F^{4m+2c}$ of dimension $\ell_1$. This implies that the column space of $B$ is a totally isotropic subspace of $F^{4m}$ of dimension $rank(B)\leq \ell_1$.
	\end{proof}
	\begin{thm}\label{thm 4.1}
		Assume that $0\leq c<2\ell$, that the column space of $B$ is totally isotropic and $rank(B)<2\ell-c$. Then, for all $(g,h)\in \Sp_{4\ell}(\BA)\times \Sp_{4m}(\BA)$,
		$$
		(1\otimes\xi_{m,\ell,v})\ast\theta^{\psi_{\delta,B}}_{\Delta(\tau,m+\ell)}(g,h)=0.
		$$
	\end{thm}
	
	\begin{proof}
		We start the proof with arbitrary $\theta_{\Delta(\tau,m+\ell)}$, and later apply the convolution by $1\otimes\xi_{m,\ell,v}$.  Assume that $rank(B)=k<2\ell-c$. By the last proposition, we may write
		$$
		B=\beta^{-1}\begin{pmatrix}0&I_k\\0&0\end{pmatrix}_{4m\times (2\ell-c)}\cdot \alpha,
		$$
		where $(\alpha,\beta)\in P_{2\ell-c-k,k}(F)\times_d Q^{4m}_k(F)\backslash GL_{2\ell-c}(F)\times \Sp_{4m}(F)$, and\\ 
		$P_{2\ell-c-k,k}(F)\times_d Q^{4m}_k(F)$ denotes the subgroup of $P_{2\ell-c-k,k}(F)\times Q^{4m}_k(F)$, consisting of the elements
		$$
		(\begin{pmatrix}\alpha_1&\ast\\0&\alpha_2\end{pmatrix},\begin{pmatrix}\alpha_2&\ast&\ast\\&\beta_2&\ast\\&&\alpha_2^*\end{pmatrix}),\ \ \alpha_1\in \GL_{2\ell-c-k}(F), \alpha_2\in \GL_k(F), \beta_2\in \Sp_{4m-2k}(F).
		$$ 	
		Then, for $b\in \Sp_{4(m+\ell)}(\BA)$,\\
		\\
		$\theta^{\psi_{\delta,B}}_{\Delta(\tau,m+\ell)}(b)=$	
		\begin{equation}\label{4.12}	
			=\int \theta_{\Delta(\tau,m+\ell)}(\begin{pmatrix}I_{2\ell-c}&0&y&x_1&x_2\\&I_c&0&x_3&x'_1\\&&I_{4m}&0&y'\\&&&I_c&0\\&&&&I_{2\ell-c}\end{pmatrix}((\hat{\alpha},\beta)b)
		\end{equation}
		$$
		\psi^{-1}(tr(w_c\delta x_3)+tr(\begin{pmatrix}0&I_k\\0&0\end{pmatrix}y))dydx:=\theta^{\psi_{\delta,k}}_{\Delta(\tau,m+\ell)}((\hat{\alpha},\beta)b).
		$$
		Next, we write the Fourier expansion of the function on $M_{k\times c}(F)\backslash M_{k\times c}(\BA)$, 
		$$
		t\mapsto \theta^{\psi_{\delta,k}}_{\Delta(\tau,m+\ell)}(\widehat{\begin{pmatrix}I_{2\ell-c-k}\\&I_k&t\\&&I_c\end{pmatrix}}b),
		$$
		at $t=0$, \\
		\\
		$\theta^{\psi_{\delta,k}}_{\Delta(\tau,m+\ell)}(b)=$
		\begin{equation}\label{4.13}
			=\sum_{E\in M_{c\times k}(F)}\int_{M_{k\times c}(F)\backslash M_{k\times c}(\BA)}\theta^{\psi_{\delta,k}}_{\Delta(\tau,m+\ell)}(\widehat{\begin{pmatrix}I_{2\ell-c-k}\\&I_k&t\\&&I_c\end{pmatrix}}b)\psi^{-1}(tr(tE))dt.
		\end{equation}
		The summand in \eqref{4.13}, corresponding to $E$, is, using \eqref{4.12}, equal to
		\begin{equation}\label{4.14}
			\int \theta_{\Delta(\tau,m+\ell)}(\begin{pmatrix}I_{2\ell-c-k}&0&0&y_1&y_2&y_3&x_1&x_2&x_3\\&I_k&t&y_4&y_5&y_6&x_4&x_5&x'_2\\&&I_c&0&0&0&x_6&x_4'&x_1'\\&&&I_k&0&0&0&y'_6&y'_3\\&&&&I_{4m-2k}&0&0&y'_3&y'_2\\&&&&&I_k&0&y'_4&y'_1\\&&&&&&I_c&t'&0\\&&&&&&&I_k&0\\&&&&&&&&I_{2\ell-c-k}\end{pmatrix}\hat{v_1}(E)b)
		\end{equation}
		
		$$ 
		\psi^{-1}(tr(w_c\delta x_6)+tr(y_4))d(x,y,t):=\theta^{\psi_{\delta,k,0}}_{\Delta(\tau,m+\ell)}(\hat{v_1}(E)b),
		$$
		where $v_1(E)=\begin{pmatrix}I_{2\ell-c}\\&I_c&-E\\&&I_k\end{pmatrix}$. From Prop. \ref{prop 1.3} it follows that  $\theta^{\psi_{\delta,k,0}}_{\Delta(\tau,m+\ell)}(b)$ is left invariant to $\hat{v}_2(u)$, where $v_2(u)=\begin{pmatrix}I_{2\ell-c-k}&0&u\\&I_k&0\\&&I_c\end{pmatrix}$, $u\in M_{(2\ell-c-k)\times c}(\BA)$.
		The reason is that in the Fourier expansion of $\theta^{\psi_{\delta,k,0}}_{\Delta(\tau,m+\ell)}$ along the subgroup of the elements $\hat{v}_2(u)$, $u\in M_{(2\ell-c-k)\times c}(F)\backslash M_{(2\ell-c-k)\times c}(\BA)$, only the trivial character contributes. The Fourier coefficients with respect to nontrivial characters of $M_{(2\ell-c-k)\times c}(F)\backslash M_{(2\ell-c-k)\times c}(\BA)$ give rise to Fourier coefficients on $\Theta_{\Delta(\tau,m+\ell)}$ which correspond to a symplectic partition of $4(m+\ell)$ which is strictly larger than $(2^{2(m+\ell)})$. We omit the details here. We have used this argument many times before.  See, for example, Theorems 8.2, 9.5 in \cite{GS21a}, where we carry out similar proofs in full detail. We conclude that\\
		\\
		$\theta^{\psi_{\delta,k,0}}_{\Delta(\tau,m+\ell)}(b)=$
		\begin{equation}\label{4.15}
			\int \theta_{\Delta(\tau,m+\ell)}(\begin{pmatrix}I_{2\ell-c-k}&0&u&y_1&y_2&y_3&x_1&x_2&x_3\\&I_k&t&y_4&y_5&y_6&x_4&x_5&x'_2\\&&I_c&0&0&0&x_6&x_4'&x_1'\\&&&I_k&0&0&0&y'_6&y'_3\\&&&&I_{4m-2k}&0&0&y'_3&y'_2\\&&&&&I_k&0&y'_4&y'_1\\&&&&&&I_c&t'&u'\\&&&&&&&I_k&0\\&&&&&&&&I_{2\ell-c-k}\end{pmatrix}b)
		\end{equation}
		
		$$ 
		\psi^{-1}(tr(w_c\delta x_6)+tr(y_4))d(x,y,t,u).
		$$
		(We take the measure of $F\backslash \BA$ to be 1.) Apply a conjugation inside $\theta_{\Delta(\tau,m+\ell)}$, in \eqref{4.15}, by $\hat{\omega}=\widehat{\begin{pmatrix}I_{2\ell-c}\\&&I_k\\&I_c\end{pmatrix}}$. Then \eqref{4.15} becomes\\
		\\
		$\theta^{\psi_{\delta,k,0}}_{\Delta(\tau,m+\ell)}(b)=$
		\begin{equation}\label{4.16}
			\int \theta_{\Delta(\tau,m+\ell)}(\begin{pmatrix}I_{2\ell-c}&y&x\\&I_{4m+2c}&y'\\&&I_{2\ell-c}\end{pmatrix}\begin{pmatrix}I_{2\ell-c+k}\\&I_c&0&z\\&&I_{4m-2k}&0\\&&&I_c\\&&&&I_{2\ell-c+k}\end{pmatrix}\hat{\omega}b)
		\end{equation}
		$$
		\psi^{-1}(tr(w_c\delta z)+tr(y\begin{pmatrix}0&I_k\\0&0\end{pmatrix}_{(4m+2c)\times (2\ell-c)})d(x,y,z).
		$$
		Repeating the last argument, several more times, using Prop. \ref{prop 1.3}, it follows that the function
		$$
		b\mapsto \int \theta_{\Delta(\tau,m+\ell)}(\begin{pmatrix}I_{2\ell-c}&y&x\\&I_{4m+2c}&y'\\&&I_{2\ell-c}\end{pmatrix}b)
		\psi^{-1}(tr(y\begin{pmatrix}0&I_k\\0&0\end{pmatrix})d(x,y)
		$$
		is left invariant to the following subgroups. First, it is invariant to 
		$$
		diag(I_{2\ell-c},\begin{pmatrix}I_k&0&u_2\\&I_{4m-2k+2c}&0\\&&I_k\end{pmatrix},I_{2\ell-c}),
		\ \ u_2\in S_k(\BA).
		$$
		Then it is invariant to $(diag(I_{2\ell-c}, U^{4m+2c}_k(\BA),I_{2\ell-c})$, and, finally, it is also invariant to $\widehat{\begin{pmatrix}I_{2\ell-c-k}&v\\&I_k\end{pmatrix}}$, $v\in M_{(2\ell-c-k)\times k}(\BA)$.
		Again, at each time, it follows that in the Fourier expansion along the corresponding subgroup, only the trivial character contributes. Going back to \eqref{4.12}, we conclude that\\
		\\
		$\theta^{\psi_{\delta,B}}_{\Delta(\tau,m+\ell)}(b)=$
		\begin{equation}\label{4.17}
			\sum_{E\in M_{c\times k}(F)}\int \theta^{U^{4(m+\ell)}_{2\ell-c+k}}_{\Delta(\tau,m+\ell)}(\begin{pmatrix}I_{2\ell-c-k}&a_1&a_2\\&I_k&a_3\\&&I_k\\&&&I_c&0&z\\&&&&I_{4m-2k}&0\\&&&&&I_c\\&&&&&&I_k&a'_3&a'_2\\&&&&&&&I_k&a'_1\\&&&&&&&&I_{2\ell-c-k}\end{pmatrix}\cdot
		\end{equation}
		$$
		\widehat{\omega v_1(E)}(\hat{\alpha}g,\beta h))\psi^{-1}(tr(w_c\delta z)+tr(a_3))d(a_1,a_2,a_3,z).
		$$
		All integrations are over variables in $F\backslash \BA$. Note that since $\tau$ is cuspidal on $\GL_2(\BA)$, for the constant term $\theta^{U^{4(m+\ell)}_{2\ell-c+k}}_{\Delta(\tau,m+\ell)}$ to be nontrivial, $2\ell-c+k$ must be even, and hence $c+k$ must be even. 
		
		By Prop. \ref{prop 1.2}, the constant term $\theta_{\Delta(\tau,m+\ell)}\mapsto \theta^{U^{4(m+\ell)}_{2\ell-c+k}}_{\Delta(\tau,m+\ell)}$ projects $\Theta_{\Delta(\tau,m+\ell)}$ into 
		$$
		\Ind_{Q_{2\ell-c+k}(\BA)}^{\Sp_{4(m+\ell)}(\BA)}(\Delta(\tau,\ell+\frac{k-c}{2})||det\cdot|^{-m-\frac{\ell}{2}+\frac{k-c}{4}}\otimes\Theta_{\Delta(\tau,m+\frac{c-k}{2})}).
		$$
		In \eqref{4.17}, we further take the constant term of $\Delta(\tau,\ell+\frac{k-c}{2})$ along $V_{2\ell-c-k,2k}$, followed by the Fourier coefficient with respect to the subgroup of elements\\ $\begin{pmatrix}I_{2\ell-c-k}\\&I_k&a_3\\&&I_k\end{pmatrix}$ and the character $\psi(tr(a_3))$. This projects $\Delta(\tau,\ell+\frac{k-c}{2})$ into
		$$
		\Ind_{P_{2\ell-c-k,2k}(\BA)}^{\GL_{2\ell-c+k}(\BA)}(\Delta(\tau,\ell-\frac{c+k}{2})|\det\cdot|^{-\frac{k}{2}}\otimes \Delta_\psi(\tau,k)|\det\cdot|^{\frac{2\ell-c-k}{4}}),
		$$
		where $\Delta_\psi(\tau,k)$ denotes the representation by right translations of $\GL_{2k}(\BA)$ in the space of functions, 
		$$
		x\mapsto \int_{M_k(F)\backslash M_k(\BA)}\varphi(\begin{pmatrix}I_k&a\\&I_k\end{pmatrix}x)\psi^{-1}(tr(a))da, \ \ \varphi\in \Delta(\tau,k).
		$$
		Denote by $\Theta^{\psi_\delta}_{\Delta(\tau,m+\frac{c-k}{2})}$ the representation by right translations of $\Sp_{4m+2c-2k}(\BA)$ in the space of functions
		$$
		y\mapsto \int_{S_c(F)\backslash S_c(\BA)}\theta_{\Delta(\tau,m+\frac{c-k}{2})}(\begin{pmatrix}I_c&0&z\\&I_{4m-2k}&0\\&&I_c\end{pmatrix}y)\psi^{-1}(tr(w_c\delta z))dz.
		$$
		Then the integral in \eqref{4.17} projects $\Theta_{\Delta(\tau,m+\ell)}$ into 
		\begin{equation}\label{4.18}
			\Ind_{Q_{2\ell-c-k,2k}(\BA)}^{\Sp_{4(m+\ell)}(\BA)}(\Delta(\tau,\ell-\frac{c+k}{2})||det\cdot|^{-m-\frac{\ell}{2}-\frac{c+k}{4}}\otimes\Delta_\psi(\tau,k)|\det\cdot|^{-m-\frac{c}{2}}\otimes\Theta^{\psi_\delta}_{\Delta(\tau,m+\frac{c-k}{2})}).	 
		\end{equation}
		Let $f$ be a function in the space of \eqref{4.18}. Denote, for $g\in \Sp_{4\ell}(\BA)$, and fixed $E_0\in M_{c\times k}(F)$, $h_0\in \Sp_{4m}(\BA)$,
		\begin{equation}\label{4.19}
			f_{E_0,h_0}(g)=f(\widehat{\omega v_1(E_0)}i(g,h_0)).
		\end{equation}
		Then by a simple check, we find that for $u\in U^{4\ell}_{2\ell-c-k}(\BA)$ and $a\in \GL_{2\ell-c-k}(\BA)$,
		\begin{equation}\label{4.20}
			f_{E_0,h_0}(u\begin{pmatrix}a\\&I_{2(c+k)}\\&&a^*\end{pmatrix}g)=|\det(a)|^{m+\frac{\ell+1}{2}+\frac{c+k}{4}}\Delta(\tau,\ell-\frac{c+k}{2})(a)f_{E_0,h_0}(g).
		\end{equation}	
		Let $f_v$ be a right $K_{4\ell,v}\times K_{4m,v}$-invariant function in the space of the $v$-component of \eqref{4.18}. For a fixed $h_0\in \Sp_{4m}(F_v)$, consider the function on $\Sp_{4\ell}(F_v)$, $f_{E_0,h_0,v}$ which is the local analog at $v$ of \eqref{4.19}, and convolve it with $\eta_{m,\ell}^\chi(\xi_{m,\ell,v})$:
		$$
		\eta_{m,\ell}^\chi(\xi_{m,\ell,v})\ast f_{E_0,h_0,v} (g_0)=\int\limits_{\Sp_{4\ell}(F_v)} \eta_{m,\ell}^\chi(\xi_{m,\ell,v})(g)f_{E_0,h_0,v}(g_0g)dg.
		$$
		We will show that 
		\begin{equation}\label{4.21}
			\eta_{m,\ell}^\chi(\xi_{m,\ell,v})\ast f_{E_0,h_0,v}=0.
		\end{equation}	 
		The proof of Theorem \ref{thm 4.1} then follows, using Theorem \ref{thm 3.1}, \eqref{4.17}, \eqref{4.19}. 
		Let $\Delta(\tau_v, \ell-\frac{c+k}{2})$ denote the local component at $v$ of $\Delta(\tau, \ell-\frac{c+k}{2})$. By our assumption on $v$, this representation is unramified. As in the proof of Lemma \ref{lem 3.1}, it is the unramified component of the representation 
		$\Ind_{P_{(\ell-\frac{c+k}{2})^2}(F_v)}^{\GL_{2\ell-c-k}(F_v)}(\chi\circ\det\otimes\chi^{-1}\circ\det)$, and, as in the proof of Lemma \ref{lem 3.1}, the unramified component at $v$ of \eqref{4.18} is the unramified component of the representation obtained by replacing in \eqref{4.18}, at the place $v$, $\Delta(\tau_v, \ell-\frac{c+k}{2})|\det\cdot|^{-m-\frac{\ell}{2}-\frac{c+k}{4}}$ by 
		$$
		\Ind_{P_{(\ell-\frac{c+k}{2})^2}(F_v)}^{\GL_{2\ell-c-k}(F_v)}(\chi\circ\det |\det\cdot|^{-m-\frac{\ell}{2}-\frac{c+k}{4}}\otimes\chi\circ\det |\det\cdot|^{m+\frac{\ell}{2}+\frac{c+k}{4}}).
		$$
		Thus, consider right $K_{4\ell,v}$-invariant functions $\varphi_v$ on $\Sp_{4\ell}(F_v)$, which satisfy the local analog of \eqref{4.20}, with $\Delta(\tau_v, \ell-\frac{c+k}{2})|\det\cdot|^{-m-\ell-\frac{c+k}{2}}$ replaced by the last representation. Then \eqref{4.20} is replaced by\\
		\\   
		$\varphi_v(u\begin{pmatrix}a_1\\&a_2\\&&I_{2(c+k)}\\&&&a_2^*\\&&&&a_1^*\end{pmatrix}g)=$
		\begin{equation}\label{4.22}
			=\chi(\det(a_1a_2))|\det(a_1)|^{m+\ell+\frac{1}{2}}|\det(a_2)|^{3m+\ell+c+k+\frac{1}{2}}\varphi_v(g),
		\end{equation}
		for $u\in U^{4\ell}_{(\ell-\frac{c+k}{2})^2}(F_v), a_1, a_2\in \GL_{\ell-\frac{c+k}{2}}(F_v)$.
		In order to show \eqref{4.21}, it is enough to show that for all $\varphi_v$, as above, satisfying \eqref{4.22},
		\begin{equation}\label{4.23}
			\eta_{m,\ell}^\chi(\xi_{m,\ell,v})\ast \varphi_v(g_0)=0,\  \  \  g_0\in \Sp_{4\ell}(F_v).
		\end{equation}
		By the Iwasawa decomposition and \eqref{4.22}, we may assume that\\
		$g_0=diag (I_{\ell-c-k},b_0,I_{\ell-c-k}):=\tilde{b}_0$, $b_0\in \Sp_{2(c+k)}(F_v)$. Also, once again, using the Iwasawa decomposition in the integration which defines \eqref{4.23}, and using \eqref{4.22}, we get, denoting, for short, $\eta_{m,\ell}^\chi(\xi_{m,\ell,v})=\eta_{\xi_v}$,\\
		\\
		$	\eta_{\xi_v}\ast \varphi_v(\tilde{b}_0)=$
		\begin{equation}\label{4.24}
			=\int \eta_{\xi_v}^{U_{(\ell-\frac{c+k}{2})^2}}(\begin{pmatrix}a_1\\&a_2\\&&b\\&&&a_2^*\\&&&&a_1^*\end{pmatrix})\chi(\det(a_1a_2))|\det(a_1)|^{m-2\ell-\frac{c+k+1}{2}}
		\end{equation}
		$$
		|\det(a_2)|^{3m-\frac{c+k+1}{2}}\varphi_v(\tilde{b}_0\tilde{b})da_1da_2db,
		$$
		where the integration is over $a_1, a_2\in \GL_{\ell-\frac{c+k}{2}}(F_v),\ b\in \Sp_{2(c+k)}(F_v)$. Consider the inner $da_1da_2$-integration in \eqref{4.24}. We claim that, for all $0\leq k<2\ell-c$ and all $b\in \Sp_{2(c+k)}(F_v)$,
		\begin{equation}\label{4.25}
			\int \eta_{\xi_v}^{U_{(\ell-\frac{c+k}{2})^2}}(\begin{pmatrix}a_1\\&a_2\\&&b\\&&&a_2^*\\&&&&a_1^*\end{pmatrix})\chi(\det(a_1a_2))|\det(a_1)|^{m-2\ell-\frac{c+k+1}{2}}
		\end{equation}
		$$
		|\det(a_2)|^{3m-\frac{c+k+1}{2}}da_1da_2=0.
		$$
		Using the Iwasawa decomposition in $\GL_{\ell-\frac{c+k}{2}}(F_v)\times \GL_{\ell-\frac{c+k}{2}}(F_v)$, the l.h.s. of \eqref{4.25} is equal to
		\begin{equation}\label{4.26}
			\varphi_{\xi_v}(b)=	\int \delta^{-\frac{1}{2}}_{Q^{4\ell}_{1^{2\ell-c-k}}}\eta_{\xi_v}^{U_{1^{2\ell-c-k}}}(\begin{pmatrix}t_1\\&t_2\\&&b\\&&&t_2^*\\&&&&t_1^*\end{pmatrix})\chi(\det(t_1t_2))|\det(t_1)|^{m-\frac{\ell}{2}-\frac{c+k}{4}}
		\end{equation}
		$$
		|\det(t_2)|^{3m+\frac{\ell}{2}+\frac{c+k}{4}}\delta^{-\frac{1}{2}}_{B_{\GL_{\ell-\frac{c+k}{2}}}}(t_1t_2)dt_1dt_2.
		$$
		Here, $t_1, t_2$ are integrated over $T_{\ell-\frac{c+k}{2}}(F_v)$. The function $\varphi_{\xi_v}$ lies in\\ $\mathcal{H}(\Sp_{2(c+k)}(F_v)/ / K_{2(c+k),v})$. It is enough show that $\hat{\mathcal{S}}(\varphi_{\xi_v})=0$.
		We have,\\
		\\
		$\hat{\mathcal{S}}(\varphi_{\xi_v})(Z_1,...,Z_{c+k})=$
		\begin{equation}\label{4.27}
			=\hat{\mathcal{S}}(\eta^\chi_{m,\ell}(\xi_{m,\ell,v}))(\chi(p_v)q_v^{-m+\ell-\frac{1}{2}},\chi(p_v)q_v^{-m+\ell-\frac{3}{2}},...,\chi(p_v)q_v^{-m+\frac{c+k+1}{2}},
		\end{equation}	
		$$
		\chi(p_v)q_v^{-3m-\frac{c+k+1}{2}},\chi(p_v)q_v^{-3m-\frac{c+k+3}{2}},...,\chi(p_v)q_v^{-3m-\ell+\frac{1}{2}},Z_1,...,Z_{c+k})=0.
		$$
		The last equality follows from \eqref{3.20}. This proves Theorem \ref{thm 4.1}.	
		
	\end{proof}
	
	Going back to \eqref{4.7}, Prop. \ref{prop 4.1} and Theorem \ref{thm 4.1} tell us that only the $\GL_{2\ell-c}(F_v)\times \Sp_{4m}(F_v)$  - orbit of $B_{2\ell-c}=\begin{pmatrix}0&I_{2\ell-c}\\0&0\end{pmatrix}$ contributes to the expansion \eqref{4.7} of $(1\otimes \xi_{m,\ell,v})\ast \theta^{\psi_\delta}_{\Delta(\tau,m+\ell)}$. Thus,
	\begin{equation}\label{4.28}
		(1\otimes \xi_{m,\ell,v})\ast \theta^{\psi_\delta}_{\Delta(\tau,m+\ell)}(g,h)=\sum_{\beta\in Q'_{2\ell-c}(F)\backslash \Sp_{4m}(F)}(1\otimes \xi_{m,\ell,v})\ast \theta^{\psi_\delta,B_{2\ell-c}}_{\Delta(\tau,m+\ell)}(g,\beta h).
	\end{equation}
	Here, $Q'_{2\ell-c}(F)$ is the subgroup of  matrices in $\Sp_{4m}(F)$ of the form $\begin{pmatrix}I_{2\ell-c}&\ast&\ast\\&b&\ast\\&&I_{2\ell-c}\end{pmatrix}$, $b\in \Sp_{4(m-\ell)+2c}(F)$.
	By \eqref{4.17},\\
	\\
	$(1\otimes\xi_{m,\ell,v})\ast\theta^{\psi_{\delta,B_{2\ell-c}}}_{\Delta(\tau,m+\ell)}(g,h)=$
	\begin{equation}\label{4.29}
		\sum_{E\in M_{c\times k}(F)}\int (1\otimes\xi_{m,\ell,v})\ast\theta^{U_{2(2\ell-c)}}_{\Delta(\tau,m+\ell)}(\begin{pmatrix}I_{2\ell-c}&a\\&I_{2\ell-c}\\&&I_c&0&z\\&&&I_{4(m-\ell)+2c}&0\\&&&&I_c\\&&&&&I_{2\ell-c}&a'\\&&&&&&I_{2\ell-c}\end{pmatrix}\cdot
	\end{equation}
	$$
	\widehat{\omega v_1(E)}(g,h))\psi^{-1}(tr(w_c\delta z)+tr(a))dadz.
	$$
	Note the case where $c=0$, i.e. $\delta=0$. We will need it later. Then $(1\otimes \xi_{m,\ell,v})\ast \theta^{\psi_0}_{\Delta(\tau,m+\ell)}$ is the constant term of $(1\otimes \xi_{m,\ell,v})\ast \theta_{\Delta(\tau,m+\ell)}$ along $i(U^{4\ell}_{2\ell}\times 1)$. By \eqref{4.28}, \eqref{4.29},\\
	\\
	$(1\otimes \xi_{m,\ell,v})\ast \theta^{U^{4\ell}_{2\ell}\times 1}_{\Delta(\tau,m+\ell)}(g,h)=$
	\begin{equation}\label{4.30}
		=\sum_{\beta\in Q'_{2\ell}(F)\backslash \Sp_{4m}(F)}\int (1\otimes\xi_{m,\ell,v})\ast\theta^{U_{4\ell}}_{\Delta(\tau,m+\ell)}(\begin{pmatrix}I_{2\ell}&a\\&I_{2\ell}\\&&I_{4(m-\ell)}\\&&&I_{2\ell}&a'\\&&&&I_{2\ell}\end{pmatrix}(g,\beta h))\cdot
	\end{equation}
	$$
	\cdot \psi^{-1}(tr(a))da.
	$$
	Denote the integral on the r.h.s. of \eqref{4.30} by $(1\otimes\xi_{m,\ell,v})\ast\theta^{U_{4\ell},\psi_{V_{(2\ell)^2}}}_{\Delta(\tau,m+\ell)}(g,\beta h)$. By an easy check, one can rewrite \eqref{4.30} as
	\begin{equation}\label{4.31}
		(1\otimes \xi_{m,\ell,v})\ast \theta^{U^{4\ell}_{2\ell}\times 1}_{\Delta(\tau,m+\ell)}(g,h)=\sum_{\beta\in Q_{2\ell}(F)\backslash \Sp_{4m}(F)}\sum_{\alpha\in \GL_{2\ell}(F)} (1\otimes\xi_{m,\ell,v})\ast\theta^{U_{4\ell},\psi_{V_{(2\ell)^2}}}_{\Delta(\tau,m+\ell)}(\hat{\alpha}g,\beta h).
	\end{equation}

	\section{ Proof of Theorem \ref{thm 2.1}: Rapid decrease of $(1\otimes \xi_{m,\ell,v})\ast \theta_{\Delta(\tau,m+\ell)}(g,h)$ in $g$}
	
	We keep the notation above. In this section we prove
	\begin{thm}\label{thm 5.1}
		For a given $h\in \Sp_{4m}(\BA)$, the function on $\Sp_{4\ell}(F)\backslash \Sp_{4\ell}(\BA)$, $g\mapsto (1\otimes \xi_{m,\ell,v})\ast \theta_{\Delta(\tau,m+\ell)}(g,h)$ is rapidly decreasing.
	\end{thm}
	Fix a Siegel domain $\mathfrak{S}_{4\ell}=\Omega_{4\ell} T_{4\ell}^+(\epsilon_0)K_{4\ell}$ in $\Sp_{4\ell}(\BA)$, where $\Omega_{4\ell}$ is a sufficiently large compact subset of $B_{\Sp_{4\ell}}(\BA)$, $\epsilon_0>0$ is sufficiently small, and $T_{4\ell}^+(\epsilon_0)$ is the subgroup of diagonal matrices $\hat{t}=diag(t_1,...,t_{2\ell},t^{-1}_{2\ell},...,t_1^{-1})$, such that each $t_i=\prod_\nu t_{i,\nu}\in \BA^*$
	satisfies $t_{i,\nu}=1$, for all $\nu<\infty$, and at the set of archimdean places, $S_\infty$, $t_{i,\nu}=a_i$, for all $\nu\in S_\infty$, where $a_i>0$, and we have, 
	\begin{equation}\label{5.1}
		\frac{a_i}{a_{i+1}}>\epsilon_0,\ \ i=1,2,...,2\ell-1;\ \ a_{2\ell}>\epsilon_0.
	\end{equation} 
	We will denote by $||\cdot||$ the norm on our adelic groups $\Sp_{2k}(\BA)$, or $\GL_k(\BA)$, as in \cite{MW95}, I.2.2. Similarly, fix a Siegel domain $\mathfrak{S}_{4m}=\Omega_{4m} T_{4m}^+(\epsilon_0)K_{4m}$ in $\Sp_{4m}(\BA)$. We will prove the rapid decrease in $g$, in Theorem \ref{thm 5.1}, for a fixed $h\in \mathfrak{S}_{4m}$. 
	
	Assume that $g\in \mathfrak{S}_{4\ell}$. Since $(1\otimes \xi_{m,\ell,v})\ast \theta_{\Delta(\tau,m+\ell)}$ is $K_{4\ell}\times 1$-finite, we may assume that $g=bt$, where $b\in \Omega$ and $\hat{t}\in T^+(\epsilon_0)$.
	We start with the Fourier expansion \eqref{4.5},\\
	\\
	$(1\otimes \xi_{m,\ell,v})\ast\theta_{\Delta(\tau,m+\ell)}(g,h)=$
	\begin{equation}\label{5.2}
		=\sum_{c=0}^{2\ell}\sum_{[\delta]\in [T_c(F)]}\sum_{\gamma\in P^\delta_{2\ell-c,c}(F)\backslash \GL_{2\ell}(F)}(1\otimes \xi_{m,\ell,v})\ast\theta^{\psi_\delta}_{\Delta(\tau,m+\ell)}((\hat{\gamma}b\hat{t},h)).
	\end{equation}
	In \eqref{5.2}, consider the terms with $1\leq c\leq 2\ell$, and let $P^\delta_{2\ell-c,c}(F)\gamma$ be a coset in $P^\delta_{2\ell-c,c}(F)\backslash \GL_{2\ell}(F)$. Write $\gamma=\begin{pmatrix}\gamma_1\\\gamma_2\end{pmatrix}$, where $\gamma_1\in M_{(2\ell-c)\times 2\ell}(F)$ and $\gamma_2\in M_{c\times 2\ell}(F)$. Denote by $\gamma_2^1$ the first column of $\gamma_2$. Although $\gamma_2^1$ depends on the representative $\gamma$, the property of being nonzero depends only on the coset $P^\delta_{2\ell-c,c}(F)\gamma$. We will prove that for each $0\leq c\leq 2\ell$, the corresponding term in \eqref{5.2} is rapidly decreasing in $g$. For $1\leq c\leq 2\ell$, we will prove this separately for the corresponding sums over $\gamma_2^1\neq 0$ (Prop. \ref{prop 5.1}) and then for $\gamma_2^1=0$ (Prop. \ref{prop 5.2}). In Prop. \ref{prop 5.3}, we will treat the case $c=0$.
	\begin{prop}\label{prop 5.1}
		Let $1\leq c\leq 2\ell$. There is $A>0$, and for every integer $N\geq 1$, there exists $k_N>0$, such that, for all $b\in \Omega_{4\ell}$, $\hat{t}\in T_{4\ell}^+(\epsilon_0)$, $h\in \Sp_{4m}(\BA)$,
		$$
		\sum_{[\delta]\in [T_c(F)]}\sum_{\gamma\in P^\delta_{2\ell-c,c}(F)\backslash \GL_{2\ell}(F); \gamma_2^1\neq 0}|(1\otimes \xi_{m,\ell,v})\ast\theta^{\psi_\delta}_{\Delta(\tau,m+\ell)}((\hat{\gamma}b\hat{t},h))|\leq k_N a_1^{-N}||h||^A	
		$$
		
	\end{prop}
	\begin{proof}
		We will prove the proposition with $\theta_{\Delta(\tau,m+\ell)}$ in place of $(1\otimes \xi_{m,\ell,v})\ast\theta_{\Delta(\tau,m+\ell)}$. Since $\theta_{\Delta(\tau,m+\ell)}$ is smooth, there are $k_1, A>0$, such that, 
		\begin{equation}\label{5.3}
			\sum_{[\delta]\in [T_c(F)]}	\sum_{\gamma\in P^\delta_{2\ell-c,c}(F)\backslash \GL_{2\ell}(F)}|\theta^{\psi_\delta}_{\Delta(\tau,m+\ell)}((\hat{\gamma}bt,h))|\leq k_1 ||\begin{pmatrix}t\\&h\\&&t^*\end{pmatrix}||^A.
		\end{equation}
		Indeed, the series of absolute values \eqref{5.3} is part of the full series of absolute values resulting from the Fourier expansion of $\theta_{\Delta(\tau,m+\ell)}$ along $U^{4\ell}_{2\ell}(F)\backslash U^{4\ell}_{2\ell}(\BA)\times I_{4m}$. By \eqref{5.1}, 
		$$
		||\begin{pmatrix}t\\&h\\&&t^*\end{pmatrix}||=max\{a_1^{\pm 1},...,a_{2\ell}^{\pm 1},||h||\}\leq k_2a_1||h||,
		$$
		where $k_2=\epsilon_0^{-4\ell}$. Thus, from \eqref{5.3}, with $k=k_1k_2^A$,
		\begin{equation}\label{5.4}
			\sum_{[\delta]\in [T_c(F)]}	\sum_{\gamma\in P^\delta_{2\ell-c,c}(F)\backslash \GL_{2\ell}(F)}|\theta^{\psi_\delta}_{\Delta(\tau,m+\ell)}((\hat{\gamma}b\hat{t},h))|\leq k (a_1||h||)^A.
		\end{equation}	
		Since $\theta_{\Delta(\tau,m+\ell)}$ is smooth, we can write it, by the lemma of Dixmier-Malliavin, as a finite sum of convolutions 
		$$
		\phi\ast \theta'_{\Delta(\tau,m+\ell)}=\int_{S_{2\ell}(\BA)}\phi(y)\rho((u_{2\ell}^{4\ell}(y),I_{4m}))\theta'_{\Delta(\tau,m+\ell)}dy
		$$
		where $\phi\in \mathcal{S}(S_{2\ell}(\BA))$ - the space of Schwartz functions on $U^{4\ell}_{2\ell}(\BA)$. Here, $\rho(x)$ denotes the right translation by $x$. Thus, let us replace $\theta_{\Delta(\tau,m+\ell)}$ by $\phi\ast \theta'_{\Delta(\tau,m+\ell)}$ and consider\\
		\\
		$\sum_{[\delta]\in [T_c(F)]}\sum_{\gamma\in P^\delta_{2\ell-c,c}(F)\backslash \GL_{2\ell}(F); \gamma_2^1\neq 0} \phi\ast\theta^{'\psi_\delta}_{\Delta(\tau,m+\ell)}((\hat{\gamma}b\hat{t},h))=$	
		\begin{equation}\label{5.5}
			=\sum_{[\delta]\in [T_c(F)]}\sum_{\gamma\in P^\delta_{2\ell-c,c}(F)\backslash \GL_{2\ell}(F); \gamma_2^1\neq 0}\int\phi(y)\theta'_{\Delta(\tau,m+\ell)}((u_{2\ell}^{4\ell}(x)\hat{\gamma}b\hat{t}u_{2\ell}^{4\ell}(y),h))
		\end{equation}
		$$
		\psi^{-1}_\delta(u_{2\ell}^{4\ell}(x))dydx,
		$$
		where $x$ is integrated along $S_{2\ell}(F)\backslash S_{2\ell}(\BA)$ and $y$ is integrated along $S_{2\ell}(\BA)$. Write $b=\hat{b}_1u_{2\ell}^{4\ell}(e)$, where $b_1\in B_{\GL_{2\ell}}(\BA)$. Simple conjugations inside \eqref{5.5} yield, with the same domains of summation and integrations,
		$$
		\sum\int\phi(y)\theta'_{\Delta(\tau,m+\ell)}((u_{2\ell}^{4\ell}(x+\gamma b_1(ty(t^*)^{-1}+e)(\gamma^*b_1^*)^{-1}\hat{\gamma}b\hat{t},h))\psi_\delta^{-1}(x)dydx.
		$$
		Switching the integrations order and changing variables in $x$, the $dy$ integration results in a Fourier transform of $\phi$ at
		$$
		(\gamma^*b_1^*t^*)^{-1}\begin{pmatrix}0&w_c\delta\\0&0\end{pmatrix}\gamma b_1t=w_{2\ell}t{}^tb_1{}^t\gamma_2\delta\gamma_2b_1t.
		$$
		Then \eqref{5.5} becomes\\
		\\
		$\sum_{[\delta]\in [T_c(F)]}\sum_{\gamma\in P^\delta_{2\ell-c,c}(F)\backslash \GL_{2\ell}(F); \gamma_2^1\neq 0} \phi\ast\theta^{'\psi_\delta}_{\Delta(\tau,m+\ell)}((\hat{\gamma}b\hat{t},h))=$	
		\begin{equation}\label{5.6}
			=\sum_{[\delta]\in [T_c(F)]}\sum_{\gamma\in P^\delta_{2\ell-c,c}(F)\backslash \GL_{2\ell}(F); \gamma_2^1\neq 0}\psi_\delta(\gamma b_1e(\gamma^*b_1^*)^{-1})\hat{\phi}(w_{2\ell}t{}^tb_1{}^t\gamma_2\delta\gamma_2b_1t)
		\end{equation}	
		$$	
		\theta^{'\psi_\delta}_{\Delta(\tau,m+\ell)}((\hat{\gamma}b\hat{t},h)).
		$$
		Since $\psi_\delta(\gamma b_1e(\gamma^*b_1^*)^{-1})$ has absolute value $1$, the above calculation also shows that\\
		\\
		$\sum_{[\delta]\in [T_c(F)]}\sum_{\gamma\in P^\delta_{2\ell-c,c}(F)\backslash \GL_{2\ell}(F); \gamma_2^1\neq 0} |\phi\ast\theta^{\psi_\delta}_{\Delta(\tau,m+\ell)}((\hat{\gamma}b\hat{t},h))|=$	
		\begin{equation}\label{5.7}
			=\sum_{[\delta]\in [T_c(F)]}\sum_{\gamma\in P^\delta_{2\ell-c,c}(F)\backslash \GL_{2\ell}(F); \gamma_2^1\neq 0}|\hat{\phi}(w_{2\ell}{}^t(\gamma_2b_1t)\delta\gamma_2b_1t)| |\theta^{'\psi_\delta}_{\Delta(\tau,m+\ell)}((\hat{\gamma}b\hat{t},h))|.
		\end{equation}
		Since $\gamma_2$ has rank $c$ and $\gamma_2^1\neq 0$, the first column of ${}^t\gamma_2\delta\gamma_2$ is nonzero.
		In \eqref{5.7}, $\hat{\phi}$ being a Schwartz function, we can estimate in \eqref{5.7}, for any integer $N\geq 1$,
		\begin{equation}\label{5.8}
			|\hat{\phi}(w_{2\ell}{}^t(\gamma_2b_1t)\delta\gamma_2b_1t)|\leq k'_N ||w_{2\ell}{}^t(\gamma_2b_1t)\delta\gamma_2b_1t||^{-N}_{max},
		\end{equation}
		for a suitable $k'_N>0$. Here, $||\cdot||_{max}=\prod_v||\cdot||_{max,v}$ denotes the product over all places of the local maximum norms. We have, for a suitable $k_3>0$, and all $t, b_1, \gamma_2$ as in \eqref{5.7},
		\begin{equation}\label{5.9}
			||w_{2\ell}{}^t(\gamma_2b_1t)\delta\gamma_2b_1t||_{max}=||{}^t(t^{-1}b_1t)t({}^t\gamma_2\delta\gamma_2)t(t^{-1}b_1t)||_{max}\geq k_3||t{}^t\gamma_2\delta\gamma_2)t||_{max}.
		\end{equation}	
		Indeed, since $b_1$ lies in a compact subset $\Omega_1$ of $B_{\GL_{2\ell}}(\BA)$, then, by \eqref{5.1}, the conjugation $t^{-1}b_1t$ shrinks the coordinates of $b_1$, and hence $t^{-1}b_1t$ lies in a compact subset $\Omega'_1\subset \Omega_1$. Clearly,
		$$
		||t{}^t\gamma_2\delta\gamma_2)t||_{max}\geq max\{a_1a_j\cdot|{}^t\gamma_2^1\delta\gamma_2^j|\}_{j\leq 2\ell}\geq \epsilon_0^{2\ell}a_1\cdot max\{|{}^t\gamma_2^1\delta\gamma_2^j|\}_{j\leq 2\ell},
		$$
		where $\gamma_2^j$ denotes the $j$-th column of $\gamma_2$. By \eqref{5.9}, there is $k_4>0$, such that, for $t, b_1, \gamma_2$ as in \eqref{5.7},
		\begin{equation}\label{5.10}
			||w_{2\ell}{}^t(\gamma_2b_1t)\delta\gamma_2b_1t||_{max}\geq k_4	a_1\cdot max\{|{}^t\gamma_2^1\delta\gamma_2^j|\}_{j\leq 2\ell}.
		\end{equation} 
		In \eqref{5.7}, since $\hat{\phi}$ is a Schwartz function, the coordinates of ${}^t(\gamma_2b_1t)\delta\gamma_2b_1t)$ must lie in fixed compact subsets $C_v\subset S_{2\ell}(F_v)$, at each finite place $v$. Since $b_1\in \Omega_1$, and since $t_v=I_{2\ell}$, for all $v<\infty$, we conclude that in \eqref{5.7}, the coordinates of ${}^t\gamma_2\delta\gamma_2$ must lie inside a lattice	$L\subset \BA_\infty$. Let $k_0=min\{|x|_\infty\}_{0\neq x\in L}$. Then $k_0>0$ and inside the support of $\hat{\phi}$ in \eqref{5.7}, we have, by \eqref{5.10},
		\begin{equation}\label{5.11}
			||w_{2\ell}{}^t(\gamma_2b_1t)\delta\gamma_2b_1t||_{max}\geq k_5	a_1,
		\end{equation} 
		for a suitable positive constant $k_5$, which depends on $\phi$ and $\epsilon_0$. Note, again, that in \eqref{5.7}, we take $\gamma_2^1\neq 0$, and hence the r.h.s. of \eqref{5.10} is nonzero. Using \eqref{5.11} in \eqref{5.8}, we get that for every positive integer $N$,
		\begin{equation}\label{5.12}
			|\hat{\phi}(w_{2\ell}{}^t(\gamma_2b_1t)\delta\gamma_2b_1t)|\leq k''_Na_1^{-N},
		\end{equation}
		for a suitable positive constant $k''_N$, which depends on $\phi$, $\epsilon_0$ and $N$. By \eqref{5.4}, \eqref{5.7}, \eqref{5.12}, we get that there is $A>0$, and for every positive integer $N$, there is $c_N>0$, such that, for all $b\in \Omega$, $\hat{t}\in T^+(\epsilon_0)$, $h\in \Sp_{4m}(\BA)$,
		$$
		\sum_{[\delta]\in [T_c(F)]}\sum_{\gamma\in P^\delta_{2\ell-c,c}(F)\backslash \GL_{2\ell}(F); \gamma_2^1\neq 0} |\phi\ast\theta'^{\psi_\delta}_{\Delta(\tau,m+\ell)}((\hat{\gamma}b\hat{t},h))|\leq k_Na_1^{-N+A}||h||^A.
		$$ 
		This completes the proof of Proposition \ref{prop 5.1}.	
		
	\end{proof}
	
	We now consider the sum over the terms in \eqref{5.2} with $\gamma_2^1=0$. This forces $c<2\ell$. For the next proposition, we take $h\in \mathfrak{S}_{4m}$, and, again, it is enough to take $h$ of the form $h=b'\hat{t'}$, where $b'\in \Omega_{4m}$, and $\hat{t'}=diag(t'_1,...,t'_{2m},(t')^{-1}_{2m},...,(t')_1^{-1})\in T^+_{4m}(\epsilon_0)$. Thus, all finite coordinates of $t'_i$ are $1$, and at $v\in S_\infty$, $t'_{i,v}=a'_i>0$, satisfying the inequalities analogous to \eqref{5.1}.
	\begin{prop}\label{prop 5.2}
		Let $1\leq c< 2\ell$. There is $A>0$, and for every integer $N\geq 1$, there exists $c_N>0$, such that, for all $b\in \Omega_{4\ell}$, $b'\in \Omega_{4m}$, $\hat{t}\in T_{4\ell}^+(\epsilon_0)$, $\hat{t'}\in T_{4m}^+(\epsilon_0)$,
		$$
		\sum_{[\delta]\in [T_c(F)]}\sum_{\gamma\in P^\delta_{2\ell-c,c}(F)\backslash \GL_{2\ell}(F); \gamma_2^1= 0}|(1\otimes \xi_{m,\ell,v})\ast\theta^{\psi_\delta}_{\Delta(\tau,m+\ell)}((\hat{\gamma}b\hat{t},b'\hat{t'}))|\leq c_N a_1^{-N}(a_1')^{A+N}	
		$$	
	\end{prop}
	\begin{proof}
		For $\gamma\in P^\delta_{2\ell-c,c}(F)\backslash \GL_{2\ell}(F)$, with $\gamma_2^1=0$, we may take $\gamma=\begin{pmatrix}1\\&\gamma'\end{pmatrix}$, $\gamma'\in P^\delta_{2\ell-1-c,c}(F)\backslash \GL_{2\ell-1}(F)$. By \eqref{4.28}, the Fourier expansion of the following function on $F^{4m}\backslash \BA^{4m}$,
		$$
		x\mapsto  (1\otimes \xi_{m,\ell,v})\ast\theta^{\psi_\delta}_{\Delta(\tau,m+\ell)}(\begin{pmatrix}1&0&x&0&\ast\\&I_{2\ell-1}&0&0&0\\&&I_{4m}&0&x'\\&&&I_{2\ell-1}&0\\&&&&1\end{pmatrix}(\hat{\gamma}b\hat{t},h)):=
		$$
		$$
		:=(1\otimes \xi_{m,\ell,v})\ast\theta^{\psi_\delta}_{\Delta(\tau,m+\ell)}(\bar{y}(x)(\hat{\gamma}b\hat{t},h))
		$$
		at $x=0$ is
		$$
		\sum_{0\neq e\in F^{4m}}\int_{F^{4m}\backslash \BA^{4m}}(1\otimes \xi_{m,\ell,v})\ast\theta^{\psi_\delta}_{\Delta(\tau,m+\ell)}(\bar{y}(x)(\hat{\gamma}b\hat{t},h))\psi^{-1}(x\cdot {}^te)dx.
		$$
		The point here is that the trivial character, i.e. $e=0$ does not contribute. Write, in the last sum, $e=(1,0,...,0){}^t\eta^{-1}$, with $\eta\in Q'_1(F)\backslash \Sp_{4m}(F)$, where $Q'_1(F)$ is the subgroup of matrices in $\Sp_{4m}(F)$ of the form $\begin{pmatrix}1&\ast&\ast\\&\eta'&\ast\\&&1\end{pmatrix}, \eta'\in \Sp_{4m-2}(F)$. Then the last Fourier expansion is\\
		\\
		$(1\otimes \xi_{m,\ell,v})\ast\theta^{\psi_\delta}_{\Delta(\tau,m+\ell)}((\hat{\gamma}b\hat{t},h))=$
		\begin{equation}\label{5.13}
			\sum_{\eta\in Q'_1(F)\backslash \Sp_{4m}(F)}\int_{F^{4m}\backslash \BA^{4m}}(1\otimes \xi_{m,\ell,v})\ast\theta^{\psi_\delta}_{\Delta(\tau,m+\ell)}(\bar{y}(x)(\hat{\gamma}b\hat{t},\eta h))\psi^{-1}(x_1)dx.
		\end{equation}
		By \eqref{5.13}, it suffices prove that, for a given $\theta'_{\Delta(\tau,m+\ell)}\in \Theta_{\Delta(\tau,m+\ell)}$, there is $A>0$, and for every integer $N\geq 1$, there exists $k_N>0$, such that, for all $b\in \Omega$, $\hat{t}\in T^+(\epsilon_0)$, $h\in \Sp_{4m}(\BA)$,
		\begin{equation}\label{5.14}
			\sum_{[\delta]\in [T_c(F)]}\sum_{\gamma\in P^\delta_{2\ell-c,c}(F)\backslash \GL_{2\ell}(F); \gamma_2^1= 0}\sum_{\eta\in Q'_1(F)\backslash \Sp_{4m}(F)}| (\theta')^{\tilde{\psi}_\delta}_{\Delta(\tau,m+\ell)}((\hat{\gamma}b\hat{t},\eta h))|	
		\end{equation}	
		$$
		\leq k_N a_1^{-N}||h||^{A+N},
		$$
		where
		$$
		(\theta')^{\tilde{\psi}_\delta}_{\Delta(\tau,m+\ell)}((\hat{\gamma}b\hat{t},\eta h))=\int_{F^{4m}\backslash \BA^{4m}} (\theta')^{\psi_\delta}_{\Delta(\tau,m+\ell)}(\bar{y}(x)(\hat{\gamma}b\hat{t},\eta h))\psi^{-1}(x_1)dx.
		$$
		As in the proof of Prop. \ref{prop 5.1}, By the lemma of Dixmier-Malliavin, it is enough to prove \eqref{5.14}, with $\theta'_{\Delta(\tau,m+\ell)}=f_{\phi,\varphi}\ast \theta''_{\Delta(\tau,m+\ell)}$, 
		where
		$$
		f_{\phi,\varphi}\ast \theta''_{\Delta(\tau,m+\ell)}=\int_{\BA^{4m+1}}\phi(u)\varphi(z)\rho(\begin{pmatrix}1&0&u&0&z\\&I_{2\ell-1}&0&0&0\\&&I_{4m}&0&u'\\&&&I_{2\ell-1}&0\\&&&&1\end{pmatrix})\theta''_{\Delta(\tau,m+\ell)}d(u,z),
		$$
		$\phi\in \mathcal{S}(\BA^{4m})$, $\varphi\in \mathcal{S}(\BA)$. We have, for $\gamma=\begin{pmatrix}1\\&\gamma'\end{pmatrix}$, as in the beginning of the proof,
		$$
		\int_{F^{4m}\backslash \BA^{4m}}(f_{\phi,\varphi}\ast \theta''_{\Delta(\tau,m+\ell)})^{\psi_\delta}(\bar{y}(x)(\hat{\gamma}b\hat{t},\eta h))\psi^{-1}(x_1)dx=
		$$
		$$
		=\alpha_\varphi\int_{F^{4m}\backslash \BA^{4m}}\int_{\BA^{4m}}\phi(u)(\theta'')^{\psi_\delta}_{\Delta(\tau,m+\ell)}(\bar{y}(x+b_{1,1}a_1(\eta h)^{-1})(\hat{\gamma}b\hat{t},\eta h))\psi^{-1}(x_1)dudx,
		$$
		where $\alpha_\varphi=\int_\BA\varphi(z)dz$. Switching order of integration and changing variables, we get
		$$
		\hat{\phi}(b_{1,1}t_1h^{-1}\eta^{-1}\cdot e_1)\int_{F^{4m}\backslash \BA^{4m}} (\theta'')^{\psi_\delta}_{\Delta(\tau,m+\ell)}(\bar{y}(x)(\hat{\gamma}b\hat{t},\eta h))\psi^{-1}(x_1)dx,
		$$
		where $e_1=\begin{pmatrix}1\\0\\ \vdots\\0\end{pmatrix}$. Thus, we can majorize the l.h.s. of \eqref{5.14} by
		\begin{equation}\label{5.15}
			\sum_{[\delta]\in [T_c(F)]}\sum_{\gamma'\in P^\delta_{2\ell-1-c,c}(F)\backslash \GL_{2\ell-1}(F)}\sum_{\eta\in Q'_1(F)\backslash \Sp_{4m}(F)}|\hat{\phi}(b_{1,1}t_1h^{-1}\eta^{-1}\cdot e_1)|\cdot
		\end{equation}	
		$$
		\cdot\int_{F^{4m}\backslash \BA^{4m}}|(\theta'')^{\psi_\delta}_{\Delta(\tau,m+\ell)}(\bar{y}(x)(\widehat{\begin{pmatrix}1\\&\gamma'\end{pmatrix}}b\hat{t},\eta h))|dx.
		$$
		Substitute $h=b'\hat{t'}$, $b'\in \Omega_{4m}$, $\hat{t'}\in T^+_{4m}(\epsilon_0)$. Exactly as in the last proof, since $\hat{\phi}$ is a Schwartz function, we conclude that the coordinates of $\eta^{-1}\dot e_1$ in \eqref{5.15} must lie inside a lattice $L\subset \BA_\infty$. In \eqref{5.15}, inside the support of $\hat{\phi}$, we have, for positive constants $k_1,...,k_4$,
		$$
		||b_{1,1}t_1(t')^{-1}(b')^{-1}\eta^{-1}\cdot e_1||_{max}\geq k_1a_1||((t')^{-1}(b')^{-1}t')(t')^{-1}\eta^{-1}\cdot e_1||_{max}\geq 
		$$
		$$
		\geq k_2a_1||(t')^{-1}\eta^{-1}\cdot e_1||_{max}\geq k_3a_1(a'_1)^{-1}||\eta^{-1}\cdot e_1||_{max}\geq k_4a_1(a'_1)^{-1}.
		$$
		Here, we used that, since $b\in \Omega_{4\ell}$, $b_{1,1}$ lies in a compact subset of $\BA^*$. Since $b'\in \Omega_{4m}$, and the archimedean coordinates $a'_i$ satisfy the analog of \eqref{5.1}, the conjugation $(t')^{-1}(b')^{-1}t')$ lies in a compact subset $\Omega'\subset \Omega_{4m}$. Finally, let $k_0=min \{|z|_\infty |\ 0\neq z\in L\}$. Then $k_0>0$ and $||\eta^{-1}\cdot e_1||_{max}\geq k_0k_0'$, where $k'_0>0$ depends on the (compact) support of $\hat{\phi}$ at the finite places. We take $k_4=k_0k_0'k_3$. Thus, for every positive integer $N$, 
		\begin{equation}\label{5.16}
			|\hat{\phi}(b_{1,1}t_1(t')^{-1}(b')^{-1}\eta^{-1}\cdot e_1)|\leq k'_Na_1^{-N}(a'_1)^N,
		\end{equation}
		for a suitable positive constant $k'_N$, which depends on $\phi$, $\epsilon_0$ and $N$. Using \eqref{5.16}, we can majorize \eqref{5.15} by
		\begin{equation}\label{5.17}
			k'_N a_1^{-N}(a'_1)^N\sum_{[\delta]\in [T_c(F)]} \sum_{\gamma'}\sum_\eta \int_{F^{4m}\backslash \BA^{4m}}|(\theta'')^{\psi_\delta}_{\Delta(\tau,m+\ell)}(\bar{y}(x)(\widehat{\begin{pmatrix}1\\&\gamma'\end{pmatrix}}b\hat{t},\eta b'\hat{t'}))|dx
		\end{equation}
		where  $\gamma'$ runs over $P^\delta_{2\ell-1-c,c}(F)\backslash \GL_{2\ell-1}(F)$ and $\eta$ - over $Q'_1(F)\backslash \Sp_{4m}(F)$. As in \eqref{5.3}, there are $k_5, k_6, A>0$, such that the last triple series is majorized by
		$$
		k_5||\begin{pmatrix}t\\&b'\hat{t'}\\&&t^*\end{pmatrix}||^A\leq k_6 a^A_1(a'_1)^A.
		$$
		This proves Proposition \ref{prop 5.2}.
	\end{proof}
	
	It remains to consider the term in \eqref{5.2} with $c=0$. This is exactly\\
	$(1\otimes \xi_{m,\ell,v})\ast \theta^{U^{4\ell}_{2\ell}\times 1}_{\Delta(\tau,m+\ell)}(g,h)$.
	\begin{prop}\label{prop 5.3}
		There is $A>0$, and for every integer $N\geq 1$, there exists $k_N>0$, such that, for all $b\in \Omega_{4\ell}$, $b'\in \Omega_{4m}$, $\hat{t}\in T_{4\ell}^+(\epsilon_0)$, $\hat{t'}\in T_{4m}^+(\epsilon_0)$,
		$$
		|(1\otimes \xi_{m,\ell,v})\ast \theta^{U^{4\ell}_{2\ell}\times 1}_{\Delta(\tau,m+\ell)}(b\hat{t},b'\hat{t'})|\leq k_N a_1^{-N}(a_1')^{A+N}	
		$$	
	\end{prop}
	\begin{proof}
		By \eqref{4.7} and Theorem \ref{thm 4.1} with $c=0$, we may use \eqref{4.14} with $k=2\ell$, and hence it suffices to prove that, for a given $\theta'_{\Delta(\tau,m+\ell)}\in \Theta_{\Delta(\tau,m+\ell)}$, there is $A>0$, and for every integer $N\geq 1$, there exists $k_N>0$, such that, for $b,t,b',t'$ as above,
		\begin{equation}\label{5.18}
			\sum|\int \theta'_{\Delta(\tau,m+\ell)}(\begin{pmatrix}I_{2\ell}&y&z\\&I_{4m}&y'\\&&I_{2\ell}\end{pmatrix}(b\hat{t},\beta b'\hat{t'})) \psi^{-1}(tr(y\cdot \begin{pmatrix}I_{2\ell}\\0\end{pmatrix}))d(y,z)|\leq 
		\end{equation}
		$$
		\leq k_N a_1^{-N}(a_1')^{A+N},	
		$$		
		where $\beta$ runs over $Q'_{2\ell}(F)\backslash \Sp_{4m}(F)$. As in the proofs of the last two propositions, by the lemma of Dixmier-Malliavin, it is enough to prove \eqref{5.18} for $\theta'_{\Delta(\tau,m+\ell)}=\varphi\ast \theta''_{\Delta(\tau,m+\ell)}$, where $\varphi\in \mathcal{S}(U_{2\ell}^{4(m+\ell)}(\BA))$, and
		$$
		\varphi\ast \theta''_{\Delta(\tau,m+\ell)}(x)=\int_{U_{2\ell}^{4(m+\ell)}(\BA)}\varphi(u)\theta''_{\Delta(\tau,m+\ell)}(xu)du.
		$$
		Denote, for $y\in M_{2\ell \times 4m}(\BA)$,
		$$
		\phi(y)=\int_{S_{2\ell}(\BA)}\varphi(\begin{pmatrix}I_{2\ell}&0&z\\&I_{4m}&0\\&&I_{2\ell}\end{pmatrix}\begin{pmatrix}I_{2\ell}&y&s(y)\\&I_{4m}&y'\\&&I_{2\ell}\end{pmatrix})dz,
		$$
		where we may take $s(y)=\frac{1}{2}yJ_{4m}{}^tyw_{2\ell}$. Then $\phi\in \mathcal{S}(M_{2\ell \times 4m}(\BA))$. Write $b=\hat{b}_1u_{2\ell}^{4\ell}(e)$, where $b_1\in B_{\GL_{2\ell}}(\BA)$. Substituting $\varphi\ast \theta''_{\Delta(\tau,m+\ell)}$ in the integral in \eqref{5.18}, the integral becomes
		$$
		\hat{\phi}((\hat{t'})^{-1}(b')^{-1}\beta^{-1}\cdot\begin{pmatrix}I_{2\ell}\\0\end{pmatrix}b_1t)\int \theta''_{\Delta(\tau,m+\ell)}(\begin{pmatrix}I_{2\ell}&y&z\\&I_{4m}&y'\\&&I_{2\ell}\end{pmatrix}(b\hat{t},\beta b'\hat{t'}))\cdot 
		$$
		$$
		\psi^{-1}(tr(y\cdot \begin{pmatrix}I_{2\ell}\\0\end{pmatrix}))d(y,z).
		$$
		Thus,we can majorize the l.h.s. of \eqref{5.18} by,
		\begin{equation}\label{5.19}
			\sum |\hat{\phi}((\hat{t'})^{-1}(b')^{-1}\beta^{-1}\cdot\begin{pmatrix}I_{2\ell}\\0\end{pmatrix}b_1t)|\int|	\theta''_{\Delta(\tau,m+\ell)}(\begin{pmatrix}I_{2\ell}&y&z\\&I_{4m}&y'\\&&I_{2\ell}\end{pmatrix}(b\hat{t},\beta b'\hat{t'}))|d(y,z).
		\end{equation}
		As in \eqref{5.16}, for every positive integer $N$, there is a $k'_N>0$, such that
		\begin{equation}\label{5.20}
			|\hat{\phi}((\hat{t'})^{-1}(b')^{-1}\beta^{-1}\cdot\begin{pmatrix}I_{2\ell}\\0\end{pmatrix}b_1t)|\leq k'_Na_1^{-N}(a'_1)^N.
		\end{equation}	
		Now, we finish the proof as in the end of Proposition \ref{prop 5.2}.
	\end{proof}
	
	The proof of Theorem \ref{thm 5.1} is now complete.
	
	Let us complete the proof of Theorem \ref{thm 2.1}. Let $m\geq 2\ell$. Assume that $\theta_{\Delta(\tau,m+\ell)}$ is right $i(K_{4\ell,v}\times K_{4m,v})$-invariant, and assume that the integral \eqref{2.2} converges absolutely. By Theorem \ref{thm 5.1}, the integral \eqref{2.2}, with $(1\otimes \xi_{m,\ell,v})\ast \theta_{\Delta(\tau,m+\ell)}$ in place of $\theta_{\Delta(\tau,m+\ell)}$ always converges absolutely. Then we have\\
	\\
	$E((1\otimes\xi_{m,\ell,v})\ast\theta_{\Delta(\tau,m+\ell)}, \theta_{\Delta(\tau,\ell)};h)=$
	$$
	=\int_{\Sp_{4\ell}(F)\backslash \Sp_{4\ell}(\BA)}(1\otimes\xi_{m,\ell,v})\ast\theta_{\Delta(\tau,m+\ell)}(i(g,h))\theta_{\Delta(\tau,\ell)}(g)dg=
	$$
	$$
	=\int_{\Sp_{4\ell}(F)\backslash \Sp_{4\ell}(\BA)}(\eta^\chi_{m,\ell}(\xi_{m,\ell,v})\otimes 1)\ast\theta_{\Delta(\tau,m+\ell)}(i(g,h))\theta_{\Delta(\tau,\ell)}(g)dg=
	$$
	$$
	=\int_{\Sp_{4\ell}(F)\backslash \Sp_{4\ell}(\BA)}\theta_{\Delta(\tau,m+\ell)}(i(g,h))(\eta^\chi_{m,\ell}(\xi_{m,\ell,v})\ast\theta_{\Delta(\tau,\ell)})(g)dg=
	$$
	\ \ \ \ \ \ \ \ \ $=c_{m,\ell,v}E(\theta_{\Delta(\tau,m+\ell)}, \theta_{\Delta(\tau,\ell)};h)$, \\
	where
	\begin{equation}\label{5.21}
		c_{m,\ell,v}=\hat{\mathcal{S}}(\eta^\chi_{m,\ell}(\xi_{m,\ell,v}))(\chi(p_v)q_v^{-\ell+\frac{1}{2}},\chi(p_v)q_v^{-\ell+\frac{3}{2}},...,\chi(p_v)q_v^{\ell-\frac{1}{2}})\neq 0,
	\end{equation}
	by \eqref{3.20}. The reason that $c_{m,\ell,v}\neq 0$ is our assumption that $m\geq 2\ell$. More generally, for any $\ell\leq m$, in the notation of Theorems \ref{thm 2.3}, \ref{thm 2.4}, the same calculation as the last one shows that if the integral
	$$
	\int\limits_{\Sp_{4\ell}(F)\backslash \Sp_{4\ell}(\BA)}\theta_{\Delta(\tau,m+\ell)}(i(g,h))E(f_{\Delta(\tau,\ell),s};g)dg
	$$
	converges absolutely for some $s$, then
	$$
	\int\limits_{\Sp_{4\ell}(F)\backslash \Sp_{4\ell}(\BA)}(1\otimes\xi_{m,\ell,v})\ast\theta_{\Delta(\tau,m+\ell)}(i(g,h))E(f_{\Delta(\tau,\ell),s};g)dg=
	$$
	$$
	=\int\limits_{\Sp_{4\ell}(F)\backslash \Sp_{4\ell}(\BA)}\theta_{\Delta(\tau,m+\ell)}(i(g,h))(\eta^\chi_{m,\ell}(\xi_{m,\ell,v})\ast E(f_{\Delta(\tau,\ell),s};g)dg=
	$$
	$$
	=P(q_v^{-s},q_v^s)\int\limits_{\Sp_{4\ell}(F)\backslash \Sp_{4\ell}(\BA)}(\theta_{\Delta(\tau,m+\ell)}(i(g,h))E(f_{\Delta(\tau,\ell),s};g)dg,
	$$
	where\\
	$P(q_v^{-s},q_v^s)=$
	$$
	\hat{\mathcal{S}}(\eta^\chi_{m,\ell}(\xi_{m,\ell,v}))(t_vq_v^{-s-\frac{\ell-1}{2}},t_vq_v^{-s-\frac{\ell-3}{2}},...,t_vq_v^{-s-\frac{1-\ell}{2}},t_vq_v^{s+\frac{1-\ell}{2}},t_vq_v^{s+\frac{3-\ell}{2}},...,t_vq_v^{s+\frac{\ell-1}{2}}).
	$$
	Here, $t_v=\chi(p_v)$. Thus,
	\begin{equation}\label{5.22}
		P(q_v^{-s},q_v^s)=\alpha_{m,\ell,v}\prod_{i=1}^\ell(q_v^{-s-\frac{\ell-(2i-1)}{2}}-q_v^{-m+\ell-\frac{1}{2}})(q_v^{s+\frac{\ell-(2i-1)}{2}}-\chi(p_v)^2q_v^{-m+\ell-\frac{1}{2}})\cdot
	\end{equation}
	$$
	\cdot \prod_{i=1}^\ell(q_v^{-s-\frac{\ell-(2i-1)}{2}}-\chi(p_v)^2q_v^{-m+\ell-\frac{1}{2}})(q_v^{s+\frac{\ell-(2i-1)}{2}}-q_v^{-m+\ell-\frac{1}{2}}).
	$$
	Here, $\alpha_{m,\ell,v}$ is defined right after \eqref{3.20}.
	\begin{lem}\label{lem 5.5}
		When $m\geq 2\ell$,
		$$
		P(q_v^{-\frac{\ell}{2}},q_v^{\frac{\ell}{2}})=c_{m,\ell,v}\neq 0.
		$$
		When $\ell\leq m\leq 2\ell-1$, recall that we assume that $v$ is such that $\chi$ is not quadratic. Then $	P(q_v^{-s},q_v^s)$ has a simple zero at $s=\frac{\ell}{2}$.
	\end{lem}
	\begin{proof}
		Substituting $s=\frac{\ell}{2}$ in \eqref{5.22}, we need to solve one of the equations
		$$
		\pm (\frac{\ell}{2}+\frac{\ell-(2i-1)}{2})=-m+\ell-\frac{1}{2},
		$$
		or one of the equations
		$$
		\pm (\frac{\ell}{2}+\frac{\ell-(2i-1)}{2})=-m+\ell-\frac{1}{2}+2\alpha,
		$$
		for $1\leq i\leq \ell$, where $|\chi(p_v)|=q_v^\alpha$. When we consider the plus sign, we get either $m=i-1<\ell$, which is impossible, or $m=i-1+2\alpha$, which implies that $2\alpha$ is an integer. Since $|2\alpha|<1$, we get that $\alpha=0$, and then $m=i-1<\ell$, which is impossible. Consider the equations with the minus sign. Then we get $m=2\ell-i$, which is possible, when $\ell\leq m\leq 2\ell-1$, exactly for one $1\leq i\leq \ell$. We also get $m=2\ell-i+2\alpha$, which, as before, forces $\alpha=0$, and hence $\chi(p_v)=q_v^{\frac{\pi k_v}{log(q_v)}i}$, for some integer $k_v$, and then $\chi(p_v)^2=1$, that is $\chi^2=1$, contrary to our assumption. This is the point where we need this assumption, to guarantee that when $m\leq 2\ell-1$, $P(q_v^{-s},q_v^s)$ has a simple zero at $s=\frac{\ell}{2}$. 
	\end{proof}

	\section{Proof of Theorem \ref{thm 2.3}}
	
	We prove, in more detail,
	\begin{thm}\label{thm 6.1}
		Assume that $\theta_{\Delta(\tau,m+\ell)}$ is $i(K_{4\ell}\times K_{4m})$-finite. Let $f_{\Delta(\tau,\ell),s}$ be a smooth, holomorphic section of $\rho_{\Delta(\tau,\ell),s}$ (see \eqref{0.9}). Then, for $Re(s)$ sufficiently large, $h\in \Sp_{4m}(\BA)$,
		\begin{equation}\label{6.1}
			\int\limits_{\Sp_{4\ell}(F)\backslash \Sp_{4\ell}(\BA)}(1\otimes\xi_{m,\ell,v})\ast\theta_{\Delta(\tau,m+\ell)}(i(g,h))E(f_{\Delta(\tau,\ell),s};g)dg=
		\end{equation}
		$$
		=\sum_{\gamma\in Q_{2\ell}(F)\backslash \Sp_{4m}(F)}F(f_{\Delta(\tau,\ell),s},(1\otimes\xi_{m,\ell,v})\ast\theta_{\Delta(\tau,m+\ell)};\gamma h),
		$$
		where
		\begin{equation}\label{6.2}
			F(f_{\Delta(\tau,\ell),s},\theta_{\Delta(\tau,m+\ell)};h)
			=\int\limits_{U_{2\ell}(\BA)\backslash \Sp_{4\ell}(\BA)}\theta^{U_{4\ell},\psi_{V_{(2\ell)^2}}}_{\Delta(\tau,m+\ell)}(i(g,h))f_{\Delta(\tau,\ell),s}(g)dg.
		\end{equation}
		(See right after \eqref{4.30} for the definition of the Fourier coefficient $\theta^{U_{4\ell},\psi_{V_{(2\ell)^2}}}_{\Delta(\tau,m+\ell)}$.)
		The integral \eqref{6.2} converges absolutely for $Re(s)$ sufficiently large and\\ $F(f_{\Delta(\tau,\ell),s},\theta_{\Delta(\tau,m+\ell)};\cdot)$ defines a $K_{4m}$-finite,
		holomorphic section of 
		$$
		\rho_{\Delta(\tau,\ell)|\det\cdot|^s\otimes \Theta_{\Delta(\tau,m-\ell)}}=\Ind_{Q_{2\ell}(\BA)}^{\Sp_{4m}(\BA)}\Delta(\tau,\ell)|\det\cdot|^s\otimes \Theta_{\Delta(\tau,m-\ell)}.	
		$$	
	\end{thm}
	\begin{proof}
		Since $g\mapsto (1\otimes\xi_{m,\ell,v})\ast\theta_{\Delta(\tau,m+\ell)}(i(g,h))$ is rapidly decreasing, for any given $h$, the integral \eqref{6.1} converges absolutely, for all $s$, which is not a pole of $E(f_{\Delta(\tau,\ell),s})$. Since we assume that $\theta_{\Delta(\tau,m+\ell)}$ is $K_{4(m+\ell)}$-finite, we may assume that $f_{\Delta(\tau,\ell),s}$ is $K_{4\ell}$-finite. (The integral \eqref{6.1} factors through such sections). Write, for $Re(s)$ sufficiently large,
		$$
		E(f_{\Delta(\tau,\ell),s};g)=\sum_{\gamma\in Q_{2\ell}(F)\backslash \Sp_{4\ell}(F)}f_{\Delta(\tau,\ell),s}(\gamma g),
		$$
		and substitute in the l.h.s. of \eqref{6.1}, which yields
		$$
		\int\limits_{Q_{2\ell}(F)U_{2\ell}(\BA)\backslash \Sp_{4\ell}(\BA)}(1\otimes\xi_{m,\ell,v})\ast\theta^{U_{2\ell}^{4\ell}\times 1}_{\Delta(\tau,m+\ell)}(i(g,h))f_{\Delta(\tau,\ell),s}(g)dg.
		$$
		
		By \eqref{4.31}, this is equal to
		\begin{equation}\label{6.3}
			\sum_{\gamma\in Q_{2\ell}(F)\backslash \Sp_{4m}(F)}\int\limits_{U_{2\ell}(\BA)\backslash \Sp_{4\ell}(\BA)}(1\otimes\xi_{m,\ell,v})\ast\theta^{U_{4\ell},\psi_{V_{(2\ell)^2}}}_{\Delta(\tau,m+\ell)}(g,\gamma h)f_{\Delta(\tau,\ell),s}(g)dg.
		\end{equation}
		This is \eqref{6.1}. In order to justify the passage to \eqref{6.3}, we will show that, for $Re(s)$ sufficiently large, 
		\begin{equation}\label{6.4}
			\int\limits_{U_{2\ell}(\BA)\backslash \Sp_{4\ell}(\BA)}|\theta^{U_{4\ell},\psi_{V_{(2\ell)^2}}}_{\Delta(\tau,m+\ell)}(g, h)f_{\Delta(\tau,\ell),s}(g)|dg<\infty,
		\end{equation}
		for all $K_{4\ell}$-finite $f_{\Delta(\tau,\ell),s}$ and $K_{4(m+\ell)}$-finite elements $\theta_{\Delta(\tau,m+\ell)}\in \Theta_{\Delta(\tau,m+\ell)}$. We will then show the assertion that $F(f_{\Delta(\tau,\ell),s},\theta_{\Delta(\tau,m+\ell)};h)$ is a section of\\ $\rho_{\Delta(\tau,\ell)|\det\cdot|^s\otimes \Theta_{\Delta(\tau,m-\ell)}}$, so that the summation over $\gamma$ in \eqref{6.3} is the usual summation defining an Eisenstein series, and hence it converges absolutely, since $Re(s)$ is sufficiently large. We may assume that $h=I_{4m}$. Using the Iwasawa decomposition in $\Sp_{4\ell}(\BA)$, it is enough to consider the integral 
		\begin{equation}\label{6.5}
			\int\limits_{\GL_{2\ell}(\BA)}|\theta^{U_{4\ell},\psi_{V_{(2\ell)^2}}}_{\Delta(\tau,m+\ell)}(\hat{a}, I_{4m})\varphi_{\Delta(\tau,\ell)}(a)| |det(a)|^{s-\ell-\frac{1}{2}}da,
		\end{equation}
		where $\varphi_{\Delta(\tau,\ell)}$ is an automorphic form in the space of $\Delta(\tau,\ell)$. By Prop. \ref{prop 1.2}, we may consider in place of \eqref{6.5}, the integral
		\begin{equation}\label{6.6}
			\int\limits_{\GL_{2\ell}(\BA)}|\varphi^{\psi_{V_{(2\ell)^2}}}_{\Delta(\tau,2\ell)}(\begin{pmatrix}a\\&I_{2\ell}\end{pmatrix})\varphi_{\Delta(\tau,\ell)}(a)| |det(a)|^{s+m-\ell}da,
		\end{equation}
		where $\varphi_{\Delta(\tau,2\ell)}$ is an automorphic form in the space of $\Delta(\tau,2\ell)$, and, for $y\in \GL_{4\ell}(\BA)$,
		$$
		\varphi^{\psi_{V_{(2\ell)^2}}}_{\Delta(\tau,2\ell)}(y)=\int\limits_{M_{2\ell}(F)\backslash M_{2\ell}(\BA)}\varphi^{\psi_{V_{(2\ell)^2}}}_{\Delta(\tau,2\ell)}(\begin{pmatrix}I_{2\ell}&x\\&I_{2\ell}\end{pmatrix}y)\psi^{-1}(tr(x))dx.
		$$
		Let us rewrite \eqref{6.6} using the Cartan decomposition. We get
		\begin{equation}\label{6.7}
			\int\limits_{K_{\GL_{2\ell}}}\int\limits_{T^+_{2\ell}(\BA)}\int\limits_{K_{\GL_{2\ell}}}|\varphi^{\psi_{V_{(2\ell)^2}}}_{\Delta(\tau,2\ell)}(\begin{pmatrix}k_1tk_2\\&I_{2\ell}\end{pmatrix})\varphi_{\Delta(\tau,\ell)}(k_1tk_2)| |det(t)|^{s+m-\ell}\gamma(t)dk_1dtdk_2,
		\end{equation}
		where $T^+_{2\ell}(\BA)$ denotes the set of diagonal matrices $t=diag(t_1,...,t_{2\ell})\in T_{2\ell}(\BA)$, such that at any place $v$, $|t_{1,v}|_v\leq |t_{2,v}|_v\leq...\leq |t_{2\ell,v}|_v$, and, for $v$ archimedean, $t_{1,v},...,t_{2\ell,v}$, are all positive. We will explicate $\gamma(a)=\prod_v\gamma_v(a_v)$ soon. 
		By Lemma 1.6 in \cite{CFK18}, we have in \eqref{6.7}
		$$
		\varphi^{\psi_{V_{(2\ell)^2}}}_{\Delta(\tau,2\ell)}(\begin{pmatrix}k_1tk_2\\&I_{2\ell}\end{pmatrix})=\varphi^{\psi_{V_{(2\ell)^2}}}_{\Delta(\tau,2\ell)}(\begin{pmatrix}tk_2\\&k_1^{-1}\end{pmatrix}).
		$$
		Hence due to the $K_{\GL_{4\ell}}$-finiteness of $\varphi_{\Delta(\tau,2\ell)}$, it is enough to consider 
		\begin{equation}\label{6.8}
			\int\limits_{T^+_{2\ell}(\BA)}|\varphi^{\psi_{V_{(2\ell)^2}}}_{\Delta(\tau,2\ell)}(\begin{pmatrix}t\\&I_{2\ell}\end{pmatrix})\xi_{\Delta(\tau,\ell)}(t)| |det(t)|^{s+m-\ell}\gamma(t)dt,
		\end{equation}
		where $\xi_{\Delta(\tau,\ell)}$ is a matrix coefficient of $\Delta(\tau,\ell)$. If $\varphi^{\psi_{V_{(2\ell)^2}}}_{\Delta(\tau,2\ell)}(\begin{pmatrix}t\\&I_{2\ell}\end{pmatrix})\neq 0$, then, for each finite place $v$,  $|t_{2\ell,v}|_v\leq A_v$, and $A_v=1$, for almost all finite $v$. The argument is standard. We take a matrix $z\in M_{2\ell}(F_v)$, close to zero. Think of $z$ as an element of $M_{2\ell}(\BA)$ having the zero coordinate at all places other than $v$. We take $z$ sufficiently close to zero, so that $\varphi_{\Delta(\tau,2\ell)}$ is fixed by the right translation by $\begin{pmatrix}I_{2\ell}&z\\&I_{2\ell}\end{pmatrix}$. Then we get, for all such $z$,
		$$
		\varphi^{\psi_{V_{(2\ell,v)^2}}}_{\Delta(\tau,2\ell)}(\begin{pmatrix}t\\&I_{2\ell}\end{pmatrix})=\psi_v(tr(t_vz))\varphi^{\psi_{V_{(2\ell)^2}}}_{\Delta(\tau,2\ell)}(\begin{pmatrix}t\\&I_{2\ell}\end{pmatrix}).
		$$
		This necessarily bounds $|t_{2\ell}|_v$ as we want. The analogue at an archimedean place $v$ is that $\varphi^{\psi_{V_{(2\ell)^2}}}_{\Delta(\tau,2\ell)}(\begin{pmatrix}t\\&I_{2\ell}\end{pmatrix})$ is rapidly decreasing as $t_{2\ell,v}$ tends to infinity. We omit the details. Using the moderate growth of $\varphi_{\Delta(\tau,\ell)}$, $\xi_{\Delta(\tau,2\ell)}$, we may bound \eqref{6.8} by
		\begin{equation}\label{6.9}
			c\int\limits_{T^+_{2\ell}(\BA)}\phi(t)||t||^A |det(t)|^{s+m-\ell}\gamma(t)dt,
		\end{equation}
		where $c, A>0$ and $\phi\in \mathcal{S}(\BA^{2\ell})$ is a positive Schwartz function. Let $t\in T^+_{2\ell}(\BA)$. Then, for a finite place $v$, $\gamma_v(t_v)=\delta^{-1}_{B_{\GL_{2\ell}}}(t_v)$ (assuming that the measure of $\GL_{2\ell}(\mathcal{O}_v)$ is 1). For an archimedean place $v$, there is $c_1>0$, such that $\gamma_v(t_v)\leq c_1\delta^{-1}_{B_{\GL_{2\ell}}}(t_v)$. Thus, the integral \eqref{6.9} is majorized by
		\begin{equation}\label{6.10}
			c_2\int\limits_{T^+_{2\ell}(\BA)}\phi(t)||t||^A |det(t)|^{s+m-\ell}\delta^{-1}_{B_{\GL_{2\ell}}}(t)dt,  \ c_2>0.
		\end{equation}
		The integral \eqref{6.10} converges for $Re(s)$ sufficiently large. This proves \eqref{6.4}. 
		
		Assume that $Re(s)$ is sufficiently large, so that the integral \eqref{6.2} converges absolutely. We will show now that $F(f_{\Delta(\tau,\ell),s},\theta_{\Delta(\tau,m+\ell)};h)$ is a section of\\ $\rho_{\Delta(\tau,\ell)|\det\cdot|^s\otimes \Theta_{\Delta(\tau,m-\ell)}}$. It is immediate to check that 
		$F(f_{\Delta(\tau,\ell),s},\theta_{\Delta(\tau,m+\ell)};h)$ is left invariant to $U_{2\ell}^{4m}(\BA)$. Let $a\in \GL_{2\ell}(\BA)$, $b\in \Sp_{4(m-\ell)}(\BA)$. Then\\
		
		$ F(f_{\Delta(\tau,\ell),s},\theta_{\Delta(\tau,m+\ell)};\begin{pmatrix}a\\&b\\&&a^*\end{pmatrix}h)=$
		\begin{equation}\label{6.11}
			=\int\limits_{U_{2\ell}(\BA)\backslash \Sp_{4\ell}(\BA)}\theta^{U_{4\ell},\psi_{V_{(2\ell)^2}}}_{\Delta(\tau,m+\ell)}(\begin{pmatrix}I_{2\ell}\\&a\\&&b\\&&&a^*\\&&&&I_{2\ell}\end{pmatrix}i(g,h))f_{\Delta(\tau,\ell),s}(g)dg.
		\end{equation}
		As in \eqref{6.6}, using Prop. \ref{prop 1.2}, the first factor of the integrand in \eqref{6.11} has the form
		\begin{equation}\label{6.12}
			q(\begin{pmatrix}I_{2\ell}\\&a\\&&b\\&&&a^*\\&&&&I_{2\ell}\end{pmatrix}i(g,h)),
		\end{equation}
		where $q$ is a function in the space of 
		\begin{equation}\label{6.13}
			\Ind_{Q_{4\ell}(\BA)}^{\Sp_{4(m+\ell)}(\BA)}\Delta_\psi(\tau,2\ell)|\det\cdot|^{-m}\otimes \Theta_{\Delta(\tau,m-\ell)}.
		\end{equation}
		Here, $\Delta_\psi(\tau,2\ell)$ denotes the representation of $\GL_{4\ell}(\BA)$, by right translations in the space of functions $\varphi^{\psi_{V_{(2\ell)^2}}}_{\Delta(\tau,2\ell)}$, where $\varphi_{\Delta(\tau,2\ell)}$ lies in the space of $\Delta(\tau,2\ell)$. By Lemma 1.6 in\cite{CFK18}, \eqref{6.12} is equal to 
		$$
		|\det(a)|^{2m+1}q(i(\hat{a}^{-1}g,\begin{pmatrix}I_{2\ell}\\&b\\&&I_{2\ell}\end{pmatrix}h)).
		$$
		Substituting in \eqref{6.11}, and changing variable $g\mapsto \hat{a}g$, we get
		\begin{equation}\label{6.14}
			|\det(a)|^{2(m-\ell)}\int\limits
			_{U_{2\ell}(\BA)\backslash \Sp_{4\ell}(\BA)}\theta^{U_{4\ell},\psi_{V_{(2\ell)^2}}}_{\Delta(\tau,m+\ell)}(i(g,\begin{pmatrix}I_{2\ell}\\&b\\&&I_{2\ell}\end{pmatrix}h))f_{\Delta(\tau,\ell),s}(\hat{a}g)dg.
		\end{equation}
		We may take $h=I_{4m}$. Using the Iwasawa decomposition in $\Sp_{4\ell}(\BA)$, the $K_{4\ell}$-finiteness of the integrand, and \eqref{6.13}, we get that \eqref{6.14} is a finite sum \\
		$\delta_{Q_{2\ell}^{4m}}^{\frac{1}{2}}(\begin{pmatrix}a\\&I_{4(m-\ell)}\\&&a^*\end{pmatrix})|\det(a)|^s\sum_{j=1}^N\theta^{(j)}_{\Delta(\tau,m-\ell)}(b)\cdot$
		\begin{equation}\label{6.15}
			\cdot \int\limits_{\GL_{2\ell}(\BA)}\varphi^{(j)}_{\Delta(\tau,\ell)}(ac)(\varphi_{\Delta(\tau,2\ell)}^{(j)})^{\psi_{V_{(2\ell)^2}}}(\begin{pmatrix}c\\&I_{2\ell}\end{pmatrix})|\det(c)|^{s+m-\ell}dc.
		\end{equation}
		Here, $\varphi^{(j)}_{\Delta(\tau,\ell)}, \varphi_{\Delta(\tau,2\ell)}^{(j)}, \theta^{(j)}_{\Delta(\tau,m-\ell)}$ are $K$-finite elements of $\Delta(\tau,\ell), \Delta(\tau,2\ell), \Theta_{\Delta(\tau,m-\ell)}$, respectively. The integral in \eqref{6.15}, as a function of $a$, defines a function in the space of $\Delta(\tau,\ell)$. Formally, this is clear. To justify this, we use the Cartan decomposition, exactly as we did in the first part of the proof, and express the integral in \eqref{6.15} as a finite sum of  the form
		\begin{equation}\label{6.16}
			\sum_{j'} \int\limits_{T^+_{2\ell}(\BA)}\xi^{(j')}_{\Delta(\tau,\ell)}(t)(\tilde{\varphi}_{\Delta(\tau,2\ell)}^{(j')})^{\psi_{V_{(2\ell)^2}}}(\begin{pmatrix}t\\&I_{2\ell}\end{pmatrix})\gamma(t) |\det(t)|^{s+m-\ell}dt\cdot \tilde{\varphi}^{(j')}_{\Delta(\tau,\ell)}(a),
		\end{equation}
		where $\xi^{(j')}_{\Delta(\tau,\ell)}$ are matrix coefficients of $\Delta(\tau,\ell)$, and $\tilde{\varphi}_{\Delta(\tau,2\ell)}^{(j')}, \tilde{\varphi}_{\Delta(\tau,\ell)}^{(j')}$ are elements of $\Delta(\tau,2\ell), \Delta(\tau,\ell)$, respectively. Recall that we have seen in the first part of the proof that the integrals in \eqref{6.16} converge absolutely, for $Re(s)$ sufficiently large. This shows that $F(f_{\Delta(\tau,\ell),s},\theta_{\Delta(\tau,m+\ell)};h)$ is a section of $\rho_{\Delta(\tau,\ell)|\det\cdot|^s\otimes \Theta_{\Delta(\tau,m-\ell)}}$, and the proof of theorem is complete.

	\end{proof}
	
	We note that without using the Cartan decomposition, we can express the integral in \eqref{6.15}, as a finite sum, as follows, instead of \eqref{6.16},
	\begin{equation}\label{6.17}
		\sum_{j} \int\limits_{\GL_{2\ell}(\BA)}\eta^{(j)}_{\Delta(\tau,\ell)}(c)(\phi_{\Delta(\tau,2\ell)}^{(j)})^{\psi_{V_{(2\ell)^2}}}(\begin{pmatrix}c\\&I_{2\ell}\end{pmatrix})|\det(c)|^{s+m-\ell}dc\cdot \phi^{(j)}_{\Delta(\tau,\ell)}(a),
	\end{equation}
	where $\eta^{(j)}_{\Delta(\tau,\ell)}$ are matrix coefficients of $\Delta(\tau,\ell)$, and $\phi_{\Delta(\tau,2\ell)}^{(j)}, \phi_{\Delta(\tau,\ell)}^{(j)}$ are elements of $\Delta(\tau,2\ell), \Delta(\tau,\ell)$, respectively.
	\vspace{0.5cm}
	
	We now address the issue of analytic continuation of $F(f_{\Delta(\tau,\ell),s},\theta_{\Delta(\tau,m+\ell)};h)$.
	We view $\theta^{U_{4\ell},\psi_{V_{(2\ell)^2}}}_{\Delta(\tau,m+\ell)}$ as an element in the space of \eqref{6.13}. Let $v$ be a place of $F$. Denote the local factor at $v$ of $\Delta(\tau,2\ell)$ by $\Delta(\tau_v,2\ell)$, realized in a space $V_{\Delta(\tau_v,2\ell)}$. Similarly, denote the local factor at $v$ of $\Theta_{\Delta(\tau,m-\ell)}$ by $\Theta_{\Delta(\tau_v,m-\ell)}$, realized in a space $V_{\Theta_{\Delta(\tau_v,m-\ell)}}$. Note that, up to scalars, there is a unique (continuous when $v$ is archimedean) linear functional $C^{\psi_v}$ on $V_{\Delta(\tau_v,2\ell)}$, such that, for all $\xi\in V_{\Delta(\tau_v,2\ell)}$ 
	$$
	C^{\psi_v}(\Delta(\tau_v,2\ell)(\begin{pmatrix}I_{2\ell}&x\\&I_{2\ell}\end{pmatrix})\xi)=\psi_v(tr(x))C^{\psi_v}(\xi).
	$$
	This is a special case of a Theorem of Gourevitch and Kaplan \cite{GK22}. We may then realize $\Delta(\tau_v,2\ell)$ in the space $C(\Delta(\tau_v,2\ell),\psi_v)$ of functions on $GL_{2\ell}(F_v)$  $a\mapsto C^{\psi_v}(\Delta(\tau_v,2\ell)(a)u_v)$, $u_v\in V_{\Delta(\tau_v,2\ell)}$. Thus, we can write the analogous local sections of $F(f_{\Delta(\tau,\ell),s},\theta_{\Delta(\tau,m+\ell)};h)$ at each place and consider them for decomposable data. For this, we fix isomorphisms $p_{\tau,2\ell}, p_{\tau,\ell}, q_{\tau,m-\ell}$ of $\otimes_v'\Delta(\tau_v,2\ell), \otimes'_v\Delta(\tau_v,\ell)$, $\otimes'_v\Theta_{\Delta(\tau_v,m-\ell)}$ with $\Delta(\tau,2\ell), \Delta(\tau,\ell), \Theta_{\Delta(\tau,m-\ell)}$, respectively. Let $S_0$ be a finite set of places of $F$, containing the infinite places, outside which $\tau$ is unramified. For $v\notin S_0$, fix unramified vectors $v^0_{\tau_v,\ell}\in V_{\Delta(\tau_v,\ell)}$, $\eta^0_{\tau_v,m-\ell}\in V_{\Theta_{\Delta(\tau_v,m-\ell)}}$, and $W^0_{\Delta(\tau_v,2\ell)}\in C(\Delta(\tau_v,2\ell),\psi_v)$, such that $W^0_{\Delta(\tau_v,2\ell)}(I_{4\ell})=1$. For each place $v$, let $f_{\Delta(\tau_v,\ell),s}$ be a $K_{4\ell,v}$-finite, holomorphic section of $\rho_{\Delta(\tau_v,\ell),s}$, 
	and let $f_{C(\Delta(\tau_v,2\ell),\psi_v),\Theta_{\Delta(\tau_v,m-\ell)}}$ be a $K_{4\ell,v}\times K_{4m,v}$-finite function in the space of
	\begin{equation}\label{6.18}
		\Ind_{Q_{4\ell}(F_v)}^{\Sp_{4(m+\ell)}(F_v)}C(\Delta(\tau_v,2\ell),\psi_v)|\det\cdot|^{-m}\otimes \Theta_{\Delta(\tau_v,m-\ell)}. 
	\end{equation}
	Denote by $f'_{C(\Delta(\tau_v,2\ell),\psi_v),\Theta_{\Delta(\tau_v,m-\ell)}}$ the function $f_{C(\Delta(\tau_v,2\ell),\psi_v),\Theta_{\Delta(\tau_v,m-\ell)}}$ composed with the evaluation at $a=I_{4\ell}$ in the $\GL_{4\ell}(F_v)$ factor.
	Let $S$ be a finite set of places, containing $S_0$. Assume that for $v\notin S$, $f_{\Delta(\tau_v,\ell),s}=f^0_{\Delta(\tau_v,\ell),s}$ is the unramified section, such that $f^0_{\Delta(\tau_v,\ell),s}(I_{4\ell})=v^0_{\tau_v,\ell}$, and, similarly, $f_{C(\Delta(\tau_v,2\ell),\psi_v),\Theta_{\Delta(\tau_v,m-\ell)}}=f^0_{C(\Delta(\tau_v,2\ell),\psi_v),\Theta_{\Delta(\tau_v,m-\ell)}}$ is the unramified vector taking the value $W^0_{\Delta(\tau_v,2\ell)}\otimes \eta^0_{\tau_v,m-\ell}$ at $I_{4(m+\ell)}$. Define, for a place $v$ of $F$, and for $h\in \Sp_{4m}(F_v)$,\\
	$	F_v(f_{\Delta(\tau_v,\ell),s},f_{C(\Delta(\tau_v,2\ell),\psi_v),\Theta_{\Delta(\tau_v,m-\ell)}};h)=$
	\begin{equation}\label{6.19}
		=\int\limits_{U_{2\ell}(F_v)\backslash \Sp_{4\ell}(F_v)} f_{\Delta(\tau_v,\ell),s}(g)\otimes f'_{C(\Delta(\tau_v,2\ell),\psi_v),\Theta_{\Delta(\tau_v,m-\ell)}}(i(g,h))dg.
	\end{equation}
	For given $g,h$, the integrand is an element of $\Delta(\tau_v,\ell)\otimes \Theta_{\Delta(\tau_v,m-\ell)}$. Assume that in \eqref{6.2}, 
	$$
	f_{\Delta(\tau,\ell),s}=p_{\tau,\ell}\circ (\otimes_v f_{\Delta(\tau_v,\ell),s}),
	$$
	$$
	\theta_{\Delta(\tau,m+\ell)}^{U_{4\ell},\psi_{V_{(2\ell)^2}}}=(p_{\tau,2\ell}\otimes q_{\tau,m-\ell})\circ (\otimes_v f_{C(\Delta(\tau_v,2\ell),\psi_v),\Theta_{\Delta(\tau_v,m-\ell)}}).
	$$
	Then
	$$
	F(f_{\Delta(\tau,\ell),s},\theta_{\Delta(\tau,m+\ell)};h)=(p_{\tau,\ell}\circ q_{\tau,m-\ell})(\otimes_v F_v(f_{\Delta(\tau_v,\ell),s},f_{C(\Delta(\tau_v,2\ell),\psi_v),\Theta_{\Delta(\tau_v,m-\ell)}};h_v)).
	$$
	In complete analogy with the last part of the proof of Theorem \ref{thm 6.1}, we see that $F_v(f_{\Delta(\tau_v,\ell),s},f_{C(\Delta(\tau_v,2\ell),\psi_v),\Theta_{\Delta(\tau_v,m-\ell)}};\cdot)$, for $Re(s)$ sufficiently large, is a $K_{4m,v}$-finite section of $\Ind_{Q_{2\ell}(F_v)}^{\Sp_{4m}(F_v)}\Delta(\tau_v,\ell)|\det\cdot|^s\otimes \Theta_{\Delta(\tau_v,m-\ell)}$. In fact,  analogously with \eqref{6.15}, \eqref{6.16}, we get, for $a\in \GL_{2\ell}(F_v)$, $b\in \Sp_{4(m-\ell)}(F_v)$, with notation analogous to \eqref{6.15}, \eqref{6.16}, that it is a finite sum\\
	$F_v(f_{\Delta(\tau_v,\ell),s},f_{C(\Delta(\tau_v,2\ell),\psi_v),\Theta_{\Delta(\tau_v,m-\ell)}};\begin{pmatrix}a\\&b\\&&a^*\end{pmatrix})=\delta_{Q_{2\ell}^{4m}}^{\frac{1}{2}}(\hat{a})|\det(a)|^s$
	\begin{equation}\label{6.20}
		\sum_r (\int\limits_{T^+_{2\ell}(F_v)}\xi^{(r)}_{\Delta(\tau_v,\ell)}(t)W^{(r)}_{\Delta(\tau_v,2\ell)}(\begin{pmatrix}t\\&I_{2\ell}\end{pmatrix})\gamma_v(t)|\det(t)|^{s+m-\ell}dt)\cdot
	\end{equation}
	$$
	\Delta(\tau_v,\ell)(a)\varphi^{(r)}_{\Delta(\tau_v,\ell)}\otimes \Theta_{\Delta(\tau_v,m-\ell)}(b)\theta^{(r)}_{\Delta(\tau_v,m-\ell)}.
	$$
	As in \eqref{6.17}, we can write the following expression instead of \eqref{6.20},\\
	$F_v(f_{\Delta(\tau_v,\ell),s},f_{C(\Delta(\tau_v,2\ell),\psi_v),\Theta_{\Delta(\tau_v,m-\ell)}};\begin{pmatrix}a\\&b\\&&a^*\end{pmatrix})=\delta_{Q_{2\ell}^{4m}}^{\frac{1}{2}}(\hat{a})|\det(a)|^s$
	\begin{equation}\label{6.21}
		\sum_{r'} (\int\limits_{\GL_{2\ell}(F_v)}\xi^{(r')}_{\Delta(\tau_v,\ell)}(g)W^{(r')}_{\Delta(\tau_v,2\ell)}(\begin{pmatrix}g\\&I_{2\ell}\end{pmatrix})|\det(g)|^{s+m-\ell}dg)\cdot
	\end{equation}
	$$
	\Delta(\tau_v,\ell)(a)\varphi^{(r')}_{\Delta(\tau_v,\ell)}\otimes \Theta_{\Delta(\tau_v,m-\ell)}(b)\theta^{(r')}_{\Delta(\tau_v,m-\ell)}.
	$$
	Denote 
	$$
	\mathcal{L}(\xi_{\Delta(\tau_v,\ell)},W_{\Delta(\tau_v,2\ell)},s)=\int\limits_{\GL_{2\ell}(F_v)}\xi_{\Delta(\tau_v,\ell)}(g)W_{\Delta(\tau_v,2\ell)}(\begin{pmatrix}g\\&I_{2\ell}\end{pmatrix})|\det(g)|^{s+m-\ell}dg. 
	$$
	We note that, for $v\notin S_0$, when data are unramified, normalized, and with the measure of $K_{4\ell,v}$ taken to be $1$, \eqref{6.21}, evaluated at the identity, becomes
	\begin{equation}\label{6.22}
		F_v(f^0_{\Delta(\tau_v,\ell),s},f^0_{C(\Delta(\tau_v,2\ell),\psi_v),\Theta_{\Delta(\tau_v,m-\ell)}};I_{4m})=	\mathcal{L}(\xi^0_{\Delta(\tau_v,\ell)},W^0_{\Delta(\tau_v,2\ell)},s)\cdot (v^0_{\tau_v,\ell}\otimes \eta^0_{\tau_v,m-\ell}).
	\end{equation}
	Here, $\xi^{0}_{\Delta(\tau_v,\ell)}$ is the spherical matrix coefficient of $\Delta(\tau_v,\ell)$. Assume, also, that, for $v\notin S$, $\psi_v$ is normalized. As in the last proof, the integrals in \eqref{6.20}, \eqref{6.21}, converge absolutely for $Re(s)$ sufficiently large. 
	\begin{thm}\label{thm 6.2}
		In the notation of \eqref{6.22}, for $Re(s)$ sufficiently large,
		\begin{equation}\label{6.23}
			\mathcal{L}(\xi^0_{\Delta(\tau_v,\ell)},W^0_{\Delta(\tau_v,2\ell)},s)=L(\Delta(\tau_v,\ell)\times \tau_v,s+m-\ell+\frac{1}{2}).
		\end{equation}
	\end{thm}
	\begin{proof}
		The proof is similar to \cite{G18}, Sec. 2.2. Let $\tau_v=\Ind_{B_{\GL_2}(F_v)}^{\GL_2(F_v)}(\chi_v\times \chi_v^{-1})$, where $\chi_v$ is an unramified character of $F_v^*$.
		Then $\Delta(\tau_v,\ell)$ is the unramified irreducible subrepresentation of $\Ind_{B_{\GL_{2\ell}(F_v)}}^{\GL_{2\ell}(F_v)}\chi_{v,\ell}$, where
		$$
		\chi_{v,\ell}=\chi_v|\cdot|^{\frac{1-\ell}{2}}\otimes \chi^{-1}_v|\cdot|^{\frac{1-\ell}{2}}\otimes \chi_v|\cdot|^{\frac{3-\ell}{2}}\otimes \chi^{-1}_v|\cdot|^{\frac{3-\ell}{2}}\otimes\cdots\otimes \chi_v|\cdot|^{\frac{\ell-1}{2}}\otimes \chi^{-1}_v|\cdot|^{\frac{\ell-1}{2}}.
		$$
		Let $f^0_{\chi_{v,\ell}}$ be the unramified function in the space of $\Ind_{B_{\GL_{2\ell}(F_v)}}^{\GL_{2\ell}(F_v)}\chi_{v,\ell}$, such that $f^0_{\chi_{v,\ell}}(I_{2\ell})=1$. Then
		$$
		\xi^{0}_{\Delta(\tau_v,\ell)}(g)=\int\limits_{\GL_{2\ell}(\mathcal{O}_v)}f^0_{\chi_{v,\ell}}(kg)dk.
		$$	
		Recall that the function $g\mapsto W^0_{\Delta(\tau_v,2\ell)}(\begin{pmatrix}g\\&I_{2\ell}\end{pmatrix})$ is bi-$\GL_\ell(\mathcal{O}_v)$ invariant. We get that, for $Re(s)$ large, \\
		\\
		\vspace{0.3cm}
		$\mathcal{L}(\xi^{0}_{\Delta(\tau_v,\ell)},W^0_{\Delta(\tau_v,2\ell)},s)=$
		$$
		=\int\limits_{\GL_\ell(\mathcal{O}_v)}\int\limits_{GL_{2\ell}(F_v)}f^0_{\chi_{v,\ell}}(g)W^0_{\Delta(\tau_v,2\ell)}(\begin{pmatrix}k^{-1}g\\&I_{2\ell}\end{pmatrix})|\det(g)|^{s+m-\ell}dkdg=	
		$$	
		\begin{equation}\label{6.24}
			=\int_{GL_{2\ell}(F_v)}f^0_{\chi_{v,\ell}}(g)W^0_{\Delta(\tau_v,2\ell)}(\begin{pmatrix}g\\&I_{2\ell}\end{pmatrix})|\det(g)|^{s+m-\ell}dg.
		\end{equation}
		By Prop. 3.5 in \cite{LM20}, we can express $W^0_{\Delta(\tau_v,2\ell)}$, in the Shalika model, in terms of the Zelevinsky model (terminology of \cite{LM20}), as follows. Let $f^0_{\tau_v,2\ell}$	be the unramified function in the space of $\Ind_{P_{2^{2\ell}(F_v)}}^{\GL_{4\ell}(F_v)}\tau_v|\det\cdot|^{\frac{1}{2}-\ell}\times \tau_v|\det\cdot|^{\frac{3}{2}-\ell}\times\cdots\times \tau_v|\det\cdot|^{\ell-\frac{1}{2}}$, and assume that $\tau_v$ is realized in its Whittaker model with respect to $\psi_v$, $W(\tau_v,\psi_v)$, so that the value $f^0_{\tau_v,2\ell}(I_{4\ell})$ is the function on $GL_2(F_v)^{\times(2\ell)}$
		$$
		(a_1,...,a_{2\ell})\mapsto W^0_{\tau_v,\psi_v}(a_1)W^0_{\tau_v,\psi_v}(a_2)\cdots W^0_{\tau_v,\psi_v}(a_{2\ell}), 
		$$
		where $W^0_{\tau_v,\psi_v}$ is the unramified element of $W(\tau_v,\psi_v)$, such that $W^0_{\tau_v,\psi_v}(I_2)=1$. 
		Similarly, for $h\in \GL_{4\ell}(F_v)$, denote the value of $f^0_{\tau_v,2\ell}(h)$ at $(I_2,...,I_2)$ ($2\ell$ times) by 
		$f^0_{\tau_v,2\ell}(h)(1)$.
		Let $\epsilon_0$ be the following Weyl element in $\GL_{4\ell}(F_v)$,
		$$
		\epsilon_0=\begin{pmatrix}e_1&0\\0&e_1\\e_2&0\\0&e_2\\\vdots\\e_{2\ell}&0\\0&e_{2\ell}\end{pmatrix},
		$$
		where $e_1,...,e_{2\ell}$ is the standard basis of row vectors in $F_v^{2\ell}$. Then
		\begin{equation}\label{6.25}
			W^0_{\Delta(\tau_v,2\ell)}(\begin{pmatrix}g\\&I_{2\ell}\end{pmatrix})=\int_{Y_{2\ell}(F_v)}f^0_{\tau_v,2\ell}(\epsilon_0 \begin{pmatrix}g&y\\&I_{2\ell}\end{pmatrix})(1)dy,
		\end{equation}
		where $Y_{2\ell}(F_v)$ is the subspace of lower triangular nilpotent matrices in $M_{2\ell}(F_v)$. The expression \eqref{6.25} is a special case of \cite{CFGK21}, Sec. 3.2.
		
		For a given $g$, the integral \eqref{6.25} stabilizes in  $Y_{2\ell}(\mathcal{O}_v)$. Note, also, as in Lemma 3.2 in \cite{LM20}, that the l.h.s. of \eqref{6.25} is supported in $\GL_{2\ell}(F_v)\cap M_{2\ell}(\mathcal{O}_v)$. Substituting in \eqref{6.24}, and using the Iwasawa decomposition, we get\\
		
		\vspace{0.1cm}
		
		$\mathcal{L}(\xi^{0}_{\Delta(\tau_v,\ell)},W^0_{\Delta(\tau_v,2\ell)},s)=$
		\begin{equation}\label{6.26}
			\int\prod_{i=1}^\ell(\chi_v(\frac{t_{2i-1}}{t_{2i}})|t_{2i-1}t_{2i}|^{3i-2-\frac{3\ell}{2}}|t_{2i}|)f^0_{\tau_v,2\ell}(\epsilon_0 \begin{pmatrix}vt&y\\&I_{2\ell}\end{pmatrix})(1)|\det(t)|^{s+m-\ell}dydvdt,
		\end{equation}
		where the $dt$ integration is over $t=diag(t_1,...,t_{2\ell})\in T_{2\ell}(F_v)$, such that $|t_j|\leq 1, j=1,...,2\ell$.	The $dv$ integration is over $v\in V_{1^{2\ell}}(F_v)$, such that $|
		v_{i,j}|\leq |t_j|^{-1}$, for all $1\leq i<j\leq 2\ell$. Finally, the $dy$ integration is over $Y_{2\ell}(\mathcal{O}_v)$. 
		Write $v=\begin{pmatrix}1&-x_1\cdot v'\\&v'\end{pmatrix}$, where $x_1=(x_{1,2},...,x_{1,2\ell})$ and $v'\in V_{1^{2\ell-1}}(F_v)$. Then we have
		\begin{equation}\label{6.27}
			f^0_{\tau_v,2\ell}(\epsilon_0 \begin{pmatrix}vt&y\\&I_{2\ell}\end{pmatrix})(1)=\psi_v(\sum_{j=2}^{2\ell}x_{1,j}y_{j,1})f^0_{\tau_v,2\ell}(\epsilon_0 \begin{pmatrix}I_{2\ell}&y\\&I_{2\ell}\end{pmatrix}\begin{pmatrix}t_1\\&v't'\\&&I_{2\ell}\end{pmatrix}(1),
		\end{equation}	
		where $t'=diag(t_2,...,t_{2\ell})$. We may change the order of integration and integrate, first in the $x_{1,j}$  variables, and then we must have $|y_{j,1}|\leq |t_j|$, for $2\leq j\leq 2\ell$. Next, write  $v'=\begin{pmatrix}1&-x_2\cdot v''\\&v''\end{pmatrix}$, where $x_2=(x_{2,3},...,x_{2,2\ell})$ and $v''\in V_{1^{2\ell-2}}(F_v)$. After a simple change of variable in $y$, the integration in the variables $x_{2,j}$ gives inside the integrand in \eqref{6.26}
		$$
		\int\limits_{|x_{2,j}|\leq |t_j|^{-1}, 3\leq j\leq 2\ell}\psi_v(\sum_{j=3}^{2\ell}x_{2,j}y_{j,2})f^0_{\tau_v,2\ell}(\epsilon_0 \begin{pmatrix}I_{2\ell}&y\\&I_{2\ell}\end{pmatrix}\begin{pmatrix}t_1\\&t_2\\&&v''t''\\&&&I_{2\ell}\end{pmatrix}(1)dx_2,
		$$
		where $t''=diag(t_3,...,t_{2\ell})$. We get that $|y_{2,j}|\leq |t_j|$, $3\leq j\leq 2\ell$. We continue in this way and get that the $dydv$ integration in \eqref{6.26} is equal to
		$$
		|t_2t_3^2\cdots t_{2\ell}^{2\ell-1}|^{-1}\int\limits_{y\in Y_{2\ell}(F_v), |y_{j,i}|\leq |t_j|, i<j}f^0_{\tau_v,2\ell}(\epsilon_0 \begin{pmatrix}t\\&I_{2\ell}\end{pmatrix}\begin{pmatrix}I_{2\ell}&t^{-1}y\\&I_{2\ell}\end{pmatrix})(1)dy.
		$$	
		In the last integral, the coordinates of $t^{-1}y$ are all in $\mathcal{O}_v$, and hence it is equal to
		$$
		f^0_{\tau_v,2\ell}(\epsilon_0 \begin{pmatrix}t\\&I_{2\ell}\end{pmatrix})(1)=f^0_{\tau_v,2\ell}(\diag(\begin{pmatrix}t_1\\&1\end{pmatrix},...,\begin{pmatrix}t_{2\ell}\\&1\end{pmatrix})(1)=
		$$	
		\begin{equation}\label{6.28}
			=|t_1|^{\ell-\frac{1}{2}}|t_2|^{\ell-\frac{3}{2}}\cdots|t_{2\ell}|^{\frac{1}{2}-\ell}\prod_{i=1}^{2\ell}W^0_{\tau_v,\psi_v}(\begin{pmatrix}t_i\\&1\end{pmatrix}).
		\end{equation}
		Substitute \eqref{6.28} in \eqref{6.26}, and we get\\
		\\ 
		$\mathcal{L}(\xi^{0}_{\Delta(\tau_v,\ell)},W^0_{\Delta(\tau_v,2\ell)},s)=$
		$$
		\prod_{i=1}^\ell\int\limits_{F_v^*}W^0_{\tau_v,\psi_v}(\begin{pmatrix}t_{2i-1}\\&1\end{pmatrix})\chi_v(t_{2i-1})|t_{2i-1}|^{s+m-\ell+\frac{2i-1-\ell}{2}}d^*t_{2i-1}\cdot
		$$
		$$
		\cdot \prod_{i=1}^\ell\int\limits_{F_v^*}W^0_{\tau_v,\psi_v}(\begin{pmatrix}t_{2i}\\&1\end{pmatrix})\chi_v^{-1}(t_{2i})|t_{2i}|^{s+m-\ell+\frac{2i-1-\ell}{2}}d^*t_{2i}.
		$$
		The factors in the last products are the Jacquet-Langlands local unramified integrals for $\tau_v\times \chi_v^{\pm 1}$, which give
		\begin{equation}\label{6.29}
			\prod_{i=1}^\ell L(\tau_v\times\chi_v,s+m-\ell+\frac{2i-\ell}{2})L(\tau_v\times\chi^{-1}_v,s+m-\ell+\frac{2i-\ell}{2})=
		\end{equation}
		$$
		=\prod_{i=1}^\ell L(\tau_v\times\tau_v,s+m-\ell+\frac{2i-\ell}{2})=L(\Delta(\tau_v,\ell)\times\tau_v,s+m-\ell+\frac{1}{2}).
		$$

	\end{proof}
	
	Similar reasoning as in the last proof can be used to prove
	\begin{thm}\label{thm 6.3}
		In the previous notation, for a given place $v$, assume that $\tau_v=\Ind_{B_{\GL_2}}^{\GL_2(F_v)}\chi_v\times\chi_v^{-1}$ ($\chi_v$ is not necessarily unramified, and $v$ is not necessarily finite). Then $\mathcal{L}(\xi_{\Delta(\tau_v,\ell)},W_{\Delta(\tau_v,2\ell)},s)$ is a meromorphic function. It is a finite sum of the form
		$$
		\sum \prod_{i=1}^\ell L(W_i,\chi_v,s+m-\ell+\frac{2i-1-\ell}{2})L(W'_i,\chi^{-1}_v,s+m-\ell+\frac{2i-1-\ell}{2}),
		$$
		where $W_i,W'_i\in W(\tau_v,\psi_v)$, and $L(W_i,\chi_v,z)$ is the analytic continuation of the Jacquet-Langlands local integral, which converges absolutely for $Re(z)$ large, 
		$$
		\int_{F_v^*}W_i(\begin{pmatrix}t\\&1\end{pmatrix})\chi_v(t)|t|^zd^*t.
		$$
	\end{thm}
	\begin{proof}
		Since the function  $g\mapsto W_{\Delta(\tau_v,2\ell)}(\begin{pmatrix}g\\&I_{2\ell}\end{pmatrix})$ is bi-$K_{\GL_\ell,v}$ finite, as in \eqref{6.24}, it is enough to consider, with similar notation, for $Re(s)$ large, integrals of the form 	
		\begin{equation}\label{6.30}
			\mathcal{L}'(f_{\chi_v,\ell},W_{\Delta(\tau_v,2\ell)},s)=\int_{GL_{2\ell}(F_v)}f_{\chi_{v,\ell}}(g)W_{\Delta(\tau_v,2\ell)}(\begin{pmatrix}g\\&I_{2\ell}\end{pmatrix})|\det(g)|^{s+m-\ell}dg.
		\end{equation}
		We need to take $f_{\chi_v,\ell}$ which lies in the image of the intertwining operator defining $\Delta(\tau_v,\ell)$, but the proof of analytic continuation of $\mathcal{L}(\xi_{\Delta(\tau_v,\ell)},W_{\Delta(\tau_v,2\ell)},s)$ will certainly follow if we take any $f_{\chi_v,\ell}$. 	
		Using the Iwasawa decomposition in \eqref{6.30}, the $K_{\GL_{2\ell},v}$-finiteness of $f_{\chi_v,\ell}$, the Dixmier-Malliavin lemma in the archimedean case, and its simple analog in the non-archimedean case, it	is enough to consider integrals of the form
		\begin{equation}\label{6.31}
			\int\prod_{i=1}^\ell(\chi_v(\frac{t_{2i-1}}{t_{2i}})|t_{2i-1}t_{2i}|^{3i-2-\frac{3\ell}{2}}|t_{2i}|)W_{\Delta(\tau_v,2\ell)}(\begin{pmatrix}vt\\&I_{2\ell}\end{pmatrix})\phi(vt)|\det(t)|^{s+m-\ell}dvdt,
		\end{equation}	
		where $\phi$ is a Schwartz function on the upper $2\ell\times 2\ell$ triangular matrices $\mathcal{B}_{2\ell}(F_v)$. As in \eqref{6.25}, with similar notation, we have an expression
		\begin{equation}\label{6.32}
			W_{\Delta(\tau_v,2\ell)}(I_{4\ell})=\int_{Y_{2\ell}(F_v)}f_{\tau_v,2\ell}(\epsilon_0 \begin{pmatrix}I_{2\ell}&y\\&I_{2\ell}\end{pmatrix})(1)dy.
		\end{equation}	
		See \cite{CFGK21}, Sec. 3.2. We interpret the integral \eqref{6.32} as a repeated integral, as follows
		\begin{equation}\label{6.33}
			\int_{F_v}\int_{F_v^2}\cdots\int_{F_v^{2\ell-2}}\int_{F_v^{2\ell-1}}f_{\tau_v,2\ell}(\epsilon_0 \begin{pmatrix}I_{2\ell}&y\\&I_{2\ell}\end{pmatrix})(1)dy^{2\ell-1}\cdots dy^2dy^1,
		\end{equation}
		where, for $y\in Y_{2\ell}(F_v)$, $y^j$ denotes the part of the $j$-th column of $y$, below the diagonal,  $j=1,2,...2\ell-1$. Each integration converges absolutely, that is
		$$	
		\int_{F_v^i}\big |\int_{F_v^{i+1}}\cdots\int_{F_v^{2\ell-1}}f_{\tau_v,2\ell}(\epsilon_0 \begin{pmatrix}I_{2\ell}&y\\&I_{2\ell}\end{pmatrix})(1)dy^{2\ell-1}\cdots dy^{i+1}\big |dy^i<\infty.
		$$	
		We may view \eqref{6.33} as the value at $z_1=\cdots=z_\ell=0$ of the analytic continuation of the following functional	on $\Ind_{P_{2^{2\ell}(F_v)}}^{\GL_{4\ell}(F_v)}\tau_v|\det\cdot|^{z_1+\frac{1}{2}-\ell}\times\cdots \tau_v|\det\cdot|^{z_\ell-\frac{1}{2}}\times \tau_v|\det\cdot|^{-z_\ell+\frac{1}{2}}\times\cdots \times \tau_v|\det\cdot|^{-z_1+\ell-\frac{1}{2}}$, given for a smooth holomorphic section $f_{\tau,2\ell,\underline{z}}$ of the last induced representation, by the following integral, which is absolutely convergent, for $Re(z_i-z_{i+1}), i=1,...,\ell-1$ and $Re(z_\ell)$ sufficiently large,
		$$
		W(f_{\tau,2\ell,\underline{z}})=\int_{Y_{2\ell}(F_v)}f_{\tau_v,2\ell,\underline{z}}(\epsilon_0 \begin{pmatrix}I_{2\ell}&y\\&I_{2\ell}\end{pmatrix})(1)dy.
		$$
		By the Dixmier-Malliavin lemma in the archimedean case, and the smoothness of the section in the finite case, we may take
		$$
		f_{\tau,2\ell,\underline{z}}(h)=\int_{V_{1^{2\ell}}(F_v)}\varphi(u)f'_{\tau,2\ell,\underline{z}}(h\begin{pmatrix}I_{2\ell}\\&u\end{pmatrix})du,
		$$
		where $\varphi\in \mathcal{S}(V_{1^{2\ell}}(F_v))$, and $f'_{\tau,2\ell,\underline{z}}$ is a smooth holomorphic section. Then, in the above domain of $\underline{z}$,
		$$
		W(f_{\tau,2\ell,\underline{z}})(I_{4\ell})=\int_{Y_{2\ell}(F_v)}f'_{\tau_v,2\ell,\underline{z}}(\epsilon_0 \begin{pmatrix}I_{2\ell}&y\\&I_{2\ell}\end{pmatrix})(1)\hat{\varphi}(y)dy,
		$$
		where $\hat{\varphi}$ is an appropriate Schwartz function on $Y_{2\ell}(F_v)$. The last integral converges absolutely everywhere and is holomorphic. Going back to \eqref{6.31}, consider, first, the integrals
		\begin{equation}\label{6.34}
			\int_{V_{1^{2\ell}}(F_v)}\int_{Y_{2\ell}(F_v)}f_{\tau_v,2\ell,\underline{z}}(\epsilon_0 \begin{pmatrix}vt&y\\&I_{2\ell}\end{pmatrix})(1)\phi(vt)dydv=
		\end{equation}
		$$
		\int_{V_{1^{2\ell}}(F_v)}\int_{\mathcal{B}_{2\ell}(F_v)\backslash M_{2\ell}(F_v)}f_{\tau_v,2\ell,\underline{z}}(\epsilon_0 \begin{pmatrix}t&y\\&I_{2\ell}\end{pmatrix})(1)\psi_v^{-1}(tr(vy))\phi(vt)dydv.
		$$
		Write $v=I_{2\ell}+u$, $u\in \mathcal{N}_{2\ell}(F_v)$- the upper $2\ell\times 2\ell$ nilpotent matrices. We may assume that $\phi(vt)=\phi(t+ut)=\phi_1(t)\phi_2(ut)$, where $\phi_1\in \mathcal{S}(F_v^{2\ell})$, $\phi_2\in \mathcal{S}(\mathcal{N}_{2\ell}(F_v))$. Change variable $u\mapsto ut^{-1}$. Then the last integral is equal to
		$$
		\phi_1(t)\eta(t)\int_{\mathcal{B}_{2\ell}(F_v)\backslash M_{2\ell}(F_v)}f_{\tau_v,2\ell,\underline{z}}(\tilde{t}\epsilon_0 \begin{pmatrix}I_{2\ell}&t^{-1}y\\&I_{2\ell}\end{pmatrix})(1)\psi_v^{-1}(tr(y))\tilde{\phi}_2(t^{-1}y)dy.
		$$
		Here, for $t=diag(t_1,...,t_{2\ell})$, $\tilde{t}=diag(t_1,1,t_2,1,...,t_{2\ell},1)$; $\eta(t)=|t_2t_3^2\cdots t_{2\ell}^{2\ell-1}|^{-1}$, and
		$$
		\tilde{\phi}_2(y)=\int_{\mathcal{N}_{2\ell}(F_v)}\phi_2(u)\psi_v^{-1}(tr(uy))du.
		$$
		Write $y=\bar{y}+\underline{y}$, where $\bar{y}\in \mathcal{B}_{2\ell}(F_v)$ and $\underline{y}\in Y_{2\ell}(F_v)$. Then $\tilde{\phi}_2(y)=\tilde{\phi}_2(\underline{y})$ depends on $\underline{y}$ only, and then this is the Fourier transform of $\phi_2$ at $\underline{y}$. Thus, the integral \eqref{6.34} is equal to 
		$$
		\phi_1(t)\eta(t)\int_{Y_{2\ell}(F_v)}f_{\tau_v,2\ell,\underline{z}}(\tilde{t}\epsilon_0 \begin{pmatrix}I_{2\ell}&t^{-1}y\\&I_{2\ell}\end{pmatrix})(1)\tilde{\phi}_2(t^{-1}y)dy=
		$$
		$$
		\phi_1(t)\int_{Y_{2\ell}(F_v)}f_{\tau_v,2\ell,\underline{z}}(\tilde{t}\epsilon_0 \begin{pmatrix}I_{2\ell}&y\\&I_{2\ell}\end{pmatrix})(1)\tilde{\phi}_2(y)dy:=\phi_1(t)f^{\tilde{\phi}_2}_{\tau_v,2\ell,\underline{z}}(\tilde{t})(1).
		$$	
		We may evaluate at $z_i=0$, $i=1,...,\ell$, and thus, \eqref{6.31} is a sum of products of integrals of the form 
		$$
		\int\prod_{i=1}^\ell(\chi_v(\frac{t_{2i-1}}{t_{2i}})|t_{2i-1}t_{2i}|^{3i-2-\frac{3\ell}{2}}|t_{2i}|)\phi_1(t)f^{\tilde{\phi}_2}_{\tau_v,2\ell,\underline{0}}(\tilde{t})(1)|\det(t)|^{s+m-\ell}dt.
		$$
		This integral is equal to a finite sum of integrals of the form
		\begin{equation}\label{6.35}
			\prod_{i=1}^\ell\int_{F_v^*}W^{2i-1}_{\tau_v,\psi_v}(\begin{pmatrix}t_{2i-1}\\&1\end{pmatrix})\varphi_{2i-1}(t_{2i-1})\chi_v(t_{2i-1})|t_{2i-1}|^{s+m-\ell+\frac{2i-1-\ell}{2}}d^*t_{2i-1}\cdot
		\end{equation}
		$$
		\cdot \prod_{i=1}^\ell\int_{F_v^*}W^{2i}_{\tau_v,\psi_v}(\begin{pmatrix}t_{2i}\\&1\end{pmatrix})\varphi_{2i}(t_{2i})\chi_v^{-1}(t_{2i})|t_{2i}|^{s+m-\ell+\frac{2i-1-\ell}{2}}d^*t_{2i},
		$$	
		where, for $1\leq j\leq 2\ell$, $W^j_{\tau_v\psi_v}$ are functions in the Whittaker model $W(\tau_v,\psi_v)$, and $\varphi_j\in \mathcal{S}(F_v)$. Each integral in \eqref{6.35} converges absolutely for $Re(s)$ large and defines a meromorphic function of $s$. Finally, note that for $W\in W(\tau_v,\psi_v)$ and $\varphi\in \mathcal{S}(F_v)$,
		$$
		W(\begin{pmatrix}t\\&1\end{pmatrix})\varphi(t)=(\int_{F_v}\hat{\varphi}(-x)\tau_v(\begin{pmatrix}1&x\\&1\end{pmatrix})Wdx)(\begin{pmatrix}t\\&1\end{pmatrix}).
		$$
		This completes the proof of Theorem \ref{thm 6.3}.
	\end{proof}
	
	The last case to consider is a finite place $v$ where $\tau_v$ is supercuspidal, and then we have
	\begin{thm}\label{thm 6.4}
		Assume that $\tau_v$ is supercuspidal. Then $\mathcal{L}(\xi_{\Delta(\tau_v,\ell)},W_{\Delta(\tau_v,2\ell)},s)$ is a rational function of $q_v^{-s}$. Moreover,
		$$
		L(\Delta(\tau_v,\ell)\times \tau_v,s+m-\ell+\frac{1}{2})\mathcal{L}(\xi_{\Delta(\tau_v,\ell)},W_{\Delta(\tau_v,2\ell)},s)\in \BC[q_v^{-s},q_v^s].
		$$
	\end{thm}
	
	\begin{proof}
		As in the previous two proofs, it is enough to prove the theorem for the following integrals, which we consider, first, for $Re(s)$ sufficiently large,
		\begin{equation}\label{6.36}	
			\int\limits_{\GL_2(F_v)^{\times \ell}}\int\limits_{V_{2^\ell}(F_v)}\xi_{\tau_v}^1(g_1)\cdots\xi^\ell_{\tau_v}(g_\ell)W_{\Delta(\tau_v,2\ell)}(\begin{pmatrix}vm(\bar{g})\\&I_{2\ell}\end{pmatrix})\alpha(\bar{g})|\det(m(\bar{g}))|^{s+m-\ell}dvd\bar{g}.
		\end{equation}	
		Here, $m(\bar{g})=diag(g_1,...,g_\ell)$, $\alpha(\bar{g})=|\det(g_1)|^{\frac{3(1-\ell)}{2}}\det(g_2)|^{\frac{3(3-\ell)}{2}}\cdots |\det(g_\ell)|^{\frac{3(\ell-1)}{2}}$. The functions $\xi^i_{\tau_v}$ are matrix coefficients of $\tau_v$. Using the Iwasawa decomposition in each copy of $\GL_2(F_v)$, \eqref{6.36} is a finite sum of integrals of the form	
		\begin{equation}\label{6.37}	
			\int\limits_{T_{2\ell}(F_v)}	\int\limits_{V_{1^\ell}(F_v)}\prod_{i=1}^\ell \xi_{\tau_v}^i(\begin{pmatrix}1&v_{2i-1,2i}\\&1\end{pmatrix}\begin{pmatrix}\frac{t_{2i-1}}{t_{2i}}\\&1\end{pmatrix})W_{\Delta(\tau_v,2\ell)}(\begin{pmatrix}vt\\&I_{2\ell}\end{pmatrix})
		\end{equation}
		$$
		\beta(t)|\det(t)|^{s+m-\ell}dvdt,
		$$	
		where $\beta(t)=\alpha(t)|\frac{t_2t_4\cdots t_{2\ell}}{t_1t_3\cdots t_{2\ell-1}}|$. Since $\xi_{\tau_v}^i$ are compactly supported, modulo the center, and since each factor in the last integrand is smooth, \eqref{6.37} is a finite linear combination with coefficients which are constant multiples of integral powers of $q_v^{-s}$ of integrals of the form 
		\begin{equation}\label{6.38}	
			\int\limits_{T_\ell(F_v)}\int\limits_{V_{2^\ell}(F_v)}W_{\Delta(\tau_v,2\ell)}(\begin{pmatrix}vt'\\&I_{2\ell}\end{pmatrix})\phi(vt')\zeta(t)|\det(t)|^{2(s+m-\ell)}dvdt,
		\end{equation}
		where, for $t=diag(t_1,...,t_\ell)\in T_\ell(F_v)$, $t'=diag(t_1 I_2,...,t_\ell I_2)$, and\\ $\zeta(t)=|t_1^{3(1-\ell)}t_2^{3(3-\ell)}\cdots t_\ell^{3(\ell-1)}|$. As in \eqref{6.31}, we replaced, as we may,\\ $W_{\Delta(\tau_v,2\ell)}(\begin{pmatrix}vt'\\&I_{2\ell}\end{pmatrix})$ by $W_{\Delta(\tau_v,2\ell)}(\begin{pmatrix}vt'\\&I_{2\ell}\end{pmatrix})\phi(vt')$, with $\phi\in \mathcal{S}(\mathcal{B'}_{2\ell}(F_v))$, where
		\begin{equation}\label{6.39}
			\mathcal{B'}_{2\ell}(F_v)=\left\{\begin{pmatrix}a_1 I_2&x_{1,2}&x_{1,3}&\cdots& x_{1,\ell}\\&a_2I_2&x_{2,3}&\cdots&x_{2,\ell}\\&&&\ddots\\&&&&a_\ell I_2\end{pmatrix}\ |\ a_k\in F_v, x_{i,j}\in M_2(F_v)\right\}.
		\end{equation}
		We may assume that $\phi=\phi_1\otimes \phi_2$, where, in the notation of \eqref{6.39}, $\phi_1$ is a Schwartz function of $(a_1,...,a_\ell)\in F_v^\ell$, and $\phi_2$ is a Schwarts function of the nilpotent radical of $\mathcal{B'}_{2\ell}(F_v)$. The same arguments leading from \eqref{6.33} to \eqref{6.35}, show, with the same notation, that \eqref{6.38} has the form
		\begin{equation}\label{6.40}	
			\int\limits_{T_\ell(F_v)}\int\limits_{Y_{2\ell}(F_v)}f_{\tau_v,2\ell,\underline{0}}(\epsilon_0\begin{pmatrix}t'&t'y\\&I_{2\ell}\end{pmatrix})(1)\hat{\phi}_2(\bar{y})\phi_1(t)\zeta(t)|\det(t)|^{2(s+m-\ell)+1}dydt.
		\end{equation}		
		Here, $\bar{y}$ is the matrix obtained from $y$ by setting $y_{i+1,i}=0$, $i=1,...,2\ell-1$, and then
		$$
		\hat{\phi}_2(\bar{y})=\int\phi_2((x_{i,j})_{1\leq i<j\leq 2\ell})\psi_v^{-1}(\sum_{i=1}^{2\ell-1}\sum_{j=i+1}^{2\ell}tr(x_{i,j}\bar{y}_{j,i}))d(x_{i,j})_{1\leq i<j\leq 2\ell}
		$$
		We conclude that \eqref{6.40} is a finite sum of integrals of the form
		\begin{equation}\label{6.41}	
			\int\limits_{T_\ell(F_v)}\int\limits_{F_v^\ell}f_{\tau_v,2\ell,\underline{0}}(\tilde{t}'\epsilon_0u(y_1,...,y_\ell))(1)\phi_1(t)\zeta(t)|\det(t)|^{2(s+m-\ell)+1}dydt,
		\end{equation}
		where 
		$$
		\tilde{t}'=diag (\begin{pmatrix}t_1\\&1\end{pmatrix},\begin{pmatrix}t_1\\&1\end{pmatrix},...,\begin{pmatrix}t_\ell\\&1\end{pmatrix}, \begin{pmatrix}t_\ell\\&1\end{pmatrix}),
		$$
		$$
		u(y_1,...,y_\ell)=\begin{pmatrix}I_{2\ell}&u'(y_1,...,y_\ell)\\&I_{2\ell}\end{pmatrix},\ u'(y_1,...,y_\ell)=\begin{pmatrix}0\\y_1&0\\&&0\\&&y_2&0\\&&&&\ddots\\&&&&&0\\&&&&&y_\ell&0\end{pmatrix}.
		$$	
		Fix a positive integer $N_0$, sufficiently large, such that $f_{\tau_v,2\ell,\underline{0}}$ is fixed by right translations by $I_{4\ell}+M_{4\ell}(\mathcal{P}_v^{N_0})$. Write the inner $dy$-integral in \eqref{6.41} as the sum
		$$
		\sum_{1\leq i_1<\cdots i_r\leq \ell}I_{i_1,...,i_r}(f_{\tau_v,2\ell,\underline{0}},t),
		$$
		where
		$$
		I_{i_1,...,i_r}(f_{\tau_v,2\ell,\underline{0}},t)=\int\limits_{\begin{matrix}|y_{i_1}|,...,|y_{i_r}|\geq q_v^{N_0}\\|y_{j_1}|,...,|y_{j_{\ell-r}}|\leq q_v^{N_0-1}\end{matrix}}f_{\tau_v,2\ell,\underline{0}}(\tilde{t}'\epsilon_0u(y_1,...,y_\ell))(1)dy,
		$$
		with $j_1<\cdots j_{\ell-r}$ denoting the complement of $i_1,...,i_r$ inside $1,...,\ell$. Each one of these integrals can be expressed as a finite sum of integrals of the form
		$$
		\int\limits_{|y_{i_1}|,...,|y_{i_r}|\geq q_v^{N_0}}f_{\tau_v,2\ell,\underline{0}}(\tilde{t}'\epsilon_0u_{i_1,...,i_r}(y_{i_1},...,y_{i_r})))(1)d(y_{i_1},...y_{i_r}),
		$$
		where $u_{i_1,...,i_r}(y_{i_1},...,y_{i_r})=u(y_1,...,y_\ell)$, with $y_{j_1}=\cdots=y_{j_r}=0$, and we have, for $|y_{i_1}|,...,|y_{i_r}|\geq q_v^{N_0}$, 
		$$
		f_{\tau_v,2\ell,\underline{0}}(\tilde{t}'\epsilon_0u_{i_1,...,i_r}(y_{i_1},...,y_{i_r})))(1)=f_{\tau_v,2\ell,\underline{0}}(\tilde{t}'d_{i_1}(y_{i_1})\cdots d_{i_r}(y_{i_r})\tilde{\epsilon}_0)(1),
		$$
		where, for $1\leq i\leq \ell$, $	d_i(y_i)=diag(I_{4i-3},y_i^{-1},y_i,I_{4(\ell-i)+1})$, and $\tilde{\epsilon}_0=w_{i_1}\cdots w_{i_r}\epsilon_0$, with $w_i=diag(I_{4i-3},\begin{pmatrix}&1\\-1&\end{pmatrix},I_{4(\ell-i)+1})$. We used the identity
		$$
		\begin{pmatrix}1&y\\0&1\end{pmatrix}\begin{pmatrix}1&0\\-y^{-1}&1\end{pmatrix}=\begin{pmatrix}0&y\\y^{-1}&1\end{pmatrix},
		$$
		for $y\neq 0$, especially, for $|y|\geq q_v^{N_0}$.	
		Thus, the contribution of $I_{i_1,...,i_r}(f_{\tau_v,2\ell,\underline{0}},t)$ to \eqref{6.41} is a finite linear combination with coefficients in $\BC[q_v^{-s},q_v^s]$ of integrals of the form
		$$
		\int\limits_{(F_v^*)^r}\int\limits_{|y_1|,...,|y_r|\geq q_v^{N_0}}\prod_{j=1}^r W_j(\begin{pmatrix}t_iy_i\\&1\end{pmatrix})W'_j(\begin{pmatrix}t_iy_i\\&1\end{pmatrix})\varphi_i(t_i)|t_i|^{2(s+m-\ell)+2i-\ell}d^*(y,t),
		$$
		where $W_j,W'_j\in W(\tau_v,\psi_v)$, $\varphi_i\in \mathcal{S}(F_v)$. Changing variables $t_i\mapsto t_iy_i^{-1}$, and using the fact that, for a Whittaker function $W\in W(\tau_v,\psi_v)$, the function $t\mapsto W(\begin{pmatrix}t\\&1\end{pmatrix})$ is in $\mathcal{S}(F_v^*)$ (since $\tau_v$ is supercuspidal), the last integral is a finite linear combination, with coeeficients in $\BC[q_v^{-s},q_v^s]$ of integrals of the form
		$$
		\prod_{j=1}^r\int\limits_{|y_i|\leq q_v^{-N_0}} \varphi'_i(y_i)|y_i|^{2(s+m-\ell)+2i-\ell}d^*y_i,
		$$
		where $\varphi'_i\in \mathcal{S}(F_v)$. Each factor in the last product is a sum of a polynomial in $\BC[q_v^{-s},q_v^s]$ and a constant multiple of the local $L$-function at $v$, $L_v(2(s+m-\ell)+2i-\ell)=L(\tau_v\times\tau_v,s+m-\ell+\frac{2i-\ell}{2})$. Thus \eqref{6.41} is a linear combination 
		$$
		\sum_{1\leq i_1<\cdots i_r\leq \ell}p_{i_1,...,i_r}(q_v^{-s},q_v^s)\prod_{j=1}^rL(\tau_v\times\tau_v,s+m-\ell+\frac{2i_j-\ell}{2})=
		$$
		$$
		=Q(q_v^{-s},q_v^s)\prod_{j=1}^\ell L(\tau_v\times\tau_v,s+m-\ell+\frac{2j-\ell}{2})=Q(q_v^{-s},q_v^s)L(\Delta(\tau_v,\ell)\times\tau_v,s+m-\ell+\frac{1}{2}).
		$$
		Here, the coefficients $p_{i_1,...,i_r}(q_v^{-s},q_v^s), Q(q_v^{-s},q_v^s)$ are in $\BC[q_v^{-s},q_v^s]$.
		This proves Theorem \ref{thm 6.4}.	
		
	\end{proof}
	
	Altogether Theorems \ref{thm 6.2} - \ref{thm 6.4} give,
	\begin{cor}\label{cor 6.5} In the notation of \eqref{6.19} - \eqref{6.22}, consider a decomposable section 
		$F(f_{\Delta(\tau,\ell),s},\theta_{\Delta(\tau,m+\ell)};h)$ corresponding to $\otimes_v F_v(f_{\Delta(\tau_v,\ell),s},f_{C(\Delta(\tau_v,2\ell),\psi_v),\Theta_{\Delta(\tau_v,m-\ell)}};h_v)$.
		Then each local section,  $F_v(f_{\Delta(\tau_v,\ell),s},f_{C(\Delta(\tau_v,2\ell),\psi_v),\Theta_{\Delta(\tau_v,m-\ell)}};h_v)$ is meromorphic in the complex plane. Moreover, $$
		L(\Delta(\tau_v,\ell)\times\tau_v,s+m-\ell+\frac{1}{2})F_v(f_{\Delta(\tau_v,\ell),s},f_{C(\Delta(\tau_v,2\ell),\psi_v),\Theta_{\Delta(\tau_v,m-\ell)}};h_v)
		$$ 
		is holomorphic, and, for $v$ finite, it lies in $\BC[q_v^{-s},q_v^s]$. 	
		For $v\notin S$, \\
		\\	
		$F_v(f^0_{\Delta(\tau_v,\ell),s},f^0_{C(\Delta(\tau_v,2\ell),\psi_v),\Theta_{\Delta(\tau_v,m-\ell)}};I_{4m})=$
		$$
		=L(\Delta(\tau_v,\ell)\times\tau_v,s+m-\ell+\frac{1}{2})\cdot (v^0_{\tau_v,\ell}\otimes \eta^0_{\tau_v,m-\ell}).
		$$
	\end{cor}
	Summarizing this section, we proved Theorem \ref{thm 2.4}. Indeed, recall from \eqref{2.10}, that we defined\\
	\\
	$\mathcal{E}(\theta_{\Delta(\tau,m+\ell)},f_{\Delta(\tau,\ell),s};h)=$
	$$
	\frac{1}{P(q_v^{-s},q_v^s)}\int\limits_{\Sp_{4\ell}(F)\backslash \Sp_{4\ell}(\BA)}(1\otimes\xi_{m,\ell,v})\ast\theta_{\Delta(\tau,m+\ell)}(i(g,h))E(f_{\Delta(\tau,\ell),s};g)dg.
	$$
	where $P(q_v^{-s},q_v^s)\in \BC[q_v^{-s},q_v^s]$ is the polynomial \eqref{5.22}. The last integral is the integral \eqref{6.1}, and we proved in Theorem \ref{thm 6.1} and Corollary \ref{cor 6.5} that\\
	$\mathcal{E}(\theta_{\Delta(\tau,m+\ell)},f_{\Delta(\tau,\ell),s};h)$ is an Eisenstein series on $\Sp_{4m}(\BA)$, corresponding to\\
	$	\Ind_{Q_{2\ell}^{4m}(\BA)}^{\Sp_{4m}(\BA)}\Delta(\tau,\ell)|\det\cdot|^s\otimes \Theta_{\Delta(\tau,m-\ell)}$. Since the Eisenstein series $E(f_{\Delta(\tau,\ell),s};g)$ has at most a simple pole at $s=\frac{\ell}{2}$, which occurs for an appropriate choice of section, it follows from Lemma \ref{lem 5.5} that at $s=\frac{\ell}{2}$, that $\mathcal{E}(\theta_{\Delta(\tau,m+\ell)},f_{\Delta(\tau,\ell),s};h)$ has at most a simple pole, when $m\geq 2\ell$, and at most a double pole, when $\ell\leq m\leq 2\ell-1$. This completes the proof of Theorem \ref{thm 2.4}.

\end{document}